\documentclass[reqno,a4paper]{amsart}  

\usepackage{latexsym,xspace,enumerate,amsfonts,amsmath,amssymb,amscd}
\usepackage{xypic,color,syntonly,slashed,paralist,stmaryrd,graphicx}
\usepackage{bbm}
\usepackage{thmtools}
\usepackage{tikz-cd}
\usepackage{etex}
\usepackage{amsfonts}
\usepackage{amscd}
\usepackage{amsbsy}
\usepackage{mathtools}
\usepackage{amsxtra}
\usepackage{epsfig}
\usepackage{epic,eepic}
\usepackage{graphics}
\usepackage{graphicx}
\graphicspath{ {fig/} }
\usepackage[all]{xy}
\usepackage{psfrag}
\usepackage{amsmath}
\usepackage{amsthm}
\usepackage{setspace}
\usepackage{url}
\usepackage{here}
\usepackage{xparse}
\usepackage{xspace}
\usepackage{faktor}
\usepackage{hyperref}
\usepackage[a]{esvect} 

\usepackage{tikz}
\usetikzlibrary{calc}
\usetikzlibrary{arrows.meta}
\usetikzlibrary{decorations.pathmorphing}
\usetikzlibrary{decorations.pathreplacing,decorations.markings}
\def\centerarc[#1](#2)(#3:#4:#5)
              { \draw[#1] ($(#2)+({#5*cos(#3)},{#5*sin(#3)})$) arc (#3:#4:#5); }

\definecolor{pgreen}{RGB}{0, 100, 0}
\definecolor{pbrown}{RGB}{125,80,23}
\definecolor{dgreen}{RGB}{0, 100, 0}
\definecolor{dred}{RGB}{255, 0, 0}
\definecolor{dpurple}{RGB}{178,0,232}
\definecolor{ddarkpurple}{RGB}{134,30,174}
\definecolor{dorange}{RGB}{233,88,0}
\definecolor{dyellow}{RGB}{233,190,0}
\definecolor{dpink}{RGB}{233,0,200}

\DeclareMathAlphabet\mathbfcal{OMS}{cmsy}{b}{n}
\usepackage{adjustbox}

\DeclareRobustCommand{\SkipTocEntry}[9]{}
\usepackage{suetterl} 

\newcommand*{\calA}{\mathcal A}

\newcommand*{\calF}{\mathcal F}

\newcommand*{\calH}{\mathcal H}

\newcommand*{\calK}{\mathcal K}
\newcommand*{\calL}{\mathcal L}

\newcommand*{\calP}{\mathcal P}

\newcommand{\calQ}{{\mathcal Q}}

\newcommand*{\calX}{\mathcal X}

\newcommand{\calZ}{{\mathcal Z}}

\newcommand{\CC}{{\mathbb{C}}}

\newcommand{\HH}{{\mathbb{H}}}

\newcommand{\RR}{{\mathbb{R}}}
\newcommand{\TT}{{\mathbb{T}}}

\newcommand{\ZZ}{{\mathbb{Z}}}

\newcommand{\NN}{{\mathbb{N}}}

\newcommand{\XX}{{\mathbb{X}}}

\newcommand{\e}{\epsilon}

\newcommand{\GL}{\mr{GL}}
\newcommand{\SL}{\mr{SL}}

\newcommand{\IET}{\mr{IET}}
\newcommand{\SC}{\mr{SC}}

\newcommand{\CBR}{\mr{CBR}}

\newcommand{\mr}{\mathrm}

\newcommand{\wt}{\widetilde}

\newcommand{\wind}{\mathop{\rm wind}\nolimits}

\newcommand{\Mor}{\mathop{\rm Mor}\nolimits}
\newcommand{\Hom}{\mathop{\rm Hom}\nolimits}

\renewcommand{\ker}{\mathop{\rm ker}\nolimits}

\renewcommand{\Re}{\mathop{\rm Re}\nolimits}
\renewcommand{\Im}{\mathop{\rm Im}\nolimits}

\renewcommand{\mod}{\mathop{\rm mod}\nolimits}

\renewcommand{\div}{\mr{div}\,}

\def\sign{\mathrm{sign}}

\newcommand{\hor}{\mr{hor}\,}
\newcommand{\ver}{\mr{ver}\,}
\newcommand{\note}[1]{\marginpar{\raggedright\if@twoside\ifodd\c@page\raggedleft\fi\fi\sf\scriptsize \red{RMK: #1}}}

\newcommand\red[1]{\textcolor{red}{#1}}

\newcommand{\be}{\begin{equation}}
\newcommand{\ben}{\begin{equation}\nonumber}
\newcommand{\ee}{\end{equation}}
\newcommand{\bp}{\begin{para}}
\newcommand{\ep}{\end{para}}
\newcommand{\bps}{\begin{paras}}
\newcommand{\eps}{\end{paras}}

\newcommand{\frakr}{\mathfrak{r}}
\newcommand{\frakg}{\mathfrak{g}}

\newcommand{\bfF}{\mathbf{F}}

\newcommand{\bfy}{{\mathbf{y}}}

\newcommand{\bfR}{{\mathbf{R}}}

\newcommand{\coeff}{a}
\newcommand{\coeffb}{b}

\newtheorem{proposition}{\textbf{Proposition}} [section]
\newtheorem{lemma}[proposition]{\textbf{Lemma}} 
\newtheorem{corollary}[proposition]{\textbf{Corollary}}
\newtheorem{remark}[proposition]{\textbf{Remark}}
\newtheorem{theorem}[proposition]{\textbf{Theorem}}

\theoremstyle{definition}
\newtheorem{definition}[proposition]{\textbf{Definition}}

\newtheorem*{example*}{\textbf{Example}}
\newtheorem*{theorem*}{\textbf{Theorem}}

\def\be{\begin{equation}}   \def\ee{\end{equation}}     \def\bes{\begin{equation*}}    \def\ees{\end{equation*}}
\def\ba{\be\begin{aligned}} \def\ea{\end{aligned}\ee}   \def\bas{\bes\begin{aligned}}  \def\eas{\end{aligned}\ees}
\def\={\;=\;}  \def\+{\,+\,}          

\newcounter{para}[section]
\newenvironment{para}[2][]{\refstepcounter{para}\noindent\ignorespaces{\bf #1\thepara. #2.} \rmfamily}{\noindent\ignorespacesafterend\bigskip}
\newenvironment{paras}[1]{\noindent\ignorespaces{\bf #1.} \rmfamily}{\noindent\ignorespacesafterend\bigskip}

\numberwithin{proposition}{section}
\numberwithin{definition}{section}

\newcommand{\abs}{\mr{abs}}

\newcommand{\sh}{\mr{sh}}
\newcommand{\lng}{\mr{long}}
\newcommand{\am}{\mr{am}}
\newcommand{\br}{\mr{br}}

\def\wh#1{\widehat{#1}}
\def\nl#1{\vv{#1}}
\def\dnl#1{\vv{#1}}

\DeclareDocumentCommand{\barPr}{O{\wh{\pi}} O{ }}{{\overline{\mathcal{P}}}_{#2}({#1})}

\DeclareDocumentCommand{\tildecJ}{O{g,n} O{d}}{{\widetilde{\mathcal J}}_{#1}^{#2}}
\DeclareDocumentCommand{\barcJ}{O{g,n} O{d}}{{\overline{\mathcal J}}_{#1}^{#2}}

\DeclareDocumentCommand{\LMS}{ O{\mu} O{g,n}} {\Xi\overline{\mathcal{M}}_{#2}(#1)}

\DeclareDocumentCommand{\kLMS}{ O{\mu} O{g,n}} {\Xi^k\!\overline{\mathcal{M}}_{#2}(#1)}

\DeclareDocumentCommand{\QMS}{ O{\mu} O{g,n}} {\overline{\mathcal{Q}}_{#2}(#1)}

\newcommand{\wq}{\scalebox{1.4}[1.4]{\upshape\suetterlin q}}

\usetikzlibrary{intersections} 
\tikzset{
	every loop/.style={very thick},
	comp/.style={circle,black,draw,thin,inner sep=0pt,minimum size=5pt,font=\tiny},
	order bottom left/.style={pos=.05,left,font=\tiny},
	order top left/.style={pos=.9,left,font=\tiny},
	order bottom right/.style={pos=.05,right,font=\tiny},
	order top right/.style={pos=.9,right,font=\tiny},
	order middle right/.style={draw,thin,shape=rectangle,inner sep=2pt,outer sep=5pt,pos=.5,right,font=\tiny},
	order middle left/.style={draw,thin,shape=rectangle,inner sep=2pt,outer sep=5pt,pos=.5,left,font=\tiny},
	order node dis/.style={text width=.75cm},
	circled number/.style={circle, draw, inner sep=0pt, minimum size=12pt},
	below left with distance/.style={below left,text height=10pt},
	below right with distance/.style={below right,text height=10pt}
}

\DeclareRobustCommand{\SkipTocEntry}[9]{}

\usepackage[style=alphabetic,
            isbn=false,
            doi=false,
            url=false,
            backend=bibtex,
            maxnames =5, 
            maxalphanames = 6, firstinits = true
           ]{biblatex}
\defbibheading{bibliography}[\bibname]{%
  \section*{References}%
}
\AtEveryBibitem{\clearlist{language}}

\DeclareFieldFormat[article]{title}{{\it #1}}
\DeclareFieldFormat{journaltitle}{{\rm #1}}
\renewbibmacro{in:}{%
  \ifentrytype{article}{}{\printtext{\bibstring{in}\intitlepunct}}}
\DeclareFieldFormat[incollection]{title}{{\it #1}}
\DeclareFieldFormat{journaltitle}{{\rm #1}}

\bibliography{Bibwallcross}

\title[Wall-crossing]
{Wall-crossing formulas via spectral networks}

\date{\today}
\author{Johannes Horn, Martin M\"oller}

\begin{document}

\begin{abstract}
  We give a self-contained proof of the Kontsevich-Soibelman wall-crossing formula entirely in the scope of quadratic differentials without relying on input from DT theory. Our approach is based on path-lifting rules for spectral networks introduced by Gaiotto, Moore and Neitzke. We provide a framework to justify the convergence of the path liftings, including the cases with spiral domains. In particular, we define path lifting rules for spectral networks associated to holomorphic quadratic differentials. As an intermediate step in the proof of the wall-crossing formula, we show that upon extending the path lifting rules to $\calA_0$-laminations we generate the hat-homology lattice.
\end{abstract}
\maketitle

\tableofcontents

\section{Introduction}

Given a quadratic differential~$q$ on a Riemann surface~$X$ understanding and counting saddle connections, the finite-length geodesics between two zeros of~$q$, is one of the most important and most subtle problems in low-dimensional geometry and dynamics. There are at least three different interpretations of this counting problem. Historically first, for finite area quadratic differentials the number of saddle connections is infinite and a central topic in dynamics is to count the asymptotic growth rate as a function of the length of the saddle connection. This is measured by the Siegel-Veech constants and we refer to \cite{emz,MasurZorich, GoujSV, cmz, ADGZZ} for an (incomplete) overview over the subject.
\par
Second, Kontsevich and Soibelman provided in \cite{KS08} an approach to Donaldson-Thomas theory to count what physicists call BPS-invariants for stable objects in a 3-Calabi-Yau category provided with an extra datum, nowadays known as a Bridgeland stability condition \cite{BrStabonTriang}. The correspondence between certain moduli spaces of stability conditions and spaces of quadratic differentials (due to Gaiotto-Moore-Neitzke \cite{GMNwallcrossing} from the physicist's viewpoint, and Bridgeland-Smith \cite{BridgelandSmith} in mathematical terms) translates this counting into our framework. The Kontsevich--Soibelman wall-crossing formula does not provide a count for any specific quadratic differential. It rather describes how the BPS-invariants change as one varies the quadratic differential by expressing that a certain infinite product is covariantly constant.
\par
The goal of this paper is to provide a complete proof of the wall-crossing formula entirely in the scope of quadratic differentials, elaborating on the ideas of spectral networks and path lifting functions by Gaiotto, Moore and Neitzke \cite{GMN_spectral} without relying on the machinery of motivic DT-invariants. We first give a complete statement of the wall-crossing formula and then provide further details on known proofs and our strategy. 
\par
Third, we mention in passing that \cite{KontSoiFloerI} provides yet another 'count' of saddle connections, in the context of holomorphic Floer theory. We do not attempt to relate these ways of counting.
\par
\medskip
\paragraph{\textbf{The setup and the wall-crossing formula}} 
Let $(X,q)$ be a half-translation surface, a Riemann surface with a quadratic differential~$q$ with singularities of signature~$\mu$ that we suppose to have only simple zeros. We allow both that~$q$ has a higher order pole (infinite area) or at most simple poles (finite area). There is a canonical double cover $\pi: \Sigma \to X$ such that $\pi^*q$ is the square of a one-form. We fix such a square root~$\lambda$ and denote by $\Gamma \subset H_1(X,\ZZ)$ the anti-invariant part under the double cover involution~$\sigma$. The $\ZZ$-modules~$\Gamma$ form a local system over the moduli space $Q_{g}(\mu)$ of quadratic differentials of signature $\mu$, which is trivialized over the moduli space $Q^\Gamma_{g}(\mu)$ of \emph{(homology)-framed quadratic differentials}, as defined in Section~\ref{sec:qdiff}. We let $\TT = \Hom(\Gamma,\CC^*)$ be the \emph{twisted torus}, where multiplication is twisted by a sign as recalled in Section~\ref{sec:twistedtorus}.
\par
A direction $\theta \in S^1$ is called \emph{active}, if there is a saddle connection in the direction~$\theta$, and \emph{saddle-free} otherwise. Suppose that the active direction~$\theta$ is \emph{rank one}, i.e., that the saddle connection in direction~$\theta$ span a rank one submodule~of $\Gamma$, generated by~$\gamma_0$.  Then we define the \emph{BPS automorphism~$\calK_\theta$} by its action on $\alpha \in \Gamma$ as
\be \label{intro:K}
\calK_\theta([\alpha]) \=   \prod_{\gamma \in \NN \wh{\gamma}_0} (1-[\gamma])^{\Omega(\gamma) \langle \gamma, \alpha \rangle} \cdot [\alpha]
\ee
where $\wh{\gamma}_0$ is the (``hat-homology'') lift of~$\gamma_0$ to~$\Gamma$ and where~$\Omega(\gamma)$ is the BPS-invariant of the direction~$\theta$ defined in Section~\ref{sec:rankonedir}. This is an automorphism of the completed twisted torus $\CC_\Delta[[\TT]]$, a profinite completion of the character group of $\TT$. In fact, these automorphisms belong to a pro-unipotent Lie group~$\wh{G}_\Delta$, see  Section~\ref{sec:twistedtorus}.
\par
We say that~$(X,q)$ is \emph{rank one}, if every direction on this half-translation surface is either saddle-free or rank one. This condition is generic in the space of quadratic differentials, since the surfaces that are not rank one form a countable union of hyperplanes, indexed essentially by the rank two sub-lattices of~$\Gamma$.
\par
Under the rank one hypothesis for any acute cone~$\Delta$ whose boundary rays $\Delta^\pm$ are saddle-free and any rank one homology framed translation surface~$(X,q,f)$ the product BPS-automorphism
$$S_\Delta \,:=\, S_\Delta(X,q,f) \,:=\, \prod_{\theta \in \Delta} \calK_\theta $$
over all active rays in~$\Delta$, taken in clockwise order, converges to an element in $\wh{G}_\Delta$. 
\par
\begin{theorem} \label{thm:introWC} Suppose that $(X,q)$ has infinite area and that~$\mu$ does not have simple poles. Then the automorphism~$S_\Delta$ is covariantly constant while the boundary rays are non-active, i.e.\ for any path $P(t) = (X_t,q_t,f_t)$ with $t \in [0,1]$ connecting two homology-framed quadratic rank one differentials such that both boundary rays are non-active for all $q_t$, there is an equality
$$S_{\Delta}(X_0,q_0) \= S_{\Delta}(X_1,q_1) $$
of the associated BPS-automorphisms in $\wh{G}_\Delta$.
\end{theorem}
\par
Note that the theorem does not require the intermediate surfaces $(X_t,q_t,f_t)$ for $t \in (0,1)$ to be rank one. In particular the BPS-automorphisms are not necessarily defined there. We excluded simple poles given the length of paper. The approximation of homology generators by path lifts (see below and Section~\ref{sec:approx}) is the main step that requires additional considerations in the presence of simple poles. Our lack of an analog of this result in the finite area case gives a conditional statement only.  
\par
\begin{proposition} \label{thm:introfiniteareaWC} Suppose that $(X,q)$ has finite area and suppose the there are sufficiently many based path lifting functions so that an analog of Proposition~\ref{prop:approxFL_based} holds. Then the automorphism~$S_\Delta$ is covariantly constant while the boundary rays are non-active.
\end{proposition}
\par
However, as we illustrate in Section~\ref{sec:FiniteArea}, the statement is almost void in the finite area case anyway as the condition of non-active boundary rays is too restrictive. In turn, our method of proof (see Section~\ref{sec:FiniteArea}) shows (provided that there are sufficiently many based path liftings) that the restriction of $S_\Delta$ to the module of paths of length at most~$L$ is covariantly constant if no boundary ray of along the $P(t) = (X_t,q_t,f_t)$ has a saddle connection of length $\leq L'$ for an explicitly given $L' = L'(L)$.
\par
\medskip
\paragraph{\textbf{Approaches to the proof}} The original formulation of the wall-crossing formula in \cite{KS08} was written in the scope of Donaldson-Thomas theory. Further details on motivic DT-theory were given by Davison-Meinhardt and \cite{DMDT}, see also \cite{ḾeinhardtIntro,DMCohoDT}. Together these references contain a complete proof along the lines of motivic DT-theory. This approach is also available in the finite area case thanks to Haiden's theorem \cite{Haiden}, which identifies also those spaces of quadratic differentials with spaces of stability conditions on some CY3-categories. Our goal is to avoid this machinery entirely.
\par
Another approach to the wall-crossing formula is to use Fock-Goncharov coordinates. These coordinates on a certain moduli space of local systems satisfy similar transformation laws as~\eqref{intro:K}. This is the viewpoint taken by Allegretti in \cite{Allegretti}, see also \cite{AlleBridge}. However, his proof gives the wall-crossing formula only over the '$B_2$-locus' with at most one saddle connection in a given direction. For other more complicated configurations (see the list in Lemma~\ref{le:r1class} for the terminology of configurations) he relies on \cite{DMDT}.
\par
Gaiotto-Moore-Neitzke also provide in \cite{GMNwallcrossing} a strategy of proof of the wall-crossing formula using  Fock-Goncharov coordinates. Their argument includes the case of a ring domain. Attempting to pursue this approach in the presence of spiral domains, i.e.\ toral ends,  highlights a major difficulty: It amounts to keep track of the homology classes of saddle connections while approaching a direction in which they accumulate. While straightforward on a ring domain, it would require bookkeeping of two nested fractional part functions for a toral end, which seems not feasible.
\par
\medskip
\paragraph{\textbf{Spectral networks and path-lifting}} Our approach elaborates on the path lifting procedures with detours along spectral networks from \cite{GMN_spectral}. The detour paths are added to the trivial lifts of a given path~$\wp$ to obtain a lifting function that is well-defined under homotoping~$\wp$ over the zeros of~$q$. Our contribution is threefold:
\begin{itemize}
\item[i)] Handling various convergences issues, see \cite[Question~7]{GMN_spectral}.
\item[ii)] Confirming the expectation form \cite[Section~6.5]{GMN_spectral} that path lifting functions are 'rich enough'.
\item[iii)] Relating the limiting path lifting rules by the corresponding wall-crossing automorphism, including the case of toral ends.
\end{itemize}
\par
To address i) we work in a profinite setting, as did \cite{KS08}, truncating at paths with longer and longer central charge~$Z$, i.e., with longer $\lambda$-period. While we are finally interested in $\CC_\Delta[[\TT_-]]$ and thus in paths with central charge in~$\Delta$, the original path we are lifting might have arbitrary central charge. We are thus forced to enlarge $\CC_\Delta[[\TT_-]]$ and work with a \emph{completed (groupoid) module of signed paths} denoted by $\wh{\Pi}_\Delta$ in the version up to homotopy and by $\wh{H}_\Delta$ in the version of homology classes. The details are given in Section~\ref{sec:prelim}.
\par
We alert the reader that on these completed modules of signed paths the path composition does not extend to a well-defined multiplication, and the automorphism~$\calK_\theta$ is not continuous. However, elementary cone geometry implies that both properties hold on certain 'tame' submodules, see Proposition~\ref{prop:rngstructure} and Lemma~\ref{le:Kgeneralcont}.
\par
In Section~\ref{sec:pathliftF} we define the path lifting functions with detour rules, closely following the idea in \cite[Question~7]{GMN_spectral}. Our addendum here is the verification of convergence. See Proposition~\ref{prop:defF} for the first of several similar instances, where this is reduced to a zippered-rectangle construction well-known in quadratic differential literature. While the application to the wall-crossing formula only requires the homological version of the path lifting function $F(\wp,\theta)$, we define homotopy version~$\bfF(\wp,\theta)$ in parallel, with applications to non-abelianization in mind.
\par
To address ii) we follow in Section~\ref{sec:A0laminations} a proposal outlined in \cite{GMNframed} and extend the path lifting functions $\bfF(\cdot,\theta)$ from paths to Fock-Goncharov's $\calA_0$-laminations~$\calL$ \cite{FG_dual_Teichmueller}. These produce paths that follow along part of the boundary (of the real blowup of the singularities), while the rules for taking detours and thus homotopy invariance stay the same. The main statement, Proposition~\ref{prop:approxFL}, shows that generators of~$\Gamma$ can be \emph{approximated} in the topology induced by the profinite limit by polynomials in lifts of $\calA_0$-laminations. Moreover, we have to control the support of these approximations to work on submodules on which $\calK_\theta$ acts continuously. In the case of quadratic differentials with double poles (Section~\ref{sec:withholes}) our approximation statements are weaker and more complicated to state, but still enough to conclude.
\par
As we approach (rank one) directions with saddle connections, the one-sided limits of path lifting rules exist in the profinite topology as we show in Proposition~\ref{prop:Fconverge}.  Addressing iii), now the whole content of the wall-crossing formula has now been boiled down to comparing these one-sided limits, i.e., to a finite verification of homology classes of detour paths in a few lemmas -- with several pictures -- in Section~\ref{sec:BPSWC}, culminating in Proposition~\ref{prop:Kformula1}.
\par
We summarize the main properties of the lift and its relation with the BPS-automorphism~$\calK$. The following theorem combines Theorem~\ref{theo:homotopy_invariance}, Proposition~\ref{prop:Fconverge}, and Theorem~\ref{thm:WCforF}.
\par
\begin{theorem} \label{thm:introF}
Let $(X,q)$ be a quadratic differential. There exists a path lifting rule that
associates to each path~$\wp$ in~$X$ disjoint from zeros and poles of~$q$ and
to each saddle-free direction~$\theta$ an element $\bfF(\wp,\theta) \in \wh{\Pi}_\Delta$
with the following properties.
\begin{itemize}
\item[i)] Composition rule: For two any composable paths $\wp$ and~$\wp'$
$$\bfF(\wp \circ \wp',\theta) \= \bfF(\wp',\theta) \cdot \bfF(\wp,\theta)$$ 
\item[ii)] Full homotopy invariance:  For two paths $\wp$ and~$\wp'$ homotopic
in the complement of the poles $\bfF(\wp,\theta) =  \bfF(\wp',\theta)$.
\end{itemize} 
\par
Moreover, for each active rank one direction~$\theta_0$ there exist path lifting
rules $F^\pm(\wp,\theta_0) \in \wh{H}_\Delta$ that also satisfy the composition rule and such that the image $F(\wp,\theta) \in \wh{H}_\Delta$ of $\bfF(\wp,\theta)$ satisfies:
\begin{itemize}
\item[iii)] One-sided limits: For any path~$\wp$
$$\lim_{\theta \to \theta_0^\pm} F(\wp,\theta) \= F^\pm(\wp,\theta_0)\,.$$
\item[iv)] Wall crossing: For the automorphism~$\calK$ given in~\eqref{intro:K} and any~$\wp$
$$F^+(\wp,\theta_0) \=\calK  F^-(q,\theta_0)\,.$$
\end{itemize} 
\end{theorem}
\par
The precise definition of rank one directions and of~$\calK$ along with the proof of the theorem will be completed in Section~\ref{sec:rankonedir}.
\par
\medskip
\paragraph{\textbf{Outlook}}
Another application of spectral networks intended by Gaiotto-Moore-Neitzke is the abelianization of $\SL(K,\CC)$-local systems \cite[§10]{GMN_spectral}. The path lifting functions $\bfF(\cdot,\theta)$ define a homotopy-invariant way of lifting paths along the branched cover $\pi:\Sigma \to X$. Hence, given a (twisted) rank 1 local system on $\Sigma$ there is an associated parallel transport operator of an $\SL(K,\CC)$-local system on $X$. For spectral networks that are finite and saddle-free this was formulated mathematically in \cite{HollandsNeitzke,IonitaMorrissey}. Our work allows to define abelianization of $\SL(2,\CC)$-local system along spectral networks with dense trajectories, in particular, for local systems on closed Riemann surfaces.  This will be the focus of a sequel work. The abelianization of $\SL(2,\CC)$-local systems on closed Riemann surfaces shall provide a missing link in Bridgeland's program of geometrization of DT-invariants \cite{Bridgeland_Joyce}, where, to this day, Joyce structures are only well-understood for holomorphic quadratic differentials \cite{Bridgeland_Joyce_quadratic}. 
\par
Spectral networks were defined in \cite{GMN_spectral} not only for quadratic differentials but tuples $(\phi_2,\ldots,\phi_K)$ with $\phi_k$ a $k$-differential, more precisely for points in the regular locus of the $\SL(K,\CC)$-Hitchin base. While our (convergence and wall-crossing) results certainly apply to configurations of trajectories in higher rank spectral networks similar to those in Section~\ref{sec:rankonedir}, new infinities appear: Higher rank spectral networks include trajectories that are born at junctions (in contrast to just the critical trajectories encountered for $K=2$). This possibly leads to infinitely many trajectories that typically accumulate at the punctures. Nevertheless, we believe that similar convergence arguments can be applied to give a full geometric proof of a wall-crossing formula for $K>2$.
%
\par
\subsection*{Acknowledgments} The authors are deeply grateful to Andy Neitzke, whose countless answers and suggestions were indispensable in our attempt to understand the wall-crossing formula. We thank the Simons Center for Geometry and Physics (SCGP) for its hospitality during the program ``Geometric, Algebraic, and Physical structures around the moduli of Meromorphic Quadratic Differentials'', which contributed significantly to the progress on this project. We also thank Dylan Allegretti, Tom Bridgeland, and Jon Chaika for useful suggestions. The authors acknowledge support by Deutsche Forschungsgemeinschaft (DFG, German Research Foundation) through the Collaborative Research Centre TRR 326 \textit{Geometry and Arithmetic of Uniformized Structures}, project number 444845124.

\section{Preliminaries} \label{sec:prelim}

This section is a reference for standard material for quadratic differentials, their trajectory structure and their moduli spaces, see e.g.\ \cite{Strebel} or
\cite{BridgelandSmith}.

\subsection{Quadratic differentials} \label{sec:qdiff}

Let $\mu=(m_1,\dots,m_k) \in \ZZ^k$, such that $\sum_{i=1}^k m_k=4g-4$. Then a (possibly) meromorphic quadratic differential $(X,q)$ is of type $\mu$, if its divisor is of the form $\div(q)=\sum_{i=1}^k m_ix_i$ for some $x_i \in X$. We denote the set of zeros of~$q$ by $\calZ=\calZ(q)=\{z_1,\dots,z_r\}$ and the set of poles $\calP=\calP(q)=\{p_1,\dots,p_s\}$, obviously with $r+s=k$. We will mainly work with the open surface without the poles $X^\star= X\setminus \calP$ or the surface without all singularities of the flat structure $X^\circ=X\setminus (\calP \cup \calZ)$.
\par
We denote the moduli space of quadratic differentials of type~$\mu$ with~$k$ (unlabeled) zeros and poles by $\calQ_{g,[k]}(\mu)$ or briefly by $\calQ_{g}(\mu)$, and write $\calQ_{g,k}(\mu)$ for the cover where all special points are labeled.
\par
In the main part of this work, we consider (possibly) meromorphic quadratic differentials $(X,q)$ with simple zeros, i.e., we assume $\mu = 1^{4g-4+|\sum p_i|}$ unless stated otherwise. Unless stated otherwise we allow both the case of \emph{finite area differentials}, which is equivalent to~$q$ having at most simple poles and \emph{infinite area differentials}, equivalent to having a \emph{higher order poles}, i.e.\ a pole of order~$\geq 2$. Here we recall some background on the canonical cover, strata of differentials, and the classification of trajectories of quadratic differentials. 
\par
\medskip
\paragraph{\textbf{Canonical cover}}
To each meromorphic quadratic differential we can associate a canonical cover $\pi: \Sigma \to X=\Sigma/\sigma$, which is the defined to be the unique Riemann surface, such that there exists a meromorphic abelian differential $\lambda$, such that $\lambda^2+\pi^*q=0$. We have $\lambda=-\sigma^*\lambda$ with respect to the involution~$\sigma$ on~$\Sigma$. Consequently, $\pi$ is branched over the zeros and poles of~$q$ of odd order. For a critical point $x \in \calZ \cup \calP$ of $q$ of order $2\ell+1$ the abelian differential $\lambda$ has a critical point of order $2\ell+2$ at $\pi^{-1}x$. For a critical point~$x$ of~$q$ of order $2\ell$ the abelian differential has a critical point of order~$\ell$ at both preimages $\{y_1,y_2\}=\pi^{-1}(x)$. We will denote by $\Sigma^\star=\pi^{-1}(X^\star)$ the canonical cover without the preimages of the poles and by $\Sigma^\circ=\pi^{-1}(X^\circ)$ the canonical cover without the preimages of poles and zeros. 
\par
\medskip
\paragraph{\textbf{Strata of Differentials} }
The moduli space $\calQ_{g}(\mu)$ of quadratic differentials of type~$\mu$ is an orbifold of dimension $2g-2+k $ if~$\mu$ is meromorphic (i.e.\ there exists some $m_i <0$) and of dimension $2g-1+k$ otherwise. This can be seen by showing that the map
\[ Z: H_1(\Sigma\setminus\calP(\lambda),\CC)^- \to \CC,\quad  \gamma \to \int_\gamma \lambda,
\]
often called \emph{central charge}, defines local coordinates, called \emph{period coordinates}. Here $()^-$ denotes the anti-invariant part with respect to the involution $\sigma$. Note that $H_1(\Sigma\setminus\calP(\lambda),\calZ(\lambda),\CC)^- = H_1(\Sigma\setminus\calP(\lambda),\CC)^-$ if all zeros have odd order. We call
\be \label{eq:chargelattice}
\Gamma := H_1(\Sigma\setminus\calP(\lambda),\ZZ)^-
\ee
the \emph{charge lattice}, or following \cite{MZ08} and \cite{BridgelandSmith}, also \emph{hat-homology group}.
\par
We also define two spaces covering $\calQ_{g}(\mu)$. First, let  $\calQ^\Gamma_{g}(\mu)$ be the moduli space of homology-framed quadratic differentials, i.e.\
triples $(X,q,f: \Gamma \to H_1(\Sigma\setminus\calP(\lambda),\ZZ)^-)$ of a quadratic differential together with an isomorphism~$f$ between an abstract charge latte~$\Gamma$ and the hat-homology group of~$(X,q)$. Second, for comparison we also define $\calQ^{\pi_1}_{g}(\mu)$ to be the moduli space of Teichmüller-framed differentials, triples $(X,q,f: \Sigma_{0} \to \Sigma)$ where now~$f$ is a diffeomorphism (up to isotopy rel marked points and boundary) of a reference smooth surface~$\Sigma_{0}$ to the canonical double cover.
\medskip
\paragraph{\textbf{Half-translation structure}}
Let $(X,q)$ be a quadratic differential. Locally on~$X$ we can choose holomorphic coordinates $z$, such that $q=z^k d z^2$. An atlas of $X\setminus \calZ \cup \calP$ by charts of special holomorphic coordinates $z$ defines a flat structure on $X$ with transitions given by translations and half-turns. Such an atlas on $X$ is referred to as half-translation structure. In particular, there is a well-defined notion of (unoriented) lines in fixed direction $\theta \in \RR/\pi\RR$ with respect to $z$-coordinates. We refer to such line on $X$ as a \textit{trajectory in direction} $\theta$. \\
Similarly, we can obtain a flat structure on $\Sigma\setminus \calZ(\lambda) \cup \calP(\lambda)$ with transitions given by translations by choosing holomorphic coordinates $w$, such that $\lambda=w^l d w$. A Riemann surface with an atlas of charts glued by translations is referred to as translation surface. For a translation surface there is a well-defined notion of oriented lines in a fixed direction $\theta \in \RR/2\pi\RR$ with respect the $w$-coordinates. We refer to such lines on~$\Sigma$ as \textit{trajectories in direction} $\theta$. 
\par
\medskip
\paragraph{\textbf{Trajectories}}
Following \cite{Strebel} for a meromorphic quadratic differential $(X,q)$ trajectories in direction $\theta$ can be classified as follows. Here a trajectory is called \emph{critical} if it emanates from a point in $\calZ(q) \subset X$. 
\begin{itemize}
\item \textbf{Closed trajectories,} also know as as periodic trajectories. They come in families parameterizing \emph{ring domains} (also called \emph{cylinders}) with respect to the flat structure. The boundary of a ring domain is a union of saddle trajectories in direction~$\theta$, except with the ring domain is adjacent to a double pole on one side. The latter is called a \emph{degenerate ring domain}, the other cases are \emph{non-degenerate ring domains.}
\item \textbf{Non-critical, non-recurrent trajectories:} These trajectories connect two infinite critical points of the quadratic differential. They come in families of \textbf{flat strips} in direction $\theta$ biholomorphic to
  \bes S\=\{z \in \CC \mid a < \Im(z) <b \}.\ees
  On each boundary of each flat strip there is at least on zero of $q$. 
  \item \textbf{Recurrent trajectories:} The closure of a recurrent trajectory is a subsurface of $Y \subset X$. In this case, the interior of $Y$  is referred to as \textit{spiral domain}. The boundary of a spiral domain is a union of saddle connections in direction $\theta$.  
  \item \textbf{Saddle trajectories:} A saddle trajectory is a critical trajectory that connects two points in $\calZ(q)$. 
\end{itemize}
\par
If~$\gamma$ is a closed trajectory, then we may orient the components of $\pi^{-1}(\gamma)$ (one if~$\gamma$ is a saddle connections and two if~$\gamma$ is a core curve of a ring domain) so that the resulting class belongs to the charge lattice~$\Gamma$. It is called the \emph{hat-homology class~$\wh\gamma$ of~$\gamma$}, well-defined up to sign. Whenever we specify a choice of the abelian differential, we specify the hat-homology class uniquely by the requirement that $Z(\wh\gamma) \in \HH \cup \RR_{<0}$ (as in \cite{BridgelandSmith}).
\par
In the following, we will consider \emph{directed quadratic differentials} $(X,q,\theta)$, that is, $(X,q)$ is a meromorphic quadratic differential, a tacitly chosen abelian differential~$\lambda$ on~$\Sigma$ among the two roots of $\pi^* q$ and $\theta \in \RR/ 2\pi\RR$. Our results will not depend on the choice of~$\lambda$ as long as this choice is made consistently under deformation of~$(X,q,\theta)$, which we assume throughout. A trajectory of a directed quadratic differential thus refers to trajectory in direction $\theta$ on~$\Sigma$ with respect to $\lambda$, unless stated explicitly that we work with the unoriented trajectories on~$X$. A directed quadratic differential $(X,q,\theta)$ is called \textit{saddle-free}, if there are no saddle trajectories (in direction $\theta$).
\par
\medskip
\paragraph{\textbf{Triangulations}}
Consider a saddle-free directed differential $(X,q,\theta)$. Then the surface is decomposed in the direction~$\theta$ into strips with a single zero on either side of the strip. Choosing one representative of the parallel (non-critical, non-recurrent) trajectories in each strip defines an (ideal) triangulation~$\TT$ of~$X$, often considered also as a triangulation of the real blowup at the poles, see Section~\ref{sec:ciliated}. For each strip, i.e., for each edge~$e \in \TT$ there is a distinguished ('standard') saddle connection~$\gamma_e$ joining the two zeros and the boundary of the strip. We often use that the set $\{\gamma_e, e \in \TT\}$ is a basis of the hat-homology lattice~$\Gamma$ (see e.g.\ \cite[Lemma~3.2]{BridgelandSmith}).

\subsection{Spectral network $=$ Critical Trajectories}
\label{sec:spectralNW}

Let $(X,q,\theta)$ be a directed quadratic differential. The spectral network $W_\theta$ is a collection of oriented lines on the spectral cover $\Sigma$. It is given by the critical trajectories on $(\Sigma,\lambda)$ in direction $\theta$. By definition the trajectories are oriented, such that $\Re(e^{-i\theta}\lambda)$ increases along the orientation. Recall that each simple zero of $q$ corresponds to a zero of order $2$ of $\lambda$, from which there are emanating $6$ critical trajectories on $\Sigma$. Three of these trajectories are oriented away from the zero and three are oriented towards the zero.  We decompose $W_{\theta}=W_{\theta}^+ \cup W^-_{\theta}$, where $W^+_\theta$ is the subset of trajectories pointing away from the zero and $W^-_{\theta}$ is the subset of trajectories pointing in direction of the zero. Both $W^\pm_{\theta}$ contain $3r$ trajectories, though if $(X,q,\theta)$ has saddle trajectories in direction~$\theta$, these subsets have non-empty intersection.

Note that even though our definition coincides with the definition of the spectral network associated to a directed quadratic differential in \cite[Section~5]{GMN_spectral}, in general it does not define a spectral network subordinate to the cover $\Sigma \to X$ as defined in \cite[Section~9.1]{GMN_spectral}. We allow recurrent trajectories as part of our spectral network, which are embeddings of an half-open interval $[0,\infty)$ instead of a closed one. Furthermore around zeros at the boundary of a spiral domain, the trajectory structure is not described by condition C3 of \cite[Section~9.1]{GMN_spectral}.
  \par

\subsection{Circle bundle lifts}
%

In the sequel we need to distinguish two ways of lifting paths on a surface
$\Sigma$ (or on the complement $\Sigma^\circ$ of the zeros and poles) to the circle bundle
$c: \wt{\Sigma}^\circ \to \Sigma^\circ$:
\par
First, given a smooth oriented curve~$\gamma$ in $\Sigma^\circ$ there is the \emph{natural lift} $\vv{\gamma}$ to $\wt{\Sigma}^\circ$ given by $(\gamma(t),\arg(\dot{\gamma}(t)))$ for a smooth parametrization $\gamma: [0,1] \to \Sigma^\circ$ inducing the orientation. This will be used below in the definition of the lift function $F$.
\par
Second, we also would like to define the lift of a homology class in $H_1(\Sigma^\circ)$. Given a chain $\gamma = \sum_{i=1}^m \gamma_i$, where $\gamma_i$ are smooth pairwise disjoint simple closed curves, the \emph{canonical lift} is defined by
\[ \wt{\gamma}\ =\sum_{i=1}^m \vv{\gamma}_i + mS 
\]
where $S$ is fiber circle of $c: \wt{\Sigma}^\circ \to \Sigma^\circ$. It is shown in \cite{Johnson_Spin} that any homology class in $H_1(X)$ can be represented in this way and that the class $[\wt{\gamma}] \in \wt{H}^\abs$ is independent of this choice. Furthermore, by \cite[Theorem~1B]{Johnson_Spin} for $\alpha, \beta \in H_1(\Sigma^\circ)$ the canonical lift differs from a homomorphism by the law 
\be \label{eq:canlifthomo}
\wt{\alpha+\beta}= \wt{\alpha} + \wt{\beta} +\langle \alpha, \beta \rangle [S], 
\ee
where $\langle\ ,\ \rangle$ denotes the intersection pairing. 

\subsection{Quadratic differentials of finite area} \label{sec:qfinarea}

The content of this section is used only in Section~\ref{sec:FiniteArea} to understand the geometric significance of the wall-crossing formula in the finite area case. We start by recalling basic facts on IETs and measured foliations.
\par
For an integer~$n$ let $\Lambda$ be the set of positive real tuples $\alpha = (\alpha_1,\ldots,\alpha_n)$ with total sum one and $\pi \in S_n$ be a permutation. An \emph{interval exchange transformation (IET)} is a bijection $T = (\alpha,\pi)$ of the unit interval, defined by permuting~$n$ intervals of lengths~$\alpha_i$ according to~$\pi$. See \cite{Yoccoz} for a general introduction. An IET has a \emph{connection} if some break points $\beta_i = \sum_{j=0}^i \alpha_j$ have the same forward image $T^k(\beta_i) = T^\ell(\beta_j)$ for some $k,\ell>0$ and some~$i\neq j$. Let $\IET_n$ be the space of interval exchange transformations on $n$ intervals (with the topology inherited from $[0,1]^n$) and $\IET_n^*$ be the subspace of IETs without a connection.
\par
\begin{proposition} \label{prop:IETcomp}
  Let $T$ be an IET without a connection. The path connected component of~$T$ in $\IET_n^*$ is the interior of a closed simplicial cone~$C(T)$. This cone is homeomorphic to the space of invariant ergodic probability measures invariant under~$T$. Its dimension is bounded by the genus of the translation surface associated with~$T$. In particular, if~$T$ is uniquely ergodic, then its path connected component in $\IET_n^*$ is a point. 
\end{proposition}
\par
\begin{proof} The cone is the intersection of the simplices given by the Rauzy path associated with~$T$. See e.g.\ \cite[Section~8]{Yoccoz} for a proof. The proof shows that each boundary of the simplex is given by IETs with a connection.
\end{proof}
\par
The horizontal leaves of a quadratic differential $(X,q)$ define a foliation on~$\Sigma$ together with a transverse measure, a measured foliation~$\calF_\hor$. Similarly, the vertical foliation defines a measured foliation $\calF_\ver$. Jointly, $(\calF_\hor,\calF_\ver)$ is \emph{filling} \cite[Lemma~5.3]{GardinerMasur}, i.e.\ the sum of intersection numbers with any other measured foliation is positive. The converse is given essentially in loc.~cit., see e.g.\ also in \cite[Theorem~1.6]{Wentworth}.
\par
\begin{proposition} \label{prop:fillingfoliations}
  If $(\calF_\hor,\calF_\ver)$ is filling, then there is a unique
  quadratic differential with horizontal and vertical foliation given by $(\calF_\hor,\calF_\ver)$. 
\end{proposition}
\par
To connect the two topics, take the canonical double cover~$(\Sigma,\lambda)$ of~$(X,q)$. Take any interval~$I$ transverse to the horizontal foliation of~$\lambda$. Then the first return map to~$I$ is an IET~$T$. To fit with our convention of IETs, the result has to be rescaled to the unit interval. (To make the interval canonical up to a finite number of choices one often takes a separatrix at a zero of~$\lambda$ in the perpendicular direction with length normalization dictated by Rauzy induction.) The dimension of the space of invariant measures is a property of the horizontal foliation and thus does not depend on the choice of~$I$. Moreover, the IET has a connection if and only if the horizontal foliation has a saddle connection.
\par
The existence of an IET which is minimal and not uniquely ergodic is a non-trivial result that goes back to Keane \cite{Keane}. A systematic study of the set of non-uniquely ergodic directions in terms of Hausdorff dimension started with work of Masur-Smillie \cite{MasurSmillie} who showed that for every surface the set of minimal and non-uniquely ergodic directions has positive Hausdorff dimension, despite being 'small' in the sense that almost every direction is uniquely ergodic. The computation of Hausdorff dimension recently culminated in the following statement of Chaika-Masur \cite{ChaikaMasur}:
\par
\begin{proposition} \label{prop:notUErare}
For every~$n$ the subspace of uniquely ergodic IETs in $\IET_n$ has Hausdorff codimension~$1/2$.
\end{proposition}

\section{Formal completions of groupoid modules} \label{sec:formal}
In this section we define the module  $\wh{\Pi} = \wh\Pi^ \sign_{\Delta}(\wt{\Sigma}^\circ)$ of \emph{completed module of signed paths}. The homological variant $\wh{H} = \wh{H}^\sign_{\Delta}(\wt{\Sigma}^\circ)$ is defined in parallel. The first goal is to show in Proposition~\ref{prop:rngstructure} that the natural multiplication by path composition makes certain 'tame' submodules into rngs (without unit). Key to the sequel is some elementary cone geometry in Lemma~\ref{le:translatebycone}, which is use in Lemma~\ref{le:Kgeneralcont} to show continuity of the BPS-automorphism~$\calK_\theta$ on an exhaustive sequence of subspaces.
\par

\subsection{The twisted torus and completed Poisson subalgebras}
\label{sec:twistedtorus}

This section is a summary of the setup of the wall-crossing formula, see \cite{KS08} or \cite[Appendix~A]{BridgelandRHDT}. Proofs of various claims here are special cases of the arguments in the generalizations in Section~\ref{sec:groupoidmod}.
\par
Let $(\Gamma, \langle \cdot,\cdot \rangle)$ be any lattice, in applications typical the charge lattice~\eqref{eq:chargelattice} or the full $H_1(\Sigma,\calP(\lambda),\ZZ)$. Associated with~$\Gamma$, there is the (usual) torus
\be
\TT_+ \= \TT_+(\Gamma) \= \Hom(\Gamma, \CC^*)
\ee
with the multiplicative group structure. We will also need the \emph{twisted torus} defined as the set of maps
\be \label{eq:twistedtorusmult}
\TT_- \= \TT_-(\Gamma) \= \{g: \Gamma \to \CC^*\,,\,\, g(\gamma_1+\gamma_2) = (-1)^{\langle \gamma_1, \gamma_2 \rangle} g(\gamma_1)g(\gamma_2)\}\,.
\ee
The twisted torus is a torsor under~$\TT_+$, i.e.\ there is a free and transitive action of $\TT_+$ on~$\TT_-$ given by $(fg)(\gamma) = f(\gamma)g(\gamma)$. This action makes~$\TT_-$ into an algebraic variety, independently of the choice of a base point.
\par
Fix an acute sector $\Delta \subset \CC$ and suppose we are moreover given a linear map ('central charge') $Z: \Gamma \to \CC$.  Inside the coordinate ring $\CC[\TT_-]$ we define the subalgebra
\be
\CC_\Delta[\TT_-] \= \bigoplus_{\alpha: Z(\alpha) \in \Delta} \CC \cdot [\alpha] \qquad \subset \CC[\TT_-]\,.
\ee We interpret $[\alpha]$ as an element of $\CC[\TT_-]$ via evaluation.
The bracket
\be
\{[\alpha],[\beta]\} \= \langle \alpha,\beta \rangle \cdot [\alpha] \cdot [\beta]
\ee
makes $\CC_\Delta[\TT_-]$ into a Poisson-subalgebra of $\CC[\TT_-]$. We define the \emph{height} of an element $a = \sum a_\alpha [\alpha] \in \CC[\TT_-]$ to be
\be \label{eq:height}
h(a) \= \inf \bigl\{ |Z(\alpha)|  \colon \alpha \,\, \text{such that $a_\alpha \neq 0$} \bigr\}\,.
\ee
For any $L>0$ the subspace $\CC_\Delta[\TT_-]_{\geq L}$ that consists of elements of height $\geq L$ is a Poisson ideal (since $\Delta$ is acute) and we can therefore form the quotients and take the projective limit as $L \to \infty$, i.e.\ we let
\ba
\CC_\Delta [[\TT_-]] &\= \varprojlim_L \, \CC_\Delta[\TT_-]_{< L}, \qquad \text{where}
\qquad \\
\CC_\Delta[\TT_-]_{< L} &\= \CC_\Delta[\TT_-]/ \CC_\Delta[\TT_-]_{\geq L}
\ea
\par
\medskip
\paragraph{\textbf{The Lie algebra $\wh\frakg_\Delta$ and the Lie group $\wh{G}_\Delta$}}.
We associate with $(\Gamma, \langle \cdot,\cdot \rangle, Z, \Delta)$ and under the standing hypothesis that $\Delta$ is acute the Lie algebra
\be
\wh\frakg_\Delta \= \{a \in \CC_\Delta [[\TT_-]] \colon h(a) > 0 \}
\ee
with  Lie bracket given by the Poisson bracket. This is a pro-nilpotent Lie algebra, since $\wh\frakg_\Delta = \varprojlim_L  \frakg_{\Delta, < L}$ where
\bes
\frakg_{\Delta \leq L} \= \{a \in \CC_\Delta [[\TT_-]]_{< L} \colon h(a) > 0 \}
\ees
and since the Lie bracket eventually increases $|Z(\cdot)|$ (compare with~\eqref{equ:inequ_tameness} below) and is thus nilpotent. Finally we define the Lie group
\be
\wh{G}_\Delta \= \{ \exp(a) \colon a \in \wh\frakg_\Delta\}
\ee
of formal symbols $\exp(a)$, with group structure derived from the Lie algebra structure on $\wh\frakg_\Delta$ via the Baker-Campbell-Hausdorff formula. This Lie group is pro-unipotent, in fact the projective limit $\wh{G}_\Delta = \varprojlim_L \wh{G}_{\Delta,<L}$ of the Lie groups generated by the formal symbols $\exp(a)$ for $a \in \frakg_{\Delta, < L}$.

\subsection{Groupoid modules} \label{sec:groupoidmod}

We generalize the previous setup in three aspects. Most important, we allow for central charges not just in an acute sector but in a translated sector. The amount of translation will be arbitrary, but uniform for each element in the completion, encoded in the notion of $\Delta$-tame paths. Second, the twists as in the torus will appear indirectly as quotients of paths in a circle bundle by some winding ideal. Third, we will simultaneously deal with the homology classes of paths, thus generalizing the twisted torus, and homotopy classes of path, a further generalization to be used for non-abelianization.
\par
\medskip
\paragraph{\textbf{Groupoids}}
Let $Y$ be a topological space. Consider the path groupoid $\calP(Y)$, whose objects are points $y \in Y$ and whose morphisms $\Mor(y_1,y_2)$ are continuous paths starting at $y_1$ and ending at $y_2$.\footnote{We use the  convention that $\gamma_2 \circ \gamma_1$ denotes the path where~$\gamma_1$ is taken first, then $\gamma_2$.}
We obtain the fundamental groupoid $\pi_{\leq 1}Y$ by quotienting the morphism by homotopy of paths relative to the endpoints. Furthermore, we define the first homology groupoid $H_{\leq 1}(Y)$ with the same objects $y \in Y$ and whose morphisms $\Mor(y_1,y_2)$ are 1-chains with boundary $y_2-y_1$ modulo 1-boundaries. Eventually we'll drop~$Y$ and simply write $\pi_{\leq 1}$ and $H_{\leq 1}$ for these groupoids.
\par
\medskip
\paragraph{\textbf{Groupoid modules and central charge truncation}} 
For each of the groupoids above we consider the sum of the~$\CC$-modules
generated by the morphisms, that is, we let
\be
\Pi_{\leq 1}(Y) \=  \bigoplus\limits_{\alpha \in \Mor_{{\pi_{\leq 1}}Y}} \CC [\alpha],
\qquad
H_{\leq 1}(Y) \=  \bigoplus\limits_{\alpha \in \Mor_{{H_{\leq 1}}(Y)}} \CC [\alpha]\,.
\ee
We write elements in these modules typically as \emph{finite} sums
\bes
a \= \sum_{\alpha \in \Mor_{{\pi_{\leq 1}}Y}} \coeff(\alpha) \cdot [\alpha], \quad
\text{resp.} \quad
a \= \sum_{\alpha \in \Mor_{{H_{\leq 1}}(Y)}} \coeff(\alpha) \cdot [\alpha],
\qquad \coeff(\alpha) \in \CC\,.
\ees The concatenation of paths extends to a rng structure on $\Pi_{\leq 1}(Y)$ respectively $H_{\leq 1}(Y)$ defined by 
\ba \label{eq:pathmult}
  [\alpha_1] \cdot [\alpha_2] \= \left\{ \begin{matrix}
[\alpha_2 \circ \alpha_1] & & \text{if $\alpha_1$  and $\alpha_2$ are composable,} \\
    0 & & \text{otherwise.}
  \end{matrix} \right. 
  \ea
 Suppose that is  a ('central charge') homomorphism of local systems on $Y \times Y$
 \ba \label{eq:central_charge_hom}
 Z: H_1(Y,\{\cdot,\cdot\},\ZZ) \to \CC_{Y \times Y}\,.
 \ea
where the left hand side is the local system whose fiber over $(y_1,y_2)$
is the relative cohomology $H_1(Y,\{y_1,y_2\},\ZZ)$ and the right hand side is the constant local system. We extend the notion of height $h$ defined in \eqref{eq:height} in the obvious way using $Z$. For each $L >0$ define the submodules of elements with height bounded below by~$L$
\[ \Pi_{\leq 1}(Y)_{\geq L} \= \Bigl\{  a = \sum_{\alpha \in \Mor_{{\pi_{\leq 1}}Y}} \coeff(\alpha)\cdot [\alpha]\in \Pi_{\leq 1}(Y)
\mid h(a) \geq L \Bigr\} 
\] (note that $0$ has infinite height) and the quotient
\[ \Pi_{\leq 1}(Y)_{<L} \= \Pi_{\leq 1}(Y) /\Pi_{\leq 1}(Y)_{\geq L}.
\] 
For $L'>L>0$ there are obvious truncation maps $t_{L'L}: \Pi_{\leq 1}(Y)_{<L'} \to \Pi_{\leq 1}(Y)_{<L}$ and hence we may define the projective limit
\[ \wh{\Pi}(Y) \= \varprojlim_L  \Pi_{\leq 1}(Y)_{<L}\,.
\]
Similarly, we define $H_{\leq 1}(Y)_{\geq L}$, the quotient $H_{\leq 1}(Y)_{< L}$ and finally
\be \wh{H}(Y) \=  \varprojlim_L  H(Y)_{<L} \,. \ee
The following useful lemma follows directly from the definition:
\par
\begin{lemma}\label{lemm:description_of_elements}
Elements of these formal completions $\wh{\Pi}(Y)$ and $\wh{H}(Y)$ are precisely the elements in the product of the modules $\CC[\alpha]$ that can be represented as formal sums 
\be
\sum_{\alpha \in \Mor_{\pi_{\leq 1}}} \coeff(\alpha) \cdot [\alpha] \in \prod\limits_{\alpha \in \Mor_{\pi_{\leq 1}}} \CC [\alpha], \quad
\text{resp.} \quad
 \sum_{\alpha \in \Mor_{H_{\leq 1}}} \coeff(\alpha)\cdot [\alpha] \in \prod\limits_{\alpha \in \Mor_{H_{\leq 1}}} \CC [\alpha],
\ee
such that for each $L>0$ there are non-zero coefficients $\coeff(\alpha) \neq 0$ for only finitely many $[\alpha]$ with $\vert Z(\alpha)\vert <L$.
\end{lemma}
\par
For an element $a = \sum_\alpha \coeff(\alpha) [\alpha]$ in  $\wh{\Pi}(Y)$ or  $\wh{H}(Y)$ we denote by
\be
S(a) \= \{ \alpha: \coeff(\alpha) \neq 0\}
\ee
the \emph{support} of~$a$.
\par
Note that path concatenation does not define a ring structure on $\wh{\Pi}(Y)$, since for $\alpha \circ \beta = \gamma$ there is no upper bound on $|Z(\alpha)|$ and $|Z(\beta)|$ in terms of $|Z(\gamma)|$. 
              
\subsection{$\Delta$-submodules and their ring structure}

We now apply the previous discussion to the case~$Y= \wt{\Sigma}^\circ$ of a circle bundle $c: \wt{\Sigma}^\circ \to \Sigma^\circ$ over a surface. Again we fix an acute sector $\Delta$. We suppose that the surface~$\Sigma^\circ$ comes with a holomorphic one-form~$\lambda$ and that the central charge~$Z$ is given by integrating~$\lambda$. We define $Z(\gamma) = Z(c \circ \gamma)$ for any path~$\gamma$ in~$\wt\Sigma^\circ$. 
\par
Inside the module $\wh{\Pi}(\wt{\Sigma}^\circ)$  we now define the submodule $\wh{\Pi}_{\Delta}(\wt{\Sigma}^\circ)$ of $\Delta$-tame paths as the set of paths whose central charge is valued in cone that is a translate of the sector~$\Delta$. More precisely, for $t \in \CC$ we define the cone $\Delta^t = t + \Delta \subset \CC$ and let
\ba\label{eq:deftame}
\wh{\Pi}_{\Delta^t}(\wt{\Sigma}^\circ) &\= \{ a \in \wh{\Pi}(\wt{\Sigma}^\circ) \,\mid\,  \coeff(\alpha) \neq 0 \Rightarrow Z(\alpha) \in \Delta^t \},  \\
\wh{H}_{\Delta^t}(\wt{\Sigma}^\circ) &\= \{ a \in \wh{H}(\wt{\Sigma}^\circ) \,\mid\,  \coeff(\alpha) \neq 0 \Rightarrow Z(\alpha) \in \Delta^t \}\,.
\ea
Note that the set of translated sectors with inclusion defines a directed set. We define the \emph{$\Delta$-tame submodules} of $\wh{\Pi}(\wt{\Sigma}^\circ)$ and $\wh{H}(\wt{\Sigma}^\circ)$ as the direct limit
\ba\label{eq:tameisdirectlimit}
\wh{\Pi}_{\Delta}(\wt{\Sigma}^\circ) &\= \varinjlim_{\Delta^t}\, \wh{\Pi}_{\Delta^t}(\wt{\Sigma}^\circ),  \\
\wh{H}_{\Delta}(\wt{\Sigma}^\circ) &\=  \varinjlim_{\Delta^t} \, \wh{H}_{\Delta^t}(\wt{\Sigma}^\circ)\,.
\ea
\par
\begin{proposition} \label{prop:rngstructure}
The natural multiplication \eqref{eq:pathmult} makes $\wh{\Pi}_{\Delta}(\wt{\Sigma}^\circ)$ and also $\wh{H}_{\Delta}(\wt{\Sigma}^\circ)$ into a rng (a ring without unit). More precisely for $t_1,t_2 \in \CC$ the multiplication induces maps with target
  \bas
  &\circ: \wh{\Pi}_{\Delta^{t_1}}(\wt{\Sigma}^\circ) \times \wh{\Pi}_{\Delta^{t_2}}(\wt{\Sigma}^\circ) \to \wh{\Pi}_{\Delta^{t_1+t_2}}(\wt{\Sigma}^\circ), \\
   &\circ: \wh{H}_{\Delta^{t_1}}(\wt{\Sigma}^\circ) \times \wh{H}_{\Delta^{t_2}}(\wt{\Sigma}^\circ) \to \wh{H}_{\Delta^{t_1+t_2}}(\wt{\Sigma}^\circ)\,.
\eas
given by elements with support in $\Delta^{t_1+t_2}$.
\end{proposition}
\par
(The unit 'should' be the sum of the trivial paths at every point, which however does not belong to these modules.)
\par
The following elementary result about the sum of vectors in a translated cone will be used in many instances below. Let $t \in \CC$ and $Z_1,Z_2 \in \Delta^t$, then
\ba  \label{equ:inequ_tameness}\max \{ \vert Z_1 \vert, \vert Z_2 \vert \}- \vert t \vert \leq \vert Z_1 + Z_2 \vert +2 \vert t \vert\,.
\ea
This follows from the inequality $\max \{ \vert Z_1 \vert, \vert Z_2 \vert \} \leq \vert Z_1 + Z_2 \vert$ for $Z_1,Z_2 \in \Delta=\Delta^0$ by translation and triangle inequality. 
\par
\begin{proof} Let~$t \in \CC$ so that~\eqref{eq:deftame} applies for both elements. We write the two elements as 
\bes  a\= \sum_{\alpha \in \Mor_{\pi_{\leq 1}}} \coeff(\alpha) \cdot [\alpha], \qquad
b\ = \sum_{\beta \in \Mor_{\pi_{\leq 1}}} \coeffb(\beta) \cdot [\beta] 
\,\in\,\wh{\Pi}_\Delta(\wt{\Sigma}^\circ)\,.
\ees
We would like to define their product to be
\[
c \,:=\,  a \cdot b  \= \sum_{\gamma \in \Mor_{\pi_{\leq 1}}} \Bigl(\sum_{[\alpha] \cdot [\beta] = [\gamma]} \coeff(\alpha)\coeffb(\beta) \Bigr)[\gamma]
\]
and we have to show first that the coefficient of each $[\gamma]$ is finite. By~\eqref{equ:inequ_tameness} the absolute values of central charge $|Z(\alpha)|,\vert Z(\beta) \vert$ are bounded by $L=|Z(\gamma)|+3\vert t\vert $. By definition of the completion there are only finitely many ~$\alpha,\beta$ with $\coeff(\alpha),\coeffb(\beta)\neq 0$ and $\vert Z(\alpha)\vert,\vert Z(\beta)\vert \leq L$. Only these finitely many elements can contribute to the prefactor of a given~$[\gamma]$. Hence the coefficient is finite. 
\par
Second, from inequality~\eqref{equ:inequ_tameness} it is clear that for each given $L>0$, there are only finitely many paths $[\gamma]$ with $Z(\gamma) \leq L$ and non-zero coefficient in the product $c=a\cdot  b$.  Hence, $c$ defines an element of the formal completion by Lemma~\ref{lemm:description_of_elements}.
\par
Finally, as the cone $\Delta$ is closed under sums, the product $c=a\cdot b$ is $\Delta$-tame. If the central charges of the paths appearing in $a$ and $b$ are valued in $\Delta^{t_1}$ and $\Delta^{t_2}$ respectively, then the central charges of paths appearing in $c$ are valued in $\Delta^{t_1+t_2}$.
\end{proof}
\par
\begin{corollary} \label{cor:infiniteproduct}
  Let $a_i = \sum_\alpha \coeff_i(\alpha) [\alpha]$ for $i \in I$ be a collection of elements in the $\Delta$-tame submodule $\wh{\Pi}_\Delta(\wt{\Sigma}^\circ)$ resp. $\wh{H}_\Delta(\wt{\Sigma}^\circ)$  , such that the sum of translates is bounded, i.e.\ $a_i$ is $\Delta$-tame for the translated cone $\Delta^{t_i}$ with $\vert \sum_{i\in I} t_i \vert < \infty$.
  The joint support $S=S(a_i\mid i \in I)$ is the multi-set defined as the (big) union $\cup_{i \in I} \{\alpha: \coeff_i({\alpha}) \neq 0\}$. Assume further that the joint support satisfies that 
\be
|\{\alpha \in S \colon |Z(\alpha)| < L \}| < \infty
\ee
for any~$L \in \RR$.  Then the infinite product $\prod_{i \in I} a_i$ is a well-defined element in $\wh{\Pi}_\Delta(\wt{\Sigma}^\circ)$ resp. $\wh{H}_\Delta(\wt{\Sigma}^\circ)$.
\end{corollary}
\par
\begin{proof} This is an immediate extension of the arguments of the previous proposition. From the condition on the support and inequality~\eqref{equ:inequ_tameness} we conclude that the product is well-defined as an element of the formal completion. From the condition on the translates the infinite product is $\Delta$-tame for the translated cone $\Delta^{\sum t_i}$.
\end{proof}
\par
\medskip
\paragraph{\textbf{Topology}} Let $Y$ be a topological space with a homomorphism $Z$ as in \eqref{eq:central_charge_hom}. We provide $\wh{\Pi}(Y)$ and $\wh{H}(Y)$ with the following topology, which is reminiscent of the $I$-adic topologies. A basis of open sets is define by
\[ \{ a + \wh{\Pi}(Y)_{> L} \} \subset \wh{\Pi}(Y) \quad \text{resp.} \quad \{ a + \wh{H}(Y)_{> L} \} \subset \wh{H}(Y),
\] where $\wh{\Pi}(Y)_{> L} \subset \wh{\Pi}(Y)$ and $\wh{H}(Y)_{> L} \subset \wh{H}(Y)$ are the submodules of elements $a$ with height $h(a) > L$.  This topology is induced by the metric
\[ d(a,b) = 2^{-h(a-b)}, 
\] where $h$ is the height defined in equation \eqref{eq:height}. We refer to this topology as \emph{$Z$-adic topology}. By Lemma \ref{lemm:description_of_elements} the projective limit $\wh{\Pi}(Y)$ is the metric completion of $\Pi_{\leq 1}(Y)$. The same holds for the homological version. Now we again specialize to $Y=\wt{\Sigma}^\circ$ and the central charge $Z$.
\par
\begin{lemma} \label{le:restmultcont}
The subspace $\wh{\Pi}_{\Delta^t}(\wt{\Sigma}^\circ)$ of the $\Delta$-tame submodule is the completion of
  \[ \{ a \in \Pi_{\leq 1}(\wt{\Sigma}^\circ) \mid \coeff(\alpha)\neq 0 \Rightarrow Z(\alpha) \in \Delta^t \} 
  \]
  with respect to the $Z$-adic metric. Addition on $\wh{\Pi}(\wt{\Sigma}^\circ)$ and the  multiplication restricted to the subspaces
\be
\wh{\Pi}_{\Delta^t}(\wt{\Sigma}^\circ)\times \wh{\Pi}_{\Delta^t}(\wt{\Sigma}^\circ) \to \wh{\Pi}_{\Delta^{2t}}(\wt{\Sigma}^\circ)\ee
are continuous in the $Z$-adic topology. Similar statements hold for the homological version.
\end{lemma}
\par
\begin{proof}
  The continuity of the addition map is clear. The sequential continuity of the multiplication map follows from inequality~\eqref{equ:inequ_tameness}. In fact, consider converging sequences $a_n\to a , b_n\to b \in \wh{\Pi}_{\Delta^t}(\wt{\Sigma}^\circ)$ with respect to the $Z$-adic metric. By definition for each $L>0$ there exists $N \in \NN$, such that for each $n>N$ the differences $a_n-a$ and $b_n-b$ have height bounded below by $L$. Hence, $a_n\cdot b_n-a\cdot b$ has height bounded below by $L- 3\vert t\vert$ by inequality~\eqref{equ:inequ_tameness}.
\end{proof}
\par
More generally, in Section~\ref{sec:A0laminations} we need to allow sets of translates. For any subset $K \subset \CC$, let
\bas \wh{\Pi}_{\Delta^K}(\wt{\Sigma}^\circ) &\= \{ a \in \wh{\Pi}(\wt{\Sigma}^\circ) \,\mid\,  \coeff(\alpha) \neq 0 \Rightarrow Z(\alpha) \in \Delta+K \subset \CC \},  \\
\wh{H}_{\Delta^K}(\wt{\Sigma}^\circ) &\= \{ a \in \wh{H}(\wt{\Sigma}^\circ) \,\mid\,  \coeff(\alpha) \neq 0 \Rightarrow Z(\alpha) \in \Delta+K \subset \CC\},
\eas
\begin{lemma} \label{le:translatebycone}
  Let $\Delta_1,\Delta_2\subset \CC$ be two cones, such that $\Delta_2$ intersects the opposite cone~$-\Delta_1$ of $\Delta_1$ only in the origin, i.e.\
\be \label{eq:oppcone}
\Delta_{1,-} \cap \Delta_2 = \{0\}, \qquad \text{where}\qquad
\Delta_{1,-}=\{ -Z \mid Z \in \Delta_1\} \subset \CC.
\ee
Let $t \in \CC$ and $K= \Delta_2^t$. Then the multiplication defined in \eqref{eq:pathmult} induces on $\wh{\Pi}_{\Delta_1^K}(\wt{\Sigma}^\circ)$ the structure of a topological $\wh{\Pi}_{\Delta_1^0}(\wt{\Sigma}^\circ)$-module.
\end{lemma}
\begin{proof} Let
  \bes  a\= \sum_{\alpha \in \Mor_{\pi_{\leq 1}}} \coeff(\alpha) \cdot [\alpha] \,\in\, \wh{\Pi}_{\Delta_1^0}(\wt{\Sigma}^\circ) \quad \text{and} \quad 
b\ = \sum_{\beta \in \Mor_{\pi_{\leq 1}}} \coeffb(\beta) \cdot [\beta], 
\,\in\, \wh{\Pi}_{\Delta_1^K}(\wt{\Sigma}^\circ)\,.
\ees
Let $[\alpha],[\beta],[\gamma] \in \Mor_{\pi_{\leq 1}}$ with $Z(\alpha) \in \Delta_1$, with $Z(\beta) \in \Delta_1^K$ and $\vert Z(\gamma) \vert \leq L$. If $[\alpha][\beta]=[\gamma]$, then $Z(\beta) \in B_L(0)+\Delta_{1,-}$, where $B_L(0)$ is closed disc of radius $L$ around $0 \in \CC$ with respect to the euclidean norm. Elementary geometry shows that~\eqref{eq:oppcone} implies that
\be \label{eq:conearg} 
\Delta_1^K \cap (B_L(0)+\Delta_{1,-}) \quad \text{is a bounded region in $\CC$.}
\ee
Hence, it is contained in another ball $B_{L'}(0)$. Now, the triangle inequality implies
\be L \geq \vert Z(\gamma)\vert \geq \vert Z(\alpha)\vert-\vert Z(\beta)\vert \geq \vert Z(\alpha)\vert-L'. \label{equ:beta_versus_gamma}
\ee In particular, there are only finitely many $[\alpha]$ with $\coeff(\alpha) \neq 0$ and $[\beta]$ with $\coeffb(\beta) \neq 0$, such that $\vert Z(\beta \circ \alpha)\vert <L$. Hence, the multiplication extends to a well-defined map. 
\par
Continuity follows along the same lines: Contrapositions yield for each~$L$ the  existence of an $L'$, such that for $\vert Z(\alpha)\vert >L'+L$ or for $\vert Z(\beta) \vert > L'$ (by inequality~\eqref{equ:beta_versus_gamma}) we have $\vert Z(\beta \circ \alpha) \vert > L$. Hence, the multiplication map is  sequentially continuous.
\end{proof}
\par

\subsection{The quotient by the winding submodule}
Let $\pi_1(c) : \pi_1(\wt{\Sigma}^\circ) \to \pi_1(\Sigma^\circ)$ be the homomorphism induced by~$c$ on fundamental groups. Its kernel is normally generated by a non-trivial loop~$S$ in the circle and we define 
\be
w : \ker(\pi_1(c)) \to \ZZ/2\ZZ, \qquad S \mapsto 1
\ee
to be the \emph{winding number mod~$2$}. We define the \emph{winding submodule~$I_{\wind}$} to be the ideal in $\Pi_{\leq 1}(\wt{\Sigma}^\circ)$ and $H_{\leq 1}(\wt{\Sigma}^\circ)$ generated by
\be
\langle [\wp_1] - (-1)^{w(\wp_1 \circ \wp_2^{-1})} [\wp_2] \mid [\wp_1],[\wp_2] \in \Mor(x_1,x_2):  \wp_1 \circ  \wp_2^{-1} \in \ker(\pi_1(c))  \rangle\,. 
\ee
(Note that e.g.\ by \cite[Lemma~2.4]{CalderonSalter} the winding number of a boundary of a subsurface is its Euler characteristic. Hence the winding number $\mod 2$ is well-defined on homology classes.) Abusing notation, we denote by $\wh{I}_{\wind}$ the $Z$-adic closure of $I_{\wind}$ and by $\wh{I}_{\wind,\Delta}$ both the intersections
\[ \wh{I}_{\wind,\Delta} := \wh{I}_{\wind} \cap \wh{\Pi}_\Delta(\wt{\Sigma}^\circ) \quad \text{resp.} \quad  \wh{I}_{\wind,\Delta} := \wh{I}_{\wind} \cap \wh{H}_\Delta(\wt{\Sigma}^\circ)\,.
\]
Taking the quotients by these submodules will be essential in the homotopy invariance of the path lifting function in Section~\ref{sec:signed_hom_inv}.
\par
\begin{definition}
  We define 
\ba
\wh{\Pi}_\Delta &\,:=\, \wh{\Pi}^{\sign}_{\Delta}(\wt{\Sigma}^\circ) := \wh{\Pi}_\Delta(\wt{\Sigma}^\circ) / \wh{I}_{\wind,\Delta}, \\
\wh{H}_\Delta&\,:=\, \wh{H}^\sign_{\Delta}(\wt{\Sigma}^\circ) := \wh{H}_{\Delta}(\wt{\Sigma}^\circ) / \wh{I}_{\wind,\Delta}
\ea
the \emph{completed homotopy (resp.\ homology) groupoid modules of signed paths}. 
\end{definition}
\par
As above we define $\wh{\Pi}_{\Delta^t}$ and $\wh{H}_{\Delta^t}$ to be the images
of $\wh{\Pi}_{\Delta^t}(\wt{\Sigma}^\circ)$ and $\wh{H}_{\Delta}(\wt{\Sigma}^\circ)$ respectively in these quotients by the winding ideal. This definition implicitly uses:
\par
\begin{lemma} \label{le:quotientring}
The submodule $\wh{I}_{\wind,\Delta} \subset \wh{\Pi}_\Delta(\wt{\Sigma}^\circ)$ is an ideal. In particular, the quotients $\wh{\Pi}_\Delta$ and $\wh{H}_\Delta$ inherit from Proposition~\ref{prop:rngstructure} the structure of a rng. Addition is continuous and for any $t \in \CC$ the restriction of the multiplication to $\wh{\Pi}_{\Delta^t}$ resp.\ to $\wh{H}_{\Delta^t}$ is continuous. \par
For $t \in \CC$ and for $K$ and the remaining notation as in Lemma~\ref{le:translatebycone}, the multiplication in the ring induces on $\wh{\Pi}_{\Delta_1^K}$ the structure of a topological $\wh{\Pi}_{\Delta_1^0}$-module and it induces on $\wh{H}_{\Delta_1^K}$ the structure of a topological $\wh{H}_{\Delta_1^0}$-module.
\end{lemma}
\par
\begin{proof} The first statement would be obvious, if multiplication in $\wh{\Pi}_\Delta(\wt{\Sigma}^\circ)$ were continuous. The continuity of the restricted multiplication in Lemma~\ref{le:restmultcont} suffices as substitute: Using the representation in Lemma~\ref{lemm:description_of_elements} any element  $i \in \wh{I}_{\wind,\Delta}$ can be written as a limit of elements in $I_{\wind}$ with support in~$\Delta^t$ for some~$t$. Similarly, we may write any element $a \in \wh{\Pi}_\Delta(\wt{\Sigma}^\circ)$ as limit of elements with support in~$\Delta^{t_2}$ for some~$t_2$. Now the restricted continuity gives the claim.
\par
The remaining statements are a consequence of the first together with Lemma~\ref{le:restmultcont} and Lemma~\ref{le:translatebycone}.
\end{proof}
\par
Finally we observe that for elements~$a$ in $\wh{\Pi}_\Delta$ or$\wh{H}_\Delta$, i.e.\ for cosets modulo the winding ideal, the support itself is not well-defined, but the set of central charges $Z(S(a)) \subset \CC$ is well-defined.
\par
\medskip
\paragraph{\textbf{Absolute paths}} For essentially every object introduced in this section, there is the corresponding subobject given or generated by paths in $\wt{\Sigma}^\circ$ whose image in~$X$ is closed and smooth (including at the end points if considered as the image of~$S^1$ rather than as the image of an interval). We denote these subobjects with a superscript 'abs' and refer to them as \emph{absolute} (as opposed to relative) paths or cycles, since they do not need to be closed by may also be paths starting at $\wt{z}$ and ending at $\sigma(\wt{z})$ where $\sigma$ denotes the lift of the involution on $\wt\Sigma$ to the circle bundle. In particular $\wh{\Pi}_\Delta^\abs$ and $\wh{H}_\Delta^\abs$ denotes the \emph{submodule of absolute paths} in the completed homotopy (resp.\ homology) groupoid modules of signed paths.
\par
The argument from Proposition~\ref{prop:rngstructure} makes $\wh{H}^\abs_\Delta$ with multiplication induced by $\wh{H}_\Delta(\wt{\Sigma}^\circ)$
\be \label{eq:closedpathmult}
[a] \cdot [b] \= [a+b], \qquad a,b \in \wh{H}_1(\wt{\Sigma}^\circ)
\ee
into a ring (with unit the trivial path at any point). Moreover the same arguments as for Proposition~\ref{prop:rngstructure} show:
\begin{lemma}
The multiplication \eqref{eq:closedpathmult} makes $\wh{H}_\Delta(\wt{\Sigma}^\circ)$ into a $\wh{H}^\abs_\Delta(\wt{\Sigma}^\circ)$-module and makes $\wh{H}_\Delta$ into a $\wh{H}^\abs_\Delta$-module.
\end{lemma}
\par
Observe that the multiplication rules~\eqref{eq:closedpathmult} and~\eqref{eq:pathmult} are compatible (making $\wh{H}_\Delta$ into a $\wh{H}^\abs_\Delta$-algebra), justifying our ambiguous notation. Finally, we pinpoint in which sense the completed groupoid modules generalize the completed coordinate rings of twisted tori.
\par
\begin{proposition} \label{prop:containstwistedtorus}
  Fix an acute sector~$\Delta$. For $\Gamma = H_1(\Sigma,\calP(\lambda),\ZZ)$ the completed algebra $\CC_\Delta[[\TT_-]]$ is a closed subalgebra of $\wh{H}^\abs_\Delta$, the subalgebra given by the trivial translate~$t=0$
in the direct limit~\eqref{eq:tameisdirectlimit} modulo the winding ideal.
\par
If $\Sigma$ admits an involution~$\sigma$ preserving $\Sigma^\circ$ with associated quotient map the double cover $p:\Sigma \to X$, then the action of $\sigma$ induces an action on $\wh{H}_\Delta$, whose invariant and anti-invariant subspaces we denote by $\wh{H}_\Delta^{\abs,\pm}$. If we consider $\Gamma = H_1(\Sigma,\calP(\lambda),\ZZ)^-$, then the completed algebra $\CC_\Delta[[\TT_-]]$ is a subalgebra of $H_\Delta^{\abs,-}$, the subalgebra given by the trivial translate~$t=0$ in the direct limit~\eqref{eq:tameisdirectlimit} modulo the winding ideal.
\end{proposition}
\par
In this proposition the algebra structure on $\CC_\Delta[[\TT_-]]$ is the one coming from the twisted torus. 
\par
\begin{proof} The main observation is that the map
\be
\CC_\Delta[\TT_-]_{<L} \to \wh{H}^\abs_\Delta, \qquad  [\alpha] \mapsto [\wt{\alpha}]
\ee
sending a path to its canonical lift induces an algebra homomorphism in the projective limit as $L \to \infty$: For any $\alpha,\beta \in \CC_\Delta[\TT_-]$ if $L$ is large enough, then the multiplication in the twisted torus in~\eqref{eq:twistedtorusmult} implies that
\be
[\alpha + \beta] \= (-1)^{\langle \alpha,\beta \rangle} [\alpha][\beta]
\ee
holds in $\CC_\Delta[\TT_-]_{<L}$, while~\eqref{eq:canlifthomo} together with the observation that $[S] = -1$  in the quotient by the winding ideal implies the same equation in the target, where the multiplication law is given by summation of cycles, see Proposition~\ref{prop:rngstructure}.
\par
The statement about decompositions into $\sigma$-eigenspaces and the corresponding morphisms are obvious.
\end{proof}

\subsection{Automorphisms of groupoid modules}

This section provides a slight generalization of the automorphisms given in~\eqref{intro:K}, with the goal of preparing the ground for similar discussion beyond the case $K=2$ of quadratic differentials. We fix a direction~$\theta$ and for any $\gamma \in \Gamma$ with $\arg(Z(\gamma)) = \theta$ we fix BPS invariants $\Omega(\gamma)$ such that $\Omega(\gamma) \neq 0$ for at most finitely many~$\gamma$.
Our goal is to check the well-definedness of automorphisms given by the law
\be \label{eq:Kgeneral}
\calK_\theta([\alpha]) \=  \prod_{\gamma \,:\, \arg(Z(\gamma))=\theta} (1-[\gamma])^{\Omega(\gamma) \, \langle \alpha, \gamma \rangle}[\alpha].
\ee
\par
\begin{lemma} \label{le:Kgeneralcont}
  Suppose that $\theta \in \Delta$ as usual. Then there is a unique automorphism $\calK_\theta: \wh{H}^\abs_\Delta \to \wh{H}^\abs_\Delta$, given on elements $[\alpha]$ by~\eqref{eq:Kgeneral} and such that the restriction to every subspace $\wh{H}^\abs_{\Delta^K}$ for every $K = \Delta_2^t$ a translated cone verifying~\eqref{eq:oppcone} and $t \in \CC$ arbitrary is continuous.
\end{lemma}
\par
\begin{proof} For each element $[\alpha] \in {H}^\abs_{\leq 1}(\wt{\Sigma}^\circ)$ the formula~\eqref{eq:Kgeneral} indeed defines an element in~$\wh{H}^\abs_\Delta$ by Corollary~\ref{cor:infiniteproduct} and thus defines by linearity $\calK_\theta$ on the direct sum submodule ${H}^\abs_{\leq 1}(\wt{\Sigma}^\circ)$.
\par
Fix $\Delta_2 \subset \CC$ verifying~\eqref{eq:oppcone} and $t \in \CC$. For any $L \in \RR$ the cone argument as in~\eqref{eq:conearg} implies that 
there is $L' = L'(L,t)$ such $\calK_\theta$ takes the submodule of ${H}^\abs_{\leq 1}(\wt{\Sigma}^\circ)_{\geq L'}$ generated by paths~$\alpha$ with support in $K=\Delta_2^t$ into  $\wh{H}_\Delta(\wt{\Sigma}^\circ)_{\geq L}$. Since the submodule of ${H}^\abs_{\leq 1}(\wt{\Sigma}^\circ)$ generated by paths~$\alpha$ with support in $K$ is dense in $\wh{H}^\abs_{\Delta^K}$, there is a unique continuous extension of~$\calK_\theta$ to $\wh{H}^\abs_{\Delta^K}$.
\par
For each generator in the winding submodule the image under $\calK_\theta$ again belongs to the winding submodule since the scalar product in the exponent of~\eqref{eq:Kgeneral} only depends on the $c$-image of the path and thus every generator of the winding ideals gets multiplied by an element $\wh{H}_\Delta^\abs(\wt{\Sigma}^\circ)$. Consequently, $\calK_\theta$ descends to a continuous endomorphism  of $\wh{H}^\abs_{\Delta^K}$.
\par
Moreover, the linearity of the scalar product in the exponent of~\eqref{eq:Kgeneral} implies the multiplicativity $\calK_\theta([\alpha]\cdot[\beta]) = \calK_\theta([\alpha]) \cdot \calK_\theta([\beta])$ on ${H}^\abs_{\Delta}(\wt{\Sigma}^\circ)$ and by continuity this extends to all of $\wh{H}^\abs_{\Delta^K}$ for all $K = \Delta_2^t$ as above. Since any two elements of $\wh{H}^\abs_{\Delta}$ to be multiplied belong to some common $\wh{H}^\abs_{\Delta^K}$, this suffices to to show that $\calK_\theta$ is an automorphism.
\end{proof}
\par
Key to this proof is the cone argument, which shows that only elements $\alpha$ with $|Z(\alpha)| \leq L'$ are be mapped to a non-zero element of the truncation $\wh{H}_\Delta(\wt{\Sigma}^\circ)_{\geq L}$ by adding a multiple of an element~$\gamma$ as in the definition of~$\calK$, i.e., by adding a element with central charge in the sector~$\Delta$. It is obviously true that adding such element~$\gamma$ with $|Z(\gamma)|>2L'$ will push any such~$\alpha$ outside the ball of radius~$L'>L$. This shows:
\par
\begin{corollary} \label{cor:calKLtruncation}
The restriction of $\calK_\theta$ to the truncation at level~$L'$ of
$\wh{H}^\abs_{\Delta^K}$ coincides with the map
\be
\calK_{\theta,L} ([\alpha]) \=  \prod_{\gamma \in \NN \wh{\gamma}_0}  (1-[\gamma])^{\Omega(\gamma) \langle \wh{\gamma}_0, c(\alpha) \rangle}\Big|_{L'} [\alpha]
\ee
where the restriction $|_{L'}$ is to be interpreted as cutting off the expansion of the prefactor (as in Lemma~\ref{lemm:description_of_elements}) at elements with central charge~$>L'$ in absolute value.
\end{corollary} 
\par
\begin{remark}\label{rema:action_on_twisted_torus}
By \cite{KS08,BridgelandRHDT}, the wall-crossing automorphisms for a ray $\ell\in \Delta$ is the exponential in $\wh{G}_\Delta$ of the series
\[ DT(\ell)=-\sum_{Z(\gamma) \in \ell}\Omega(\gamma)\sum_{n \geq 1}\frac{[n\gamma]}{n^2}, 
\]  which is an element of $\wh{\frakg}_\Delta$ by the support property. $\wh{G}_\Delta$ acts faithfully on $\CC_\Delta[[\TT_-]]$ by exponentiating the Poisson derivation $a \mapsto \{ a , \cdot \}$. Explicitly, $\exp([\alpha]) \in \wh{G}_\Delta$ acts by $[\beta] \mapsto \exp(\langle \alpha,\beta \rangle [\alpha]) [\beta]$. We recover the formula for $\calK$ in \eqref{eq:Kgeneral} by computing the action of $\exp(DT(\ell))$
\begin{align*}
  \exp(DT(\ell))[\beta] &\= \Bigl(-\sum_{Z(\gamma) \in \ell}\Omega(\gamma)\sum_{n \geq 1} \exp(\frac{\langle \gamma, \beta \rangle}{n}[\gamma]^n\Bigr)[\beta]\\
  &\= \Bigl(\prod_{Z(\gamma) \in \ell} \exp(-\sum_{n \geq 1} \frac{[\gamma]^n}{n})^{\Omega(\gamma)\langle \gamma, \beta \rangle}\Bigr)[\beta]\\
  &\= \Bigl(\prod_{Z(\gamma) \in \ell}(1-[\gamma])^{\Omega(\gamma)\langle \gamma, \beta \rangle}\Bigr)[\beta].
\end{align*}
\end{remark}
\par
\begin{remark} Besides the profinite approach to convergence taken here (and e.g.\ already in \cite{KS08}) there is also an analytic approach. The two are compared in the appendices of \cite{BridgelandRHDT}. Pursuing the latter in our context we need a norm on $H_{\leq 1}(Y)$ where paths with large absolute values of central charge are negligible, for example we might give $a = \sum \coeff(\alpha) [\alpha]$ the norm
\bes
||a||_R \= \sum_\alpha |\coeff(\alpha)| e^{-R|Z(\alpha)|}
\ees
as proposed in \cite[Footnote~20, p.~32]{GMN_spectral}. Consider the Banach space given by the completion with respect to~$||\cdot||_R$. Each path lifting (see Proposition~\ref{prop:defF}) is an element of this Banach space for $R$ large enough. The locally convex space given by their union over all~$R$ contains all detours. Moreover, there are natural 'tame' subspaces where the multiplication as in Proposition~\ref{prop:rngstructure} is well-defined. \emph{However the automorphism~$\calK_\theta$ from~\eqref{intro:K} is not continuous} in this norm in general, since for any element~$\beta$ with $e := \Omega(\gamma)\langle \gamma, \beta \rangle >0$ we find that that the operator norm of $\calK_\theta$ has lower bound
\bas
\frac{||\calK_\theta(n\beta)||_R}{n ||[\beta]||_R}
&\= \frac{|| (1+ [\gamma])^{en}[{n\beta}]||_R }{n ||\beta||_R} \\
&\= \Bigl| \sum_{k=0}^n \binom{n}{k} \exp\Bigl(-R \bigl( |Z(n\beta+ ke\gamma)| - |Z(n\beta)|\bigr)\Bigr)\Bigr| \\
&\,\geq\,  \bigl| \sum_{k=0}^n \binom{n}{k}  \exp(-Rke|Z(\gamma)|)\bigr| \= (1+ \exp(-R|Z(\gamma)|)^n
\eas
tending to infinity and $n \to \infty$. This appears to be problematic when trying to control the effect of~$\calK_\theta$ by approximations, a given e.g.\ in  Proposition~\ref{prop:approxFL} below.
  \end{remark}

\section{The path lifting rule $F$} \label{sec:pathliftF}

In this section we describe a path lifting rule from paths~$\wp$ in~$X^\circ$ with poles and zeros removed to paths in the circle bundle over $\Sigma^\circ$. The goal of the rule, as stated in Theorem~\ref{theo:homotopy_invariance}, is that the classes of the lifts in the completed groupoid of signed paths is well-defined under homotopies over the zeros, i.e.\ in $X^\star$. For this lift we give a version up to homotopy (denoted by~$\bfF$) and its shadow up to homology (denoted by~$F$). The idea for the lifting rule is entirely the one in \cite{GMN_spectral} and it depends on the choice of a direction~$\theta$: Beside the trivial lift one considers all possible detours following trajectories of the spectral network. After homotoping over a zero, the number of detours does not stay constant, but a pair of paths cancels after passage modulo the winding ideal (Lemma~\ref{lem:homot_over_zero}).
\par
The contribution here complementing \cite{GMN_spectral} is to justify the convergence of the lift in the completed groupoid, even if the spectral network is dense in subsurfaces.
\par
The detour rules and the convergence statements are given both for saddle-free directions, as point of departure, and then for directions with saddle connections. There, we define (as did \cite{GMN_spectral}) two lifts $\bfF^\pm$ that will later be shown to be the one-sided limits, and again prove convergence in the completed groupoid.
\par
We continue to fix an acute sector~$\Delta$ and assume throughout that $\theta \in \Delta$.

\subsection{Definition of the lift function} \label{sec:deflift}

We fix a surface with a quadratic differential $(X,q)$. 
The goal of this section is to formalize the lift function $F(\wp,\theta)$ in \cite[Section~4.4]{GMN_spectral} of a path~$\wp$ on~$X$  with respect to the direction~$\theta$ for arbitrary rank one directions~$\theta$, allowing spiral domains. Roughly, this lift is the sum of the lifts of~$\wp$ to spectral cover~$\Sigma^\circ$ and further using its tangent vectors to the unit tangent bundle~$\wt{\Sigma}^\circ$, allowing any number of elementary detours, defined as follows.
\par
First, assume that $(X,q,\theta)$ is saddle-free. Consider the spectral network~$W_\theta$ on~$\Sigma$ as defined above. Let $\wp_{1},\wp_2$ be the two lifts to~$\Sigma$ of the path~$\wp$ on~$X$, supposed to be disjoint from the zeroes of~$q$ and labeled in arbitrary order. We abbreviate $\wp_{1,2} = \wp_1 \cup \wp_2$. Let $W_\theta^-(\wp) \subset \wp_{1,2}$ the (possibly countable, possibly dense on subintervals) subset of points where $W_\theta^-$ intersects $\wp_{1,2}$. Note that if a point $y$ belongs to $W_\theta^-(\wp)$, then the other point $\sigma(y)$ in the same fiber over~$X$ does not belong to~$W_\theta^-(\wp)$ under our assumption of $(X,q,\theta)$ being saddle-free .
\par
An \emph{elementary detour} starts near a point~$y \in W_\theta^-(\wp)$ follows $W_\theta^-$ to the zero (with infinitesimally positive imaginary part), turns by the angle $2\pi$ around the zero and returns (still with infinitesimally positive imaginary part) along $W_\theta^+$  to $\sigma(y)$, in fact along the $\sigma$-image of the ray used for getting to the zero, compare \cite[Figure~8]{GMN_spectral} and the figures below, starting with Figure~\ref{fig:AmBritDetours}.
\par
A \emph{detour~$D_\bfy$ of~$\wp$} is determined by a sequence of points $\bfy = (y_1,\ldots,y_N)$ in $W_\theta^-(\wp)$ that alternate between the sheets~$\wp_1$ and~$\wp_2$ and whose images in~$\wp$ are monotone (i.e. equal to $\wp(t_i)$ for an increasing sequence of times $t_i$). Moreover we require that if the path~$\wp$ follows the spectral network along a subinterval~$\wp(I)$ and~$I=[i_0,i_1] \subset [0,1]$ is a maximal subinterval with this property, then only the left endpoint~$\wp(i_0)$ may occur as the start of a detour. That is, we prohibit repeating the same detour starting at several points of $\wp(I)$. The detour~$D_\bfy$ itself is given as the concatenation of the segment along~$\wp_{1,2}$ from the start to~$y_1$ (on the sheet of~$y_1$), the elementary detour at~$y_1$, the segment along~$\wp_{1,2}$ from~$\sigma(y_1)$ to~$y_2$ on the other sheet, the elementary detour at~$y_2$ etc until we reach~$y_n$, make an elementary detour to~$\sigma(x_n)$ and continue to the endpoint of either of~$\wp_{1,2}$. More precisely, at each start and end of the detour we replace the $90^\circ$ turns in the above definition by a smooth
approximation. Finally we let $\nl{D}_{\bfy}$ be the natural circle bundle lift of~$D_\bfy$ to~$\wt\Sigma$. 
\par
\begin{proposition} \label{prop:defF} Let $(X,q,\theta)$ be a saddle-free directed quadratic differential. Let $\wp$ be a smooth path in $X^{\ast}=X\setminus \calZ \cup \calP$. For any~$\bfy$ let~$\nl{D}_\bfy$ be the detours in direction~$\theta$. Then the path lifting 
\be
\bfF(\wp,\theta) \= \sum_{N \in \NN}\,\, \sum_{\bfy = (y_1,\ldots,y_N)} \dnl{D}_{\bfy}
\ee
is well-defined as an element of $\wh \Pi_\Delta(\wt{\Sigma}^\circ)$ and thus the passage to homology classes defines an element $F(\wp,\Delta) \in \wh H_\Delta(\wt{\Sigma}^\circ)$.
\end{proposition}
\par
Note that even if the path~$\wp$ is closed, all the lifts $\bfF(\wp,\theta)$ and $F(\wp,\theta)$ may depend on the starting point (and hence end point) of the path, as the set of possible detours may depend on that choice.
\par
Starting with the homotopy invariance in Section~\ref{sec:signed_hom_inv} we are only interested in the signed classes $\bfF(\wp,\theta) \in \wh \Pi_\Delta = \wh \Pi_\Delta^{\sign}(\wt{\Sigma}^\circ)$ and $F(\wp,\Delta) \in \wh H_\Delta =  \wh H_\Delta^{\sign}(\wt{\Sigma}^\circ)$ in the quotients by the winding ideal, which we denote by the same letters. If $\wp$ is closed, then by definition $\bfF(\wp,\theta) \in \wh{\Pi}^\abs_\Delta$ and $F(\wp,\theta) \in \wh{H}_\Delta^\abs$. 
\par
\begin{proof} For each elementary detour~$\gamma$ the central charge $Z(\gamma) \in e^{i\theta} \RR_+$. Since the sum of the path segments along~$\wp_{1,2}$ of the total detour have bounded integrals there is some~$t \in \CC$ such that every summand of $\bfF(\wp,\theta)$ has  $Z(\dnl{D}_{\bfy}) \in t + \Delta$. 
\par
It remains to show that for every~$L$ there are only finitely many summands
of $\bfF(\wp,\theta)$ with $|Z(\dnl{D}_{\bfy})| < L$. Since all the central charges of elementary detours belong to $e^{i\theta} \RR_+$ and add up under composition, and since $|Z(\dnl{D}_{\bfy})|$ is smaller than the sum of the detour periods at most by the period of~$\wp$, it suffices to show that there only finitely many elementary detours~$\gamma$ with $|Z(\gamma)| < L$. Since there are finitely many zeros it suffices to prove this for elementary detours that follow one of the three rays~$\frakr$ of~$W_\theta^-$ of one of the zeros, say~$z_1$.
\par
To describe an elementary detour consider the ray segment~$\frakr_0$ of~$\frakr$ from~$z_1$ to its first hit~$w_1$ with the union~$\wp_{1,2}$, then the ray segment~$\frakr_1$ from $w_1$ until the ray hits again~$\wp_{1,2}$ at~$w_2$, etc. Each of the ray segments~$\frakr_k$ is a first return for the directional flow (in direction~$-\theta$) to~$\wp_{1,2}$, starting at~$w_k$ and ending at~$w_{k+1}$.
Each elementary detour~$\gamma$ follows $\frakr_k,\frakr_{k-1},\ldots,\frakr_0$ in the inverse direction, then turns around the zero and follows those segments on the other sheet in the other direction. We conclude that
\be
Z(\gamma) \= 2 \Bigl(Z(\frakr_0) + \sum_{i=1}^k Z(\frakr_i)\Bigr) \quad \in e^{i\theta} \RR_+\,.
\ee
\par
We claim that there is a lower bound $|Z(\frakr_i)| > C=C(\wp,\theta) > 0$ for all $i \in \NN$. To prove this, observe that there a finite number of points (in fact at most~$3r$ where~$r$ is the number of zeros)  where the trajectory emanating from~$\wp_\pm$ hits a zero. (The reader may compare with the construction of suspensions over interval exchange transformations e.g.\ in \cite{Yoccoz} for a visualization of the idea, but we require the 'base interval' here neither to be a straight line nor connected.)
In the complement of those finitely many rays  the return 'time' is a locally constant function, hence constant on a finite number of intervals, possibly taking the value~$+\infty$ on some of them. In any case, $Z(\frakr_i)$ takes a finite number of values in $e^{i\theta} \RR_+$, thus proving the claim.
\end{proof}
\par
\medskip
In the next step we deal with the case of a direction~$\theta$ in which~$q$ has (one or more) saddle connections. In this case we define two functions~$\bfF^\pm(\wp,\theta)$ and $F^\pm(\wp,\theta)$. They will be shown in Section~\ref{sec:varyingtheta} to correspond to the one-sided limits of directions approaching $\theta$ from above or below. The geometric considerations in Section~\ref{sec:varyingtheta} together with \cite[Section~6.2]{GMN_spectral} will provide the motivation for the definition. Here we verify that these belong to the completed module of signed paths.
\par
In the presence of saddle connections the rays along $W^-_{\theta}$ may pass through other zeros of $\lambda$. When we approach a zero of $\lambda$ along $W^-_{\theta}$, we can turn by $\pi$ in either direction to follow another trajectory in $W^-_{\theta}$ opposite to its orientation. A ray does a \emph{$+$-turn}, if it turns counterclockwise by $\pi$ or equivalently it passes the zero parallel to a straight line with constant infinitesimally negative imaginary part with respect to the flat structure. A ray does a \emph{$-$-turn} at a zero if it turns clockwise by $\pi$ or equivalently if it passes the zero parallel to a straight line with constant infinitesimally positive imaginary part.
\par
Consider a path~$\frakr \subset \Sigma$ starting at some zero~$z_i$ following the negative spectral network $W^-_{\theta}$ opposite to its orientation and ending at a point in~$\wp_\pm$, hence in $W_\theta^-(\wp)$. We say that~$\frakr$ is a \emph{$+$-admissible (resp.\ minus-admissible) ray} if at every zero passed by $\frakr$ it does a $+$-turn (resp. $-$-turn).
\par
\begin{figure}
  \resizebox{6cm}{!}{
%
%


\begin{tikzpicture}
\tikzset{
  on each segment/.style={
    decorate,
    decoration={
      show path construction,
      moveto code={},
      lineto code={
        \path [#1]
        (\tikzinputsegmentfirst) -- (\tikzinputsegmentlast);
      },
      curveto code={
        \path [#1] (\tikzinputsegmentfirst)
        .. controls
        (\tikzinputsegmentsupporta) and (\tikzinputsegmentsupportb)
        ..
        (\tikzinputsegmentlast);
      },
      closepath code={
        \path [#1]
        (\tikzinputsegmentfirst) -- (\tikzinputsegmentlast);
      },
    },
  },
  mid arrow/.style={postaction={decorate,decoration={
        markings,
        mark=at position .5 with {\arrow[#1]{stealth}}
      }}},
}
  

\node[draw, circle, fill=black, inner sep=0.0pt] (A) at (0, 0) {}; 
\node[draw, circle, fill=black, inner sep=0.0pt] (B) at (2, 0) {};
\node (Cp) at (3, 1) {};
\node (C0) at (3, 0) {};
\node (Cm) at (3, -1) {};
\node (Dp) at (-1, 1) {};
\node (D0) at (-1, 0) {};
\node (Dm) at (-1, -1) {};
\node[color=red, draw, circle, fill=red, inner sep=0.0pt](P1) at (2.75, 0.63) {};
\node[color=red, draw, circle, fill=red, inner sep=0.0pt](P2) at (2, -0.1) {};
\node[color=red, draw, circle, fill=red, inner sep=0.0pt](P3) at (0.1, -0.1) {};
\node[color=dpurple, draw, circle, fill=dpurple, inner sep=0.0pt](M1) at (2.75, -0.94) {};
\node[color=dpurple, draw, circle, fill=dpurple, inner sep=0.0pt](M2) at (2, -0.2) {};
\node[color=dpurple, draw, circle, fill=dpurple, inner sep=0.0pt](M3) at (0.1, -0.2) {};

\draw[color=blue, postaction={on each segment={mid arrow=blue}}] (B) -- (A)  {};
\draw[color=blue, postaction={on each segment={mid arrow=blue}}] (Cp) -- (B)  {};
\draw[color=blue, postaction={on each segment={mid arrow=blue}}] (Cm) -- (B)  {};
\draw[color=blue, postaction={on each segment={mid arrow=blue}}] (A) -- (Dp)  {};
\draw[color=blue, postaction={on each segment={mid arrow=blue}}] (A) -- (Dm)  {};

\draw[color=red] (P1) -- (P2)  {};
\draw[color=red] (P2) -- (P3)  {};
\draw[color=red] (P3) -- (A)  {};

\draw[color=dpurple] (M1) -- (M2)  {};
\draw[color=dpurple] (M2) -- (M3)  {};
\draw[color=dpurple] (M3) -- (A)  {};

\draw[decorate, decoration={zigzag, segment length=4, amplitude=1.5}] (C0) -- (B); 
\draw[decorate, decoration={zigzag, segment length=4, amplitude=1.5}] (A) -- (D0); 

\node (Ab) at (0, 0) {$\times$};
\node (Bb) at (2, 0) {$\times$};

\end{tikzpicture}
    }
  \caption{The spectral network at a saddle connection between two distinct zeros. In red a $+$-admissible ray that is not minus-admissible. In purple a minus-admissible ray that is not $+$-admissible. (Branch cut as zig-zag line.)}
  \label{fig:plusminusadmissible}
\end{figure}
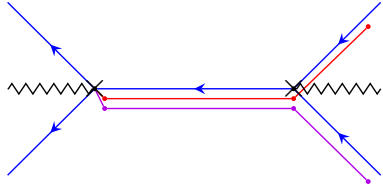
\par
An elementary \emph{elementary detour} for~$\theta^+$ starts at a point~$y \in W_\theta^-(\wp)$ follows some $+$-admissible ray~$\frakr$ to a zero (with infinitesimally positive imaginary part), passes the zero, turns by the angle $2\pi$ around the zero and returns along $W_\theta^+$  to $\sigma(y)$. Elementary detours for~$\theta^-$ are defined similarly using minus-admissible rays . (Both $\pm$-admissible detours start at points in $W_\theta^-(\wp)$ and end in $W_\theta^+(\wp)$.) Note that now the elementary detour is not determined by~$y$ alone: Consider for example a saddle connection at the boundary of a ring domain intersecting~$\wp_+$ in~$y$. The ray~$\frakr$ may loop any number of times along the saddle connection(s) bounding the cylinder before ending in~$y$. For convenience we often write $\frakr^y$ for a $\pm$-admissible ray ending in~$y$.
\par
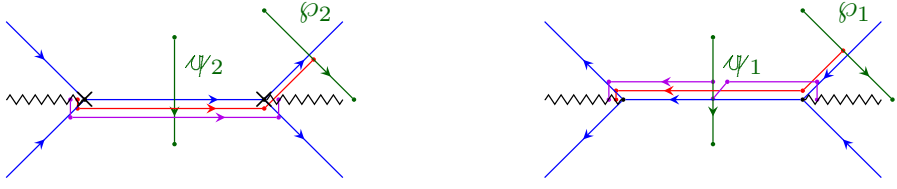
\begin{figure}
   \resizebox{12.5cm}{!}{
%
%


\begin{tikzpicture}
\tikzset{
  on each segment/.style={
    decorate,
    decoration={
      show path construction,
      moveto code={},
      lineto code={
        \path [#1]
        (\tikzinputsegmentfirst) -- (\tikzinputsegmentlast);
      },
      curveto code={
        \path [#1] (\tikzinputsegmentfirst)
        .. controls
        (\tikzinputsegmentsupporta) and (\tikzinputsegmentsupportb)
        ..
        (\tikzinputsegmentlast);
      },
      closepath code={
        \path [#1]
        (\tikzinputsegmentfirst) -- (\tikzinputsegmentlast);
      },
    },
  },
  mid arrow/.style={postaction={decorate,decoration={
        markings,
        mark=at position .5 with {\arrow[#1]{stealth}}
      }}},
  threequarter arrow/.style={postaction={decorate,decoration={
        markings,
        mark=at position .75 with {\arrow[#1]{stealth}}
      }}},
}
  

\node[draw, circle, fill=black, inner sep=0.0pt] (A) at (0, 0) {}; 
\node[draw, circle, fill=black, inner sep=0.0pt] (B) at (2, 0) {};
\node (Cp) at (3, 1) {};
\node (C0) at (3, 0) {};
\node (Cm) at (3, -1) {};
\node (Dp) at (-1, 1) {};
\node (D0) at (-1, 0) {};
\node (Dm) at (-1, -1) {};

\node[draw, circle, fill=black, inner sep=0.0pt] (AA) at (6, 0) {}; 
\node[draw, circle, fill=black, inner sep=0.0pt] (BB) at (8, 0) {};
\node (CCp) at (9, 1) {};
\node (CC0) at (9, 0) {};
\node (CCm) at (9, -1) {};
\node (DDp) at (5, 1) {};
\node (DD0) at (5, 0) {};
\node (DDm) at (5, -1) {};

\node[color=dgreen, draw, circle, fill=dgreen, inner sep=0.0pt](P1) at (8, 1) {};
\node[color=dgreen, draw, circle, fill=dgreen, inner sep=0.0pt](P2) at (9, 0) {};
\node[color=dgreen, draw, circle, fill=dgreen, inner sep=0.0pt](Ps1) at (2, 1) {};
\node[color=dgreen, draw, circle, fill=dgreen, inner sep=0.0pt](Ps2) at (3, 0) {};

\node[color=dgreen, draw, circle, fill=dgreen, inner sep=0.0pt](Q1) at (7, 0.7) {};
\node[color=dgreen, draw, circle, fill=dgreen, inner sep=0.0pt](Q2) at (7, -0.5) {};
\node[color=dgreen, draw, circle, fill=dgreen, inner sep=0.0pt](Qs1) at (1, 0.7) {};
\node[color=dgreen, draw, circle, fill=dgreen, inner sep=0.0pt](Qs2) at (1, -0.5) {};

\node[color=red, draw, circle, fill=red, inner sep=0.0pt] (PD1) at (8.45, 0.55) {};
\node[color=red, draw, circle, fill=red, inner sep=0.0pt] (PD2) at (8, 0.1) {};
\node[color=red, draw, circle, fill=red, inner sep=0.0pt] (PD3) at (5.92, 0.1) {};
\node[color=red, draw, circle, fill=red, inner sep=0.0pt] (PD4) at (5.92, 0.00) {};
\node[color=red, draw, circle, fill=red, inner sep=0.0pt] (PD5) at (-0.08, 0.00) {};
\node[color=red, draw, circle, fill=red, inner sep=0.0pt] (PD6) at (-0.08, -0.1) {};
\node[color=red, draw, circle, fill=red, inner sep=0.0pt] (PD7) at (2, -0.1) {};
\node[color=red, draw, circle, fill=red, inner sep=0.0pt] (PD8) at (2.55, 0.45) {};

\node[color=dpurple, draw, circle, fill=dpurple, inner sep=0.0pt] (QD1) at (7, 0.2) {};
\node[color=dpurple, draw, circle, fill=dpurple, inner sep=0.0pt] (QD2) at (5.84, 0.2) {};
\node[color=dpurple, draw, circle, fill=dpurple, inner sep=0.0pt] (QD3) at (5.84, 0.0) {};
\node[color=dpurple, draw, circle, fill=dpurple, inner sep=0.0pt] (QD4) at (-0.16, 0.0) {};
\node[color=dpurple, draw, circle, fill=dpurple, inner sep=0.0pt] (QD5) at (-0.16, -0.2) {};
\node[color=dpurple, draw, circle, fill=dpurple, inner sep=0.0pt] (QD6) at (2.16, -0.2) {};
\node[color=dpurple, draw, circle, fill=dpurple, inner sep=0.0pt] (QD7) at (2.16, 0.00) {};
\node[color=dpurple, draw, circle, fill=dpurple, inner sep=0.0pt] (QD8) at (8.16, 0.00) {};
\node[color=dpurple, draw, circle, fill=dpurple, inner sep=0.0pt] (QD9) at (8.16, 0.2) {};
\node[color=dpurple, draw, circle, fill=dpurple, inner sep=0.0pt] (QD10) at (7.16, 0.2) {};
\node[color=dpurple, draw, circle, fill=dpurple, inner sep=0.0pt] (QD11) at (7, 0.01) {};

%

\draw[color=blue, postaction={on each segment={threequarter arrow=blue}}] (A) -- (B)  {};
\draw[color=blue, postaction={on each segment={mid arrow=blue}}] (B) -- (Cp)  {};
\draw[color=blue, postaction={on each segment={mid arrow=blue}}] (B) -- (Cm)  {};
\draw[color=blue, postaction={on each segment={mid arrow=blue}}] (Dp) -- (A)  {};
\draw[color=blue, postaction={on each segment={mid arrow=blue}}] (Dm) -- (A)  {};

\draw[color=blue, postaction={on each segment={threequarter arrow=blue}}] (BB) -- (AA)  {};
\draw[color=blue, postaction={on each segment={threequarter arrow=blue}}] (CCp) -- (BB)  {};
\draw[color=blue, postaction={on each segment={mid arrow=blue}}] (CCm) -- (BB)  {};
\draw[color=blue, postaction={on each segment={mid arrow=blue}}] (AA) -- (DDp)  {};
\draw[color=blue, postaction={on each segment={mid arrow=blue}}] (AA) -- (DDm)  {};

\draw[color=dgreen, postaction={on each segment={threequarter arrow=dgreen}}] (P1) -- (P2)  node[near start, above right] {$\wp_1$};
\draw[color=dgreen, postaction={on each segment={threequarter arrow=dgreen}}] (Ps1) -- (Ps2)  node[near start, above right] {$\wp_2$};
\draw[color=dgreen, postaction={on each segment={threequarter arrow=dgreen}}] (Q1) -- (Q2)  node[near start, right] {$\wq_1$};
\draw[color=dgreen, postaction={on each segment={threequarter arrow=dgreen}}] (Qs1) -- (Qs2)  node[near start, right] {$\wq_2$};


%
\draw[color=red] (PD1) -- (PD2)  {};
\draw[color=red, postaction={on each segment={threequarter arrow=red}}] (PD2) -- (PD3)  {};
\draw[color=red] (PD3) -- (PD4)  {};
\draw[color=red] (PD5) -- (PD6)  {};
\draw[color=red, postaction={on each segment={threequarter arrow=red}}] (PD6) -- (PD7)  {};
\draw[color=red] (PD7) -- (PD8)  {};

\draw[color=dpurple, postaction={on each segment={mid arrow=dpurple}}] (QD1) -- (QD2)  {};
\draw[color=dpurple] (QD2) -- (QD3)  {};
\draw[color=dpurple] (QD4) -- (QD5)  {};
\draw[color=dpurple, postaction={on each segment={threequarter arrow=dpurple}}] (QD5) -- (QD6)  {};
\draw[color=dpurple] (QD6) -- (QD7)  {};
\draw[color=dpurple] (QD8) -- (QD9)  {};
\draw[color=dpurple] (QD9) -- (QD10)  {};
\draw[color=dpurple] (QD10) -- (QD11)  {};

%

\draw[decorate, decoration={zigzag, segment length=4, amplitude=1.5}] (C0) -- (B); 
\draw[decorate, decoration={zigzag, segment length=4, amplitude=1.5}] (A) -- (D0); 

\node (Ab) at (0, 0) {$\times$};
\node (Bb) at (2, 0) {$\times$};

\draw[decorate, decoration={zigzag, segment length=4, amplitude=1.5}] (CC0) -- (BB); 
\draw[decorate, decoration={zigzag, segment length=4, amplitude=1.5}] (AA) -- (DD0); 

\node (AAb) at (0, 0) {$\times$};
\node (BBb) at (2, 0) {$\times$};

\end{tikzpicture}
     }
\caption{Examples of special detours appearing in $\bfF(\wp,\theta^+)$: The red detour (of the path $\wp$) follows a $+$-admissible ray, which does not turn at the first zero it encounters. The purple detour (of the path $\wq$) takes two detours at the same point of the spectral network, following the American driving rule. (Superpose the two sheets for the mnemonic.)
}
  \label{fig:AmBritDetours}
\end{figure}
In the saddle-free case detour rules prohibited multiple detours at the same point or same subinterval of~$\wp \cap W_\theta^-$. In a saddle connection direction we slightly relax this condition, again in view of the result in Section~\ref{sec:varyingtheta}, and allow repeated detours from the same point, if this point is on a saddle connection. We use the following American and British driving rules mnemonics. Let $y\in W_\theta^-(\wp)$ be on a saddle connection and let $\frakr$ be a $+$-admissible ray ending at $y$, then $\frakr^{y}$ satisfies the \emph{American driving rule}, if the elementary detour along $\frakr^{y}$ turns right at $y$, i.e.\ the tuple of tangent vectors given by the initial tangent vector of the elementary detour and the terminal tangent vector of the trivial lift of $\wp$ at $y$ induces the orientation at $y$. In turn, a minus-admissible ray $\frakr^y$ satisfies the \emph{British driving rule} if the  detour turns left at $y$, i.e.\ the tuple of tangent vectors induces the opposite orientation.
\par
A \emph{detour} $D_{(\bfy,\bfR)}$  for~$\theta^+$ is determined by a sequence of $\bfy = (y_1,\ldots,y_N)$ in $W_\theta^-(\wp)$ together with a collection $\bfR = (\frakr^{y_1}, \ldots,\frakr^{y_N})$ of $+$-admissible rays ending at~$y_i$, such that whenever, $y_{i+1}=\sigma(y_i)$ the ray $\frakr^{y_i}$ satisfies the American driving rule. Note that whenever $\wp$ crosses a saddle connection at $x \in X$, then  $\pi^{-1}(x)=\{y,\sigma(y)\} \subset  W_\theta^-(\wp)$. Now $D_{(\bfy,\bfR)}$ is defined as above using the elementary detours along the rays~$\frakr^{y_i}$ and smoothing the turns. The definition of detours for~$\theta^-$ is analogous, using minus-admissible rays and the British driving rule. Note that the driving rules only apply to detours containing two elementary detours following one and the same saddle connection in opposite directions. The driving rule says that you should first follow the lane first crossed on an American resp. British street. 
\par
The same argument as in Proposition~\ref{prop:defF} shows that the detours $D_{(\bfy,R)}$ both for~$\theta^+$ and for~$\theta^-$ have support in $\Delta + t$ for some $t \in \CC$.
\par
\begin{proposition} \label{prop:Fpm} Let $\wp$ be a smooth path in $X^{\ast}$. The path liftings 
\bas
\bfF^+(\wp,\theta) \= \sum_{N \in \NN}\,\, \sum_{\bfy = (y_1,\ldots,y_N)} \sum_{\bfR \text{ $+$-admissible} \atop  \text{American driving rule}} \dnl{D}_{\bfy,\bfR} \\
\bfF^-(\wp,\theta) \= \sum_{N \in \NN}\,\, \sum_{\bfy = (y_1,\ldots,y_N)} \sum_{\bfR  \text{ minus-admissible,} \atop \text{British driving rule}} \dnl{D}_{\bfy,\bfR} \\
\eas
are well-defined as elements of $\wh\Pi_\Delta(\wt{\Sigma}^\circ)$ and thus define elements of the homological version $F^\pm(\wp,\theta) \in \wh{H}_\Delta(\wt{\Sigma}^\circ)$.
\par
\end{proposition}
\par
Again, subsequently we only care about the signed classes $\bfF^\pm(\wp,\theta) \in \wh \Pi_\Delta$ and $F^\pm(\wp,\Delta) \in \wh H_\Delta$ in the quotients by the winding ideal. Again, if $\wp$ is closed, then  $\bfF^\pm(\wp,\theta) \in \wh{\Pi}^\abs_\Delta$ and $F^\pm(\wp,\theta) \in \wh{H}_\Delta^\abs$.
\par
\begin{proof} Each of the rays~$\frakr^{y_i}$ starts with one of the $3r$ outgoing trajectories in direction~$\theta$. The proof of Proposition~\ref{prop:defF} works here as well replacing~$C$ with the minimum of return times to~$\wp_{1,2}$ and of the lengths of saddle connections.
\end{proof}
\par
The well-definedness statements in Proposition~\ref{prop:defF} and Proposition~\ref{prop:Fpm} together with the ring structure in Lemma~\ref{le:quotientring} and the construction of the path lifts as detours at all possible collection of points imply:
\par
\begin{corollary} \label{cor:composition}
The path lifting functions $\bfF(\wp,\theta)$ and $F(\wp,\theta)$ for~$\theta$ saddle free as well as $F^\pm(\wp,\theta)$ for directions~$\theta$ with saddle connections satisfy the composition rule of Theorem~\ref{thm:introF}.
  \end{corollary}

\subsection{Signed homotopy invariance} \label{sec:signed_hom_inv}
Let $(X,q,\theta)$ be a directed quadratic differential and $\theta \in \Delta$.
\par
\begin{theorem}\label{theo:homotopy_invariance}
  The path lifting function $\bfF$ has the following properties: 
  \begin{itemize}
    \item[i)] For smooth paths $\wp_1,\wp_2$ with $[\wt{\wp}_1]=[\wt{\wp}_2] \in \Pi_{\leq 1}^\sign(\wt{X}^\star)$ homotopic in the surface with only poles removed, the lifts agree in the signed homotopy groupoid, i.e.\ $\bfF(\wp_1,\theta)=\bfF(\wp_2,\theta) \in \wh{\Pi}_\Delta(\Sigma^\star)$ of the surface with only the poles removed. 
    \item[ii)] For composable paths $\wp_1,\wp_2$, the path lifting behaves like a homomorphism $\bfF(\wp_1,\theta)\cdot \bfF(\wp_2,\theta)=\bfF(\wp_2 \circ \wp_1,\theta) \in \wh{\Pi}_\Delta(\Sigma^\star)$. 
  \end{itemize}
As a consequence,  the path lifting function $F$ valued in the signed homology groupoid $\wh{H}_\Delta$ has these properties, too.
\end{theorem}
\par
This theorem, in particular part i) is the reason for dropping the surface argument in the final notation of the signed groupoids and writing just $\bfF(\wp,\theta) \in \wh{\Pi}_\Delta$ and $F(\wp,\theta) \in \wh{H}_\Delta$.
\par
A \emph{smooth homotopy} on $X$ is a homotopy of $C^1$-paths through $C^1$-curves. Clearly any such homotopy lifts uniquely to a homotopy on $\wt{X}$.
The following two lemmas are from \cite[Section~5.6]{GMN_spectral}. We briefly give the proofs, which do no require addressing convergence issues. 
\par
\begin{lemma}\label{lem:homot_over_zero}
  Let $(D,z)$ be a coordinate disc on $X$, such that $q=z d z^{\otimes 2}$. Assume the only critical trajectories intersecting the disc are the ones emanating from the zero at $z=0$. Let $x_1,x_2 \in D\setminus D \cap \pi(W_\theta)$ and $\wp,\wq$ be two smoothly homotopic paths in $D$ with endpoints $x_1,x_2$. Then $\bfF(\wp,\theta) =\bfF(\wq,\theta) \in \wh{\Pi}_\Delta$.
\end{lemma}
\begin{proof} If the $\wp,\wq$ are homotopic in $D\setminus\{0\}$ the statement is clear. Using this and the concatenation rule it is enough to prove the result for two paths $\wp,\wq$ as drawn in the central picture of Figure \ref{fig:homot_NS}. On the left of Figure \ref{fig:homot_NS} we see all five lifts of~$\wp$. The trivial lift from $x_1^-$ to $x_2^-$ is homotopic to the detour on the lower left that turns twice up to unit of winding. Hence these two paths cancel in $\wt{\Pi}_\Delta$. Now, it is easy to see that three remaining paths on to left are pairwise homotopic to a lift of~$\wq$ on the right. 
  \begin{figure}
   \input{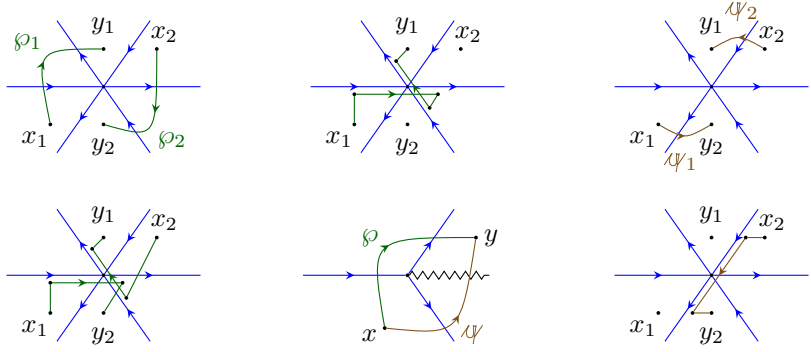}
    \caption{The lift $\bfF(\wp,\theta)$ for the two homotopic paths $\wp$ and $\wq$ around a zero: On top and left all five lifts of $\wp$ and on the right all three lifts of $\wq$. Note that $\wp_1$ and the detour on top in the middle cancels each other.}\label{fig:homot_NS}
  \end{figure}
\end{proof}
\begin{lemma}\label{lem:homot_over_traj}
  Let $(D,z)$ be a coordinate disc on $X$, such that $q=d z^{\otimes 2}$ and the only critical trajectory is the real axis $\Im(z)=0$. Let $x_1,x_2 \in D$ with $\Im(x_i)<0$ and $\wp,\wq$ be two smoothly homotopic paths in $D$ with endpoints $x_1,x_2$. Then $\bfF(\wp,\theta) =\bfF(\wq,\theta) \in \wt{\Pi}_\Delta$.
\end{lemma}
\begin{proof} By a similar argument as in the beginning of the proof of the previous lemma it is enough to show the following: The lift $\bfF(\wp,\theta)$ of a path $\wp$ with a single loop passing over the trajectory as shown on the right of Figure \ref{fig:homot_traj} is equal to its trivial lift. The two extra detours drawn in blue and pink on the left of Figure \ref{fig:homot_traj} cancel in the quotient by the winding ideal. Hence the assertion follows. 
  \begin{figure}
%
%


\begin{tikzpicture}
\tikzset{
  on each segment/.style={
    decorate,
    decoration={
      show path construction,
      moveto code={},
      lineto code={
        \path [#1]
        (\tikzinputsegmentfirst) -- (\tikzinputsegmentlast);
      },
      curveto code={
        \path [#1] (\tikzinputsegmentfirst)
        .. controls
        (\tikzinputsegmentsupporta) and (\tikzinputsegmentsupportb)
        ..
        (\tikzinputsegmentlast);
      },
      closepath code={
        \path [#1]
        (\tikzinputsegmentfirst) -- (\tikzinputsegmentlast);
      },
    },
  },
  mid arrow/.style={postaction={decorate,decoration={
        markings,
        mark=at position .5 with {\arrow[#1]{stealth}}
}}},
    ->-/.style={decoration={markings, mark=at position 0.6 with {\arrow{>}}}, postaction={decorate}},
}
  

\node (A) at (0, 0) {};
\node[draw, circle, fill=black, inner sep=0.0pt] (B) at (2, 0) {};
\node (C) at (4, 0) {};
\node (D1) at (1.3, 1) {};
\node (D2) at (2.7, 1) {};
\node (D3) at (1.3, -1) {};
\node (D4) at (2.7, -1) {};

\node[draw, circle, fill=black, inner sep=0.0pt] (Pl1) at (2.5, 0.5) {};
\node[draw, circle, fill=black, inner sep=0.0pt] (Pl2) at (3.5, 0.5) {};
\node[draw, circle, fill=black, inner sep=0.0pt] (Pl3) at (2.45, 0.05) {};

\draw[dgreen, ->-] (Pl1).. controls (2.8,-0.8) and (3.2,-0.8) ..  (Pl2);
    \node[dgreen] at (2.8,0.4) {$\wp_1$};

\node[draw, circle, fill=black, inner sep=0.0pt] (Ql1) at (1.5, -0.5) {};
\node[draw, circle, fill=black, inner sep=0.0pt] (Ql2) at (0.5, -0.5) {};
\node[draw, circle, fill=black, inner sep=0.0pt] (Ql3) at (0.47, 0.1) {};
\node[draw, circle, fill=black, inner sep=0.0pt] (Ql4) at (3.3, 0.1) {};
\node[draw, circle, fill=black, inner sep=0.0pt] (Ql5) at (1.5, 0.05) {};

\draw[dgreen, ->-] (Ql1).. controls (1.2,0.8) and (0.8,0.8) ..  (Ql2);
    \node[dgreen] at (0.8,-0.4) {$\wp_2$};
\draw[dpurple, ->-] (Ql1).. controls (1.25,0.88) and (0.75,0.88) ..  (Ql3);
\draw[dpurple, ->-] (Ql3) --  (Ql4);
\draw[dpurple, ->-] (Ql4) -- (Pl2);

\draw[dorange, ->-] (Ql1) -- (Ql5);
\draw[dorange, ->-] (Ql5) -- (Pl3);
\draw[dorange, ->-] (Pl3).. controls (2.8,-0.85) and (3.3,-0.85) ..  (Pl2);

\draw[color=blue, postaction={on each segment={mid arrow=blue}}] (A) -- (B)  {};
\draw[color=blue, postaction={on each segment={mid arrow=blue}}] (B) -- (C)  {};
\draw[color=blue, postaction={on each segment={mid arrow=blue}}] (B) -- (D1)  {};
\draw[color=blue, postaction={on each segment={mid arrow=blue}}] (D2) -- (B)  {};
\draw[color=blue, postaction={on each segment={mid arrow=blue}}] (B) -- (D3)  {};
\draw[color=blue, postaction={on each segment={mid arrow=blue}}] (D4) -- (B)  {};

  \draw[->, thick] (4.5,0) -- (5.5,0);

\node[draw, circle, fill=black, inner sep=0.0pt] (E) at (7, 0) {};
\node (E1) at (6.3, 1) {};
\node (E2) at (6.3, -1) {};
\node (F) at (9, 0) {};
\node (G) at (6, 0) {};
\node[draw, circle, fill=black, inner sep=0.0pt] (P1) at (7.5, 0.5) {};
\node[draw, circle, fill=black, inner sep=0.0pt] (P2) at (8.5, 0.5) {};

\draw[dgreen, ->-] (P1).. controls (7.8,-0.8) and (8.2,-0.8) ..  (P2);
    \node[dgreen] at (8.4,-0.6) {$\wp$};

\draw[color=blue, postaction={on each segment={mid arrow=blue}}] (E) -- (F)  {};
\draw[color=blue, postaction={on each segment={mid arrow=blue}}] (E1) -- (E)  {};
\draw[color=blue, postaction={on each segment={mid arrow=blue}}] (E2) -- (E)  {};
\draw[decorate, decoration={zigzag, segment length=4, amplitude=1.5}] (E) -- (G);

\end{tikzpicture}
    \caption{The lift $\bfF(\wp,\theta)$ for a path $\wp$ crossing a single trajectory twice. Note that the orange and the purple detour differ up to smooth homotopy only by their winding number.}\label{fig:homot_traj}
  \end{figure}
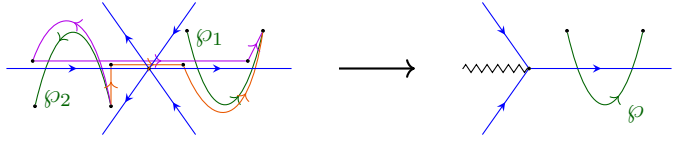
\end{proof}
\par
The main technical tool in several proofs is the \emph{truncated spectral network} $W_{\theta,L}$ depending on a direction~$\theta$ and a cut off parameter~$L \in R$. It is the subset of the spectral network~$W_\theta$ consisting of all ray segments such that the distance to some zero is at most~$L$. Obviously, when homotoping a path over an endpoint of truncated spectral network, the detours differ, but in a very controlled way:
\par
\begin{lemma}\label{lem:cross_end_of_traj}
  Consider the truncated spectral network $W_{\theta,L}$ and assume that we chose a coordinate neighborhood $(D,z)$ on $X$ such that the only truncated trajectory intersecting $D$ is given by the negative real line, in other words
  \[ \pi(W_{\theta,L}) \cap D=\{ z \in D \mid \Re(z)\leq 0, \Im(z)=0\}.
  \] Let $\frakr$ on $\Sigma$ be the corresponding truncated ray. Let $\wp_{\leq 0}$ and $\wp_{\geq 0}$ be two smoothly homotopic paths in $D^\times$ from $-i$ to $i$ through the half-planes $\Re(z)\leq 0$ respectively $\Re(z) \geq 0$. Then $\bfF(\wp_{\geq 0},\theta) - \bfF(\wp_{\leq 0},\theta)=\nl{D}_\frakr$ 
\end{lemma}
\par
\begin{proof}[Proof of Theorem \ref{theo:homotopy_invariance}] For i) consider a smooth homotopy $H: [0,1]^2 \to X, (s,t) \mapsto \wp_s(t)$ and let~$B$ be the supremum of the flat lengths of all paths $\wp_s$ for $s \in [0,1]$. To prove the statement it suffices that the two lifts agree in all truncations $\Pi_{\leq 1}(\Sigma^\star)_{<L}$ for any~$L$ up to elements of the winding ideal. Fix $L' = L +B$. Then the $<L$-truncations of the lift depend by Lemma~\ref{lem:cross_end_of_traj} only on detours at points in $W_{\theta,L'}$. This truncated network is compact. Since the image of~$H$ in $\Sigma^\star$ is compact, too, we may cover~$H$ by finitely many open sets, such that in each open set a crossing with the truncated spectral network is one of the two situation in Lemma~\ref{lem:homot_over_zero} or Lemma~\ref{lem:homot_over_traj}. Those lemmas give equality up to the winding ideal in both cases.
\par
The statement ii) is obvious from the rule for concatenation of detours.
\end{proof}

\section{Extending the path-lifting to $\calA_0$-laminations}\label{sec:A0laminations}

In Section~\ref{sec:deflift} we defined the lifts~$\bfF(\wp,\theta)$ of a path $\wp$ with endpoints $x_1,x_2 \in X$. Those lifts are in general not dense in $\wh{H}_\Delta^\abs$. For example in the $A_n$-quiver case, the genus zero strata $\calQ(1^{n+1},-n-5)$ with one higher order pole only, the surface $X^\star$ is simply connected and thus by Theorem~\ref{theo:homotopy_invariance} all lifts $F(\wp,\theta)$ are trivial. For this reason we follow the proposal outlined in \cite{GMNframed}, see Section~9 and Section~10 for examples in the $A_n$-quiver case, $n=2,3$. We extend in this section the definition of the lift functions~$\bfF$ (and also the homology groupoid version) to a more general class of paths with endpoints on the punctures. To encode where exactly the paths end we mark points on the real blowup of the surface with quadratic differential at the higher order poles. This will be encoded in the notion of ciliated surface and $\calA_0$-laminations.
\par
We continue with the standing assumption that $\theta \in \Delta$, a fixed acute sector.

\subsection{Ciliated surface associated to a quadratic differentials}
\label{sec:ciliated}

A \emph{ciliated surface} $\XX$ is compact oriented surface with boundary together with a finite collection of marked points. A marked point on the boundary is called \emph{cilium} and marked point in the interior is called \emph{hole}. A triangulation of a ciliated surface is a triangulation with vertices given by a subset of the marked points.
\par
Given a meromorphic quadratic differential $(X,q,\theta)$ of infinite area one can associate a ciliated surface~$\XX$ as follows. First we mark all poles of order exactly~$2$. These will be interior points, i.e.\ holes. We denote by~$\calH$ the set of holes in~$\XX$. Then we perform a real oriented blow-up at all poles of order $k>2$ and mark $k-2$ points on the resulting boundary. This yields a ciliated surface associated to the quadratic differential. Performing a real oriented blowup at all the poles of the abelian differential~$\lambda$, the spectral cover~$\pi: \Sigma \to X$ extends to a spectral cover $\pi: \Sigma \to \XX$ of ciliated surfaces, abusively denoted by the same letter. For each boundary component of~$\XX$ with~$k-2$ marked points and $k$ even, there are two boundary components on~$\Sigma$, each with $k-2$ marked points. For each boundary component of~$\XX$ with $k$ odd, there is one boundary component with $2(k-2)$ marked points on~$\Sigma$.
\par
If $\theta$ is a saddle-free direction on $(X,q)$ we want to be more precise with the location of the cilia: Recall that choosing one representative of a trajectory in each flat strip gives (together with the boundary curves) a triangulation of~$\XX$. In this case we mark the end points of the horizontal trajectories and refer to them as cilia. We also need to mark the mid points between any two cilia on the boundary curves. We refer to them as \emph{orange marks}, following the pictures in \cite{GMNframed}. Every ciliated surface associated with a saddle-free quadratic differential has thus $k$ orange marks at the boundary component associated with each pole of order $k>2$.
\par

\subsection{$\calA_0$-laminations and coordinates}

We need two versions of laminations that differ by a minor ``shift'' of the conventions at the boundary. The first version is taken from \cite{FG_dual_Teichmueller}. It can be prescribed by coordinates on the intersection with a triangulation. The second version is convenient for the path lifting. On an abstract ciliated surface both versions coincide:
\par
\begin{definition} \label{def:A0lam}
Given a ciliated surface $\XX$ a \emph{($\calA_0$)-lamination}~$\calL$ is an isotopy class
class of a finite collection of disjoint, non-selfintersecting smooth curves~$\gamma_i$ with integer weights~$w(\gamma_i)$, such that
  \begin{itemize}
  \item[i)] all curves are either closed or connect two points on the boundary disjoint from the cilia,
  \item[ii)]  all weights are positive unless the curve goes exactly around a cilium or is contractible to a hole,
  \item[iii)] and for a hole or a boundary segment between two cilia the sum of weights of the curves ending there is $0$.  
  \end{itemize}
  up to the following relations:
 \begin{itemize}
  \item[iv)] A contractible curve with any weight or a curve with weight $0$ can be removed.
  \item[v)] Two homotopic curves can be replaced by one upon adding their weights.
 \end{itemize}  
\end{definition}
(The qualifier $\calA_0$ is used in \cite{FG_dual_Teichmueller} to distinguish them from $\calX$-laminations. $\calX$-laminations do not appear in this paper and thus we usually omit the $\calA_0$.)
\par
Denote by $\calL(\XX)$ the set of laminations on the ciliated surface $\XX$. Given a triangulation~$\TT$ of the ciliated surface we can obtain coordinates on $\calL(\XX)$ as follows:
\par
Given a lamination $\calL \in \calL(\XX)$, homotope the curves so that each curve intersects the triangulation in the minimal number of points. Then we associate to each edge~$e$ of the triangulation one half times the weighted number of intersection points of~$e$ and~$\calL$, each intersection number being weighted with the weight of the path in~$\calL$. By assumption all coordinates associated to edges homotopic to boundary segments are~$0$. This gives a map $\iota: \calL(\XX) \to \ZZ^{E^i}$, where $E^i$ is the set of \emph{interior edges}, i.e., edges that are not homotopic to a boundary segment. (Integrality of the image follows from iii) in the definition). This map is an isomorphism: 
\par
\begin{proposition} \label{prop:FGalgo}
Suppose that the ciliated surface~$\XX$ is provided with a triangulation~$\TT$ with vertex set contained in the set of cilia. Then for any tuple $(n_e) \in \ZZ^{E^i}$ there is a lamination~$\calL$ such that $\iota(\calL) = 
(n_e)_{e \in E^i}$. 
\end{proposition}
\par
The proof is algorithmic, taken from \cite[Section~3.2]{FG_dual_Teichmueller}, see Figure~7 there: For every $s \in \NN$ sufficiently large the shifted numbers $\wt{n_e} = n_e + s$ satisfy the inequalities
\be
|\wt{n_a} - \wt{n_b}| \leq \wt{n_c} \leq \wt{n_a } + \wt{n_b}
\ee
for every triangle with sides $a,b,c \in \TT$. For such a choice of~$s$ mark $2\wt{n_e}$ points on each edge of~$\TT$. Inside each triangle $abc$, connect the points on the boundary edges in the unique way by $\wt{n_a } + \wt{n_b} + \wt{n_c}$ non-intersecting and non-self-intersecting arcs. Finally add around each cilium and around each hole a close curve with weight~$-s$. The class of the resulting lamination does not depend on the choice of~$s$.
\par
Suppose that the ciliated surface~$\XX$ is associated with a directed quadratic differential $(X,q,\theta)$. Then Definition~\ref{def:A0lam} gives a notion of a lamination using the cilia as defined in Section~\ref{sec:ciliated}. Alternatively, we get a notion of \emph{shifted ($\calA_0$)-lamination} by replacing in 
Definition~\ref{def:A0lam} everywhere 'cilium' by 'orange mark'. Laminations and shifted laminations are obviously in bijection in two ways depending on the orientation of the shift. The bijection consistent with path lifting rules below takes any lamination and performs in a tubular neighborhood of every boundary component (viewed locally as the outer boundary circle) and performs a $1/2k$-th Dehn twist 'clockwise' (see Figure~\ref{cap:FG1} for an illustration of our convention), taking each cilium into the adjacent orange mark.

\subsection{The path lifting rule for laminations}

We fix a surface with a quadratic differential $(X,q)$ of infinite area. The goal of this section is to formalize the lift function $F(\calL,\theta)$ of a shifted lamination on the ciliated surface~$\XX$  associated with the direction~$\theta$. The lift of a lamination uses the same procedure as in Section~\ref{sec:deflift}, take all possible elementary detours, but with a restricted behavior at the boundary curves due to a preferred-sheet rule. The lifting rule appears also in \cite{GMNframed}, in particular to the examples in Section~10.
\par
For the definition of $F(\calL,\theta)$ we restrict to \emph{shifted laminations without closed loops}. Those will suffice for our purposes if~$\XX$ has no holes and we avoid complications due to the lack of an orientation of the closed curves in the laminations. (See Section~\ref{sec:withholes} for the general case.) We also may use the equivalence relations~iv) and~v) to take are representative of $\calL$ such that the weight of any curve is in $\{\pm 1\}$. In the absence of closed curves the orientation will be restored by the preferred-sheet rule as follows.
\par
We fix as usual $(X,q,\theta)$ and an abelian differential~$\lambda$ on~$\Sigma$. For each boundary component of~$\Sigma$, each segment in between two orange marks is labeled \emph{incoming} if one (hence all) of the $\lambda$-trajectories in direction~$\theta$ ending at the segment are incoming (i.e.\ positive real axis). It is labeled \emph{outgoing} if the trajectories ending at the segment are outgoing, i.e.\ negative real axis.
\par
A \emph{loop decomposition $\{\wp_1,\ldots,\wp_n\}$} of a lamination~$\calL$ is a set of (unoriented) closed paths~$\wp_i$, whose support is contained in the support of the lamination together with the boundary segments between the orange marks and such that for each curve~$\gamma$ of~$\calL$ the number of loops~$\wp_i$ passing through~$\gamma$ equals~$|w(\gamma)|$. Such decompositions always exists because of the condition iii).
\par
We now fix a loop decomposition $\{\wp_1,\ldots,\wp_n\}$ of~$\calL$ and assume first that~$(X,q,\theta)$ is saddle-free. A \emph{detour~$D_{\bfy_j}$ of~$\wp_j$} is determined by the choice of an orientation of $\wp_j$ and a sequence of points $\bfy_j = (y_1,\ldots,y_N)$ in $W_\theta^-(\wp)$, whose images in~$\wp$ are monotone, such that the $y_i$ alternate between the sheets for every path segment of~$\wp_j$ in the interior of~$\XX$. We still prohibit multiple detours starting at any given subinterval where $\wp_j$ and $W_\theta^-(\wp)$ intersect. As above, the detour~$D_{\bfy_j}$ is the path concatenating a lift of a segment of~$\wp_j$ from the start to~$y_1$, the elementary detour at~$y_1$, a lift of a segment of~$\wp_j$ from~$\sigma(y_1)$ to~$y_2$ etc. 
\par
However, we only allow detours~$D_{\bfy_j}$ that obey the \emph{preferred-sheet rule}:
\begin{itemize}
\item At each point where $D_{\bfy_j}$ leaves a boundary segment into a curve in~$\calL$ with positive weight, the boundary segment has to be incoming.
\item At each point where $D_{\bfy_j}$ enters into a boundary segment from a curve in~$\calL$ with positive weight, the boundary segment has to be outgoing.
\item For curves in~$\calL$ with negative weight the above rules apply with the role of incoming and outgoing reversed.
\end{itemize}
\par
As above we turn detours into smooth paths by smoothing the corners when entering an elementary detour and entering or exiting a boundary segment. We can moreover isotope the segments along the boundary into the interior and consider them as smooth paths on~$\Sigma$ (without real oriented blowup, thus justifying the abuse of notation.) Finally we define $\dnl{D}_{\bfy_j}$ to be the natural lift of~$D_{\bfy_j}$.
We can now state the analog of Proposition~\ref{prop:defF}.
\par
\begin{proposition} \label{prop:defFlam}
  Let $(X,q,\theta)$ be a saddle-free directed quadratic differential of infinite area. Let $\calL$ be a lamination without closed loops and with support in the smooth part $\XX^{\circ}= \XX \setminus \calZ \cup \calP$. Then the path lifting 
\be
\bfF(\calL,\theta) \= \sum_{\calL = \{\wp_1,\ldots,\wp_n\}} \sum_{j=1}^n \sum_{N \in \NN}\,\, \sum_{\bfy_j = (y_1,\ldots,y_N)} \dnl{D}_{\bfy_j}\,,
\ee
where the first sum is over all loop decompositions of the lamination and the last sum is over all detours of~$\wp_j$ obeying the preferred-sheet rule, is well-defined as an element of $\wh\Pi_\Delta(\wt{\Sigma}^\circ)$ and thus defines an element $F(\wp,\theta) \in \wh{H}_\Delta(\wt{\Sigma}^\circ)$. In fact, by construction $F(\wp,\theta) \in \wh{H}_\Delta^\abs$.
\end{proposition}
\par
\begin{proof} Choose a representation of $\calL$ by arcs on $\XX \setminus \calZ \cup \calP$. Since~$\calL$ has finitely many segments with finite weights, there are at most finitely many path decompositions. For each of them the convergence proof as in Proposition~\ref{prop:defF} applies, since an infinite number of detours may occur only in spiral domains, which form a compact subset of~$\XX$ disjoint from the boundary.
\end{proof}
\par
We now proceed to a direction~$\theta$ in which~$q$ has (one or more) saddle connections. We simply combine the path decomposition of laminations with the $\pm$-admissible detour rules from Section~\ref{sec:deflift} and obtain with the same proof as above:
\par
\begin{proposition} \label{prop:FpmLam}
Let $\calL$ be a lamination without closed loops and with support in the smooth part $\XX^{\circ}$. The path liftings 
\bas
\bfF^+(\calL,\theta) \=  \sum_{\calL = \{\wp_1,\ldots,\wp_n\}} \sum_{j=1}^n  \sum_{N \in \NN}\,\, \sum_{\bfy = (y_1,\ldots,y_N), \bfR \atop \text{$+$-admissible}} \dnl{D}_{\bfy,\bfR} \\
\bfF^-(\calL,\theta) \=  \sum_{\calL = \{\wp_1,\ldots,\wp_n\}} \sum_{j=1}^n  \sum_{N \in \NN}\,\, \sum_{\bfy = (y_1,\ldots,y_N), \bfR \atop \text{$-$-admissible}} \dnl{D}_{\bfy,\bfR}\,, \\
\eas
where the first sum is over all loop decompositions of the lamination and the last sum is over all detours of~$\wp_j$ obeying the preferred-sheet rule, are well-defined as elements of $\wh\Pi_\Delta(\wt{\Sigma}^\circ)$ and thus define elements $F^\pm(\wp,\theta) \in \wh{H}_\Delta^\abs \subset \wh{H}_\Delta(\wt{\Sigma}^\circ)$.
\end{proposition}
\par
Finally we extend the homotopy-invariance in Theorem~\ref{theo:homotopy_invariance} to laminations.
\par
\begin{proposition} \label{prop:boundaryhomotopy}
The lift $\bfF(\calL,\theta)$ is independent of the representations of $\calL$ by a collection of arcs up to homotopy in $\XX \setminus \calP$.
\end{proposition}
\par
  \begin{figure}


\begin{tikzpicture}
\tikzset{
  on each segment/.style={
    decorate,
    decoration={
      show path construction,
      moveto code={},
      lineto code={
        \path [#1]
        (\tikzinputsegmentfirst) -- (\tikzinputsegmentlast);
      },
      curveto code={
        \path [#1] (\tikzinputsegmentfirst)
        .. controls
        (\tikzinputsegmentsupporta) and (\tikzinputsegmentsupportb)
        ..
        (\tikzinputsegmentlast);
      },
      closepath code={
        \path [#1]
        (\tikzinputsegmentfirst) -- (\tikzinputsegmentlast);
      },
    },
  },
  mid arrow/.style={postaction={decorate,decoration={
        markings,
        mark=at position .5 with {\arrow[#1]{stealth}}
      }}},
  other arrow/.style={postaction={decorate,decoration={
        markings,
        mark=at position .75 with {\arrow[#1]{stealth}}
      }}},
}

\node[draw, rectangle, fill=orange, inner sep=1pt] (OrA) at (-0.7, 2.41) {};
\node[draw, rectangle, fill=orange, inner sep=1pt] (OrA1) at (1.58, 2.69) {};
\node[draw, circle, fill=black, inner sep=1pt] (A) at (0, 2.5) {};
\node[inner sep=0.0pt] (B) at (3, 1) {};
\node[inner sep=1pt] (D) at (0.9, 0.2) {};

\node[draw, circle, fill=blue, inner sep=0.5pt] (Z1) at (1.4, 0.6) {};
\node[draw, circle, fill=blue, inner sep=0.5pt] (Z3) at (0.3, 0.5) {};
\node[draw, circle, fill=blue, inner sep=0.5pt] (Z4) at (2, 2) {};

\node[draw,circle,color=pgreen,fill=pgreen,inner sep=0pt] (G1) at (-0.4, 2.45) {};
\node[draw,circle,color=pbrown,fill=pbrown,inner sep=0.0pt] (G2) at (0.5, 2.56) {};
\node[inner sep=0.0pt] (G3) at (-0.45,1){};
\node[inner sep=0.0pt] (G4) at (-0.35,0.9){};
\node[inner sep=0.0pt] (G5) at (-0.4,0.5){};
\node[inner sep=0.0pt] (G6) at (-0.45,0.5){};

\draw (G2) node[above left] {out};

\draw (A) -- (B);
\draw (D) -- (A);

\draw[color=blue, postaction={on each segment={mid arrow=blue}}] (A) -- (Z1);
\draw[color=blue, postaction={on each segment={mid arrow=blue}}] (A) -- (Z3);
\draw[color=blue, postaction={on each segment={mid arrow=blue}}] (A) -- (Z4);

 \draw plot [smooth, tension=1.1] coordinates {(-1.5, 2.3) (A) (1.7, 2.7)};

 \draw[color=pgreen] plot [smooth, tension=1.2] coordinates {(G1)(G3)(G6)};
 \draw[color=pbrown] plot [smooth, tension=0.5] coordinates {(G2)(G4)(G6)};
\draw (G3) node[below right,color=pbrown] {$\wq$};
\draw (G4) node[left,color=pgreen] {$\wp$};
\draw node[circle,color=pgreen,fill=pgreen,inner sep=0.5pt] at (G1){};
\draw node[circle,color=pbrown,fill=pbrown,inner sep=0.5pt] at (G2){};

\end{tikzpicture}

    \caption{Homotoping an arc along the boundary: The two arcs of a lamination $\wp$ and $\wq$ contribute homotopic detours satisfying the preferred sheet rule.}\label{fig_homot_boundary}
  \end{figure}
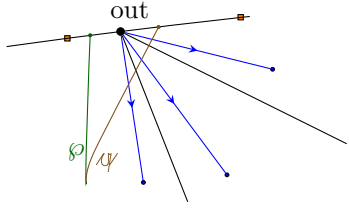
\par
\begin{proof}
  It follows from Theorem~\ref{theo:homotopy_invariance} that $\bfF(\calL,\theta)$ is independent of the collection of arcs up to homotopy in the $\XX \setminus \calZ \cup \calP$  with fixed end points on the boundary. So we are left with showing that the lift does not change when moving the end points along the boundary as long as we are not crossing an orange mark. Here the preferred-sheet rule comes into play.
  \par
  We consider the two preimages $b_i,b_o$ in $\Sigma$ of a fixed boundary component~$b$ of~$\XX$. Here $b_i$ supposed to be incoming and $b_o$ outgoing. Let $\wp$ and $\wq$ be two arcs oriented towards $b$ that differ by a homotopy relative to the boundary in a tubular neighborhood of $b$ containing no critical trajectories asymptotic to boundary components different from $b$. Hence, $\wp$ and $\wq$ differ by crossing several critical trajectories $t_1,\dots,t_k$  asymptotic to $b$ as depicted in Figure \ref{fig_homot_boundary}. Given a detour $D_\wp$ of the arc $\wp$ (i.e.\ a collection of points $\bfy$ in $W_\theta(\wp)$) it can either satisfy the preferred-sheet rule or not. If it does, it is ends on the outgoing boundary component $b_o$. Clearly, there is a $D_{\wq}$ satisfying the preferred-sheet rule that is homotopic to $D_\wp$. However, there could be extra detours for $\wq$ by following one of the trajectories $t_1,\dots,t_k$. This is only possible on the sheet bounded by $b_o$ by the detour rule. After doing such detour we end up on the boundary component $b_i$ and hence none of these extra detours will satisfy the preferred sheet rule. Hence, there are no extra detours contributing to the lift of $\wq$ in comparison with $\wp$. The case of $\wp$ and $\wq$ differing by a homotopy along their initial boundary segment is treated by a similar argument. 
\end{proof}

\subsection{Approximation by lifts of laminations} \label{sec:approx}

Suppose throughout this subsection that $(X,q,\theta)$ is a saddle-free directed quadratic differential of infinite area. For expository purposes we moreover assume for the moment that~$q$ has no double poles and return to the general case later in this subsection. 
\par
The goal of this section is to show that any monomial in the standard generators of the charge lattice~$\Gamma$ can be approximated in the $Z$-adic topology by rational multiples of monomials in $F(\calL_i,\theta)$ for appropriately chosen laminations~$\calL_i$. In the next subsection we show a slightly weaker statement in the presence of double poles.
\par
Recall that (core curves of the) the strips of~$(X,q,\theta)$ define interior edges of a triangulation~$\TT$. For each such strip (or the corresponding edge)~$e$ we let~$\gamma_e$ be the unique 'crossing' saddle connection joining the zeros at the two boundaries of the strip and we let $\wh{\gamma}_e$ be the hat-homology class of the lift, oriented such that $\Im(\wh{\gamma_e})>0$. Finally we let $\dnl\gamma_e$ be the canonical lift of $\wh\gamma_e$. Observe that in the proof of the following lemma we cut and reglue path segments, so that the statement is only a homological version.
\par
\begin{figure}
\begin{minipage}{0.45\textwidth}


\begin{tikzpicture}
\tikzset{
  on each segment/.style={
    decorate,
    decoration={
      show path construction,
      moveto code={},
      lineto code={
        \path [#1]
        (\tikzinputsegmentfirst) -- (\tikzinputsegmentlast);
      },
      curveto code={
        \path [#1] (\tikzinputsegmentfirst)
        .. controls
        (\tikzinputsegmentsupporta) and (\tikzinputsegmentsupportb)
        ..
        (\tikzinputsegmentlast);
      },
      closepath code={
        \path [#1]
        (\tikzinputsegmentfirst) -- (\tikzinputsegmentlast);
      },
    },
  },
  mid arrow/.style={postaction={decorate,decoration={
        markings,
        mark=at position .5 with {\arrow[#1]{stealth}}
      }}},
  other arrow/.style={postaction={decorate,decoration={
        markings,
        mark=at position .75 with {\arrow[#1]{stealth}}
      }}},
}

\node[draw, rectangle, fill=orange, inner sep=1pt] (OrA) at (-0.58, 2.39) {};
\node[draw, circle, fill=black, inner sep=1pt] (A) at (0, 2.5) {};
\node[inner sep=0.0pt] (A1) at (-1, 1.35) {};
\node[draw, circle, fill=black, inner sep=1pt] (B) at (3, 1.5) {};
\node[draw, rectangle, fill=orange, inner sep=1pt] (OrB) at (2.95, 1.75) {};
\node[draw, circle, fill=black, inner sep=1pt] (C) at (2.5, -1) {};
\node[draw, rectangle, fill=orange, inner sep=1pt] (OrC) at (2.85, -0.89) {};
\node[inner sep=0.0pt] (C1) at (3.7, -0.05) {};
\node[inner sep=0.0pt] (C2) at (3.3, 0.15) {};
\node[draw, circle, fill=black, inner sep=1pt] (D) at (-0.5, 0) {};
\node[draw, rectangle, fill=orange, inner sep=1pt] (OrD) at (-0.35, -0.5) {};

\node[draw, circle, fill=blue, inner sep=0.5pt] (Z1) at (1, 1.2) {};
\node[draw, circle, fill=blue, inner sep=0.5pt] (Z2) at (1.5, 0.3) {};
\node[draw, circle, fill=blue, inner sep=0.5pt] (Z3) at (-0.8, 1.8) {};
\node[draw, circle, fill=blue, inner sep=0.5pt] (Z4) at (-0.6, 1.4) {};

\node[draw, circle, fill=blue, inner sep=0.5pt] (Z5) at (3.0, 0) {};
\node[draw, circle, fill=blue, inner sep=0.5pt] (Z7) at (3.4, -0.1) {};

\node[draw, circle, color=pgreen, fill=pgreen, inner sep=0.0pt] (G1) at (-0.9, 2.34) {};
\node[draw, circle, color=pgreen, fill=pgreen, inner sep=0.0pt] (G2) at (2.06, -1.3) {};
\node[draw, circle, color=pgreen, fill=pgreen, inner sep=0.0pt] (G3) at (3.3, -0.81) {};
\node[draw, circle, color=pgreen, fill=pgreen, inner sep=0.0pt] (G4) at (0.5, 2.69) {};

\draw (G1) node[above left] {in};
\draw (G2) node[right] {out};
\draw (G3) node[below right] {in};
\draw (G4) node[above left] {out};

\draw (A) -- (B) ;
\draw (A) -- (A1) node[below] {$\scriptstyle{e_2}$};
\draw (B) -- (C) node[midway, above left] {$\scriptstyle{f_1}$};
\draw (C) -- (C2) node[above] {$\scriptstyle{f_2}$};
\draw (C) -- (D);
\draw (D) -- (A) node[midway, below right] {$\scriptstyle{e_1}$};
\draw (B) -- (D) node[midway, below right] {$\scriptstyle{e_0}$};

\draw[decorate, decoration={zigzag, segment length=4, amplitude=1.5}] (Z1) -- (Z2); 

\draw[color=blue, postaction={on each segment={mid arrow=blue}}] (A) -- (Z1)  {};
\draw[color=blue, postaction={on each segment={mid arrow=blue}}] (Z1) -- (B)  {};
\draw[color=blue, postaction={on each segment={mid arrow=blue}}] (Z1) -- (D)  {};
\draw[color=blue, postaction={on each segment={mid arrow=blue}}] (Z2) -- (B)  {};
\draw[color=blue, postaction={on each segment={mid arrow=blue}}] (Z2) -- (D)  {};

\draw[color=blue, postaction={on each segment={mid arrow=blue}}] (A) -- (Z3);
\draw[color=blue, postaction={on each segment={mid arrow=blue}}] (A) -- (Z4);
\draw[color=blue, postaction={on each segment={mid arrow=blue}}] (C) -- (Z2);
\draw[color=blue, postaction={on each segment={mid arrow=blue}}] (C) -- (Z5);
\draw[color=blue, postaction={on each segment={mid arrow=blue}}] (C) -- (Z7);

\draw plot [smooth, tension=1.1] coordinates {(-1.5, 2.3) (A) (0.7, 2.9)};
\draw plot [smooth, tension=1.2] coordinates {(3, 2.5) (B) (3.5, 0.8)};
\draw plot [smooth, tension=1.2] coordinates {(3.7, -0.8) (C) (2, -1.5)};
\draw plot [smooth, tension=1.2] coordinates { (-1, 0.8) (D) (-0.3, -1)};

\draw[color=pgreen, postaction={on each segment={other arrow=pgreen}}] (G1) -- (G2);
\draw[color=pgreen, postaction={on each segment={other arrow=pgreen}}] (G3) -- (G4);
\centerarc[pgreen](OrA)(8:-170:0.15)
\centerarc[pgreen](OrC)(15:195:0.15)

\end{tikzpicture}

\end{minipage}
\hfill 
\begin{minipage}{0.45\textwidth}


\begin{tikzpicture}
\tikzset{
  on each segment/.style={
    decorate,
    decoration={
      show path construction,
      moveto code={},
      lineto code={
        \path [#1]
        (\tikzinputsegmentfirst) -- (\tikzinputsegmentlast);
      },
      curveto code={
        \path [#1] (\tikzinputsegmentfirst)
        .. controls
        (\tikzinputsegmentsupporta) and (\tikzinputsegmentsupportb)
        ..
        (\tikzinputsegmentlast);
      },
      closepath code={
        \path [#1]
        (\tikzinputsegmentfirst) -- (\tikzinputsegmentlast);
      },
    },
  },
  mid arrow/.style={postaction={decorate,decoration={
        markings,
        mark=at position .5 with {\arrow[#1]{stealth}}
      }}},
  other arrow/.style={postaction={decorate,decoration={
        markings,
        mark=at position .75 with {\arrow[#1]{stealth}}
      }}},
}

\node[draw, rectangle, fill=orange, inner sep=1pt] (OrA) at (-0.72, 2.44) {};
\node[draw, circle, fill=black, inner sep=1pt] (A) at (0, 2.5) {};
\node[draw, circle, fill=black, inner sep=1pt] (B) at (3, 1.5) {};
\node[draw, rectangle, fill=orange, inner sep=1pt] (OrB) at (2.95, 1.75) {};
\node[draw, circle, fill=black, inner sep=1pt] (C) at (2.5, -1) {};
\node[draw, rectangle, fill=orange, inner sep=1pt] (OrC) at (3, -0.82) {};
\node[draw, circle, fill=black, inner sep=1pt] (D) at (-0.5, 0) {};
\node[draw, circle, fill=black, inner sep=0.0pt] (D1) at (0.45,-1.3) {};
\node[draw, rectangle, fill=orange, inner sep=1pt] (OrD) at (-0.35, -0.5) {};

\node[draw, circle, fill=blue, inner sep=0.5pt] (Z1) at (1, 1.2) {};
\node (Z1aux) at (0.8, 1.4) {};
\node[draw, circle, fill=blue, inner sep=0.5pt] (Z2) at (1.5, 0.3) {};

\node[draw, circle, fill=blue, inner sep=0.5pt] (Z3) at (2.6, 2.4) {};
\node[draw, circle, fill=blue, inner sep=0.5pt] (Z4) at (0.25, -1.5) {};
\node[draw, circle, fill=blue, inner sep=0.5pt] (Z5) at (0.42, -1) {};

\node[draw, circle, color=pgreen, fill=pgreen, inner sep=0.0pt] (G1) at (-0.75, 0.5) {};
\node[draw, circle, color=pgreen, fill=pgreen, inner sep=0.0pt] (G2) at (2.93, 2.2) {};
\node[draw, circle, color=pgreen, fill=pgreen, inner sep=0.0pt] (G3) at (3.29, 0.96) {};
\node[draw, circle, color=pgreen, fill=pgreen, inner sep=0.0pt] (G4) at (-0.3, -0.8) {};

\draw (D) node[left] {in};
\draw (G2) node[right] {out};
\draw (G3) node[above right] {in};
\draw (G4) node[left] {out};

\draw (A) -- (B) node[midway, above] {$\scriptstyle{g_1}$};
\draw (B) -- (C)  ;
\draw (D) -- (D1) node[near end, below] {$\scriptstyle{h_2}$};
\draw (C) -- (D)  node[midway, below] {$\scriptstyle{h_1}$};
\draw (D) -- (A) ;
\draw (B) -- (D) node[midway, right] {$\scriptstyle{e_0}$};

\draw[decorate, decoration={zigzag, segment length=4, amplitude=1.5}] (Z1) -- (Z2); 

\draw[color=blue, postaction={on each segment={mid arrow=blue}}] (A) -- (Z1)  {};
\draw[color=blue, postaction={on each segment={mid arrow=blue}}] (Z1) -- (B)  {};
\draw[color=blue, postaction={on each segment={mid arrow=blue}}] (Z1) -- (D)  {};
\draw[color=blue, postaction={on each segment={mid arrow=blue}}] (Z2) -- (B)  {};
\draw[color=blue, postaction={on each segment={mid arrow=blue}}] (Z2) -- (D)  {};
\draw[color=blue, postaction={on each segment={mid arrow=blue}}] (C) -- (Z2)  {};

\draw[color=blue, postaction={on each segment={mid arrow=blue}}] (Z3)--(B);
\draw[color=blue, postaction={on each segment={mid arrow=blue}}] (Z4)--(D);
\draw[color=blue, postaction={on each segment={mid arrow=blue}}] (Z5)--(D);

\draw plot [smooth, tension=1.2] coordinates {(-1, 2.5) (A) (1, 3)};
\draw plot [smooth, tension=1.2] coordinates {(3, 2.5) (B) (3.5, 0.8)};
\draw plot [smooth, tension=1.2] coordinates {(3.2, -0.8) (C) (1.5, -1.5)};
\draw plot [smooth, tension=1.2] coordinates { (-1, 0.8) (D) (-0.3, -1)};

\draw[color=pgreen, postaction={on each segment={mid arrow=pgreen}}] (G1) -- (G2);
\draw[color=pgreen, postaction={on each segment={mid arrow=pgreen}}] (G3) -- (G4);

\centerarc[pgreen](OrB)(100:280:0.15)
\centerarc[pgreen](OrD)(-80:105:0.15)

\end{tikzpicture}

\end{minipage}
    \caption{The algorithm of Proposition~\ref{prop:FGalgo} applied to~$e_0$ with sign $+1$ (left) and $-1$ (right)}\label{cap:FG1}
  \end{figure}
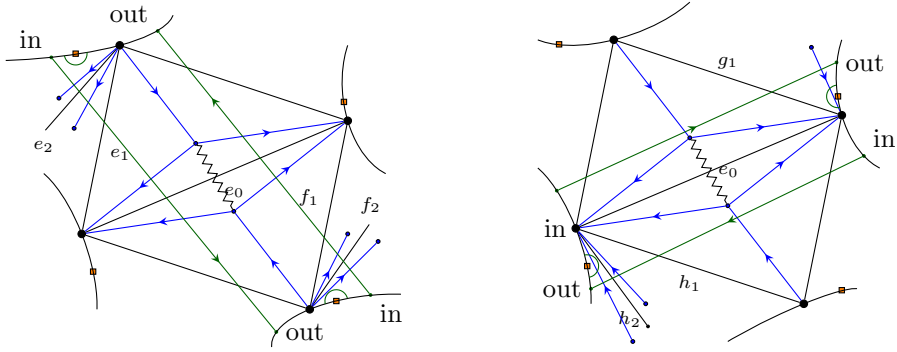

\begin{lemma} \label{le:FGpm1}
Fix an edge~$e_0$ of the triangulation~$\TT$. Let $\calL_{\pm e}$ the lamination obtained by applying Proposition~\ref{prop:FGalgo} to the tuple $n_e = \pm 1$ and $n_j=0$ for all other edges. Then
\bas
\frac12 F(\calL_{+e_0},\theta) &\= [\dnl\gamma_{e_0}] + \sum_{i=1}^k
(1+[\dnl\gamma_{e_0}]) \prod_{m=1}^i[\dnl\gamma_{e_j}]^{-1} + \sum_{j=1}^\ell
(1+[\dnl\gamma_{e_0}]) \prod_{n=1}^i[\dnl\gamma_{f_j}]^{-1} 
\\
& \qquad + \sum_{i,j=1}^{\ell,k} [\dnl\gamma_{e_0}](1+[\dnl\gamma_{e_0}]^{-1})^2 \prod_{m=1}^i 
[\dnl\gamma_{e_m}]^{-1}\prod_{n=1}^j[\dnl\gamma_{f_n}]^{-1} \qquad \text{and}\\
\frac12 F(\calL_{-e_0},\theta) &\= [\dnl\gamma_{e_0}]^{-1} + \sum_{i=1}^k
[\dnl\gamma_{e_0}]^{-1} \prod_{m=1}^i[\dnl\gamma_{g_m}]^{-1} + \sum_{i=1}^\ell
[\dnl\gamma_{e_0}]^{-1}\prod_{n=1}^i[\dnl\gamma_{h_n}]^{-1}
\\
& \qquad + \sum_{i,j=1}^{\ell,k} [\dnl\gamma_{e_0}]^{-1} \prod_{m=1}^i 
[\dnl\gamma_{g_m}]^{-1}\prod_{n=1}^j[\dnl\gamma_{h_n}]^{-1}\,,
\eas
where~$e_i,g_i$  ($i = 1,\ldots,k$) and $f_i,h_i$ ($i=1,\ldots,\ell$) are the edges adjacent as in Figure~\ref{cap:FG1} to the quadrilateral of the two triangles adjacent to $e_0$.

\end{lemma}
\par

\begin{proof} Thanks to homotopy invariance, including Proposition~\ref{prop:boundaryhomotopy} we may work with the segments of $\calL_+$ as drawn on one sheet in Figure~\ref{cap:FG1}, and their $\sigma$-images on the other sheet. There is an obvious lift that stays on this sheet, as drawn,  homotopic to $\dnl\gamma_e$. For the one that stays on the other sheet the preferred sheet rule implies the same orientation, hence the factor two.
\par
Our goal is to determine the possible detours along the segments with positive weight that are compatible with the preferred sheet rule. For the downward path on the left we can take a detour just before or just after crossing~$e_0$ (but not both). In order to be incoming at the next boundary, we must take one more detour between $e_i$ and $e_{i+1}$. This leads to the factors
$\prod_{j=1}^i[\dnl\gamma_{e_j}]^{-1}$, with or without a prefactor $[\dnl\gamma_{e_0}])$. (Recall that $[\dnl\gamma_{-e_j}]= [\dnl\gamma_{e_j}]^{-1}$.)
The same reasoning applies to the upward path on the right, leading to detours that contribute products of $[\dnl\gamma_{f_j}]^{-1}$, again with or without a prefactor $[\dnl\gamma_{e_0}])$ depending on where the first detour has been taken.
\end{proof}
\par
The important feature in the shape of the lifts in Lemma~\ref{le:FGpm1} is that it gives the desired term $[\dnl\gamma_{e_0}]^{\pm 1}$ up to a polynomial in only the inverses of the standard hat-homology generators $[\dnl\gamma_{e}]^{- 1}$ for $e \in \TT$. This gives a straightforward algorithm for approximation in the profinite completion since the companion terms in inverses eventually have large central charge. This is made precise in the following proposition. We need an extra observation since we later need  the polynomials to have central charge within a cone translate $K = \Delta_2^t$ such that the elements of $\wh{G}_\Delta$ act continuously on the corresponding submodule of the module of tame paths.
\par
\begin{proposition} \label{prop:approxFL} Suppose that $\theta = 0$ is the horizontal direction and the sign of the hat-homology classes is chosen such that $Z(\wh\gamma) \in \HH$, as usual. Then for any $e_0 \in \TT$ there are sequences $P^\pm_n$ of polynomials in the lifts $F(\calL_{+e},0)$ and $F(\calL_{-e}, 0)$ where $e \in \TT$, such that $\sum_{i=1}^n P_n^\pm \to [\dnl\gamma_{e_0}]^{\pm 1}$ in the topology of $\wh{H}_\Delta^\abs$.
\par
Moreover, there exists a cone~$\Delta_2$ strictly contained in the lower half plane and~$t \in \CC$ with the following property: Any polynomial~$P$ in $F(\calL_{+e},0)$ and $F(\calL_{-e},0)$ whose total degree in the lifts of positive generators $F(\calL_{+e},0)$ is at most one, is an element of $\wh{H}_\Delta^\abs$ with $Z(S(P)) \subset \Delta_2^t$.
\end{proposition}
\par
\begin{proof}
We start with the second statement. Any polynomial~$P$ in $[\wh{\gamma}_e]^{-1}$ for all $e \in \TT$ has central charges supported in some cone~$\Delta_2$ strictly contained in $\HH^-$ by the normalization of hat-homology classes and since $|\TT|$ is finite. By Lemma~\ref{le:FGpm1} this takes care of polynomials in $F(\calL_{-e},0)$ and one monomial in $F(\calL_{+e},0)$ involves again by Lemma~\ref{le:FGpm1} at most a single term $[\wh{\gamma}_e]$, which is taken care of by the choice of~$t$.
\par
For the first statement we will also use polynomials~$P_n^\pm$ whose total degree in the lifts of positive generators $F(\calL_{+e},0)$ is at most one and observe moreover that the joint support of the central charges all such polynomials (in $\Delta_2^t$) is a discrete subset of this cone, again and since $|\TT|$ is finite. 
\par
As algorithm for approximation we start with $P_1^\pm = F(\calL_{\pm e_0},0)$. In the subsequent steps we (still) work with the total degree in the variables $[\wh\gamma_{\pm e}]$, i.e., in the $n$-th step we consider the degree~$n$ homogeneous part~$Q$ of $[\dnl\gamma_{\pm e_0}]- P_{n-1}^\pm$ and set $P_n = P_{n-1} + \wt{Q}$, where in~$\wt{Q}$ we replaced each variable $[\wh\gamma_{\pm e}]$ by the corresponding $F(\calL_{\pm e},0)$. Observe that the shape of the lifts Lemma~\ref{le:FGpm1} imply inductively that~$Q$ has indeed degree at most one in the positive generators.
\par
Finally we observe that the absolute value of the central charge of a homogeneous polynomial~$P$ with degree at most one in the positive generators has a lower bound that is a linear function in the degree of~$P$. In fact,
\be
\Im Z(P) \geq (n-1) \min_{e \in \TT} \Im Z(\wh\gamma_e) - \max_{e \in \TT} \Im Z(\wh\gamma_e)\,.
\ee
This implies that $[\dnl\gamma_{\pm e_0}]- P_{n}^\pm$ indeed converges to zero in the topology of $\wh{H}_\Delta^\abs$. 
\end{proof}

\subsection{Variations of the lift function for ciliated surfaces with holes}
\label{sec:withholes}

We describe a version of the lift function that additionally depends on the choice of a base point on $\wt{\Sigma}$. This is necessary to approximate the closed paths associated to the edges of an triangulation by lifts in the presence of holes.
\par 
Let $\calL$ be a lamination with loop decomposition $\{\wp_1,\ldots,\wp_n\}$. Now, we treat the general case, that is, some of the $\wp_j$ might just be closed loops on $\Sigma^\circ$. Note that the differential of $c \circ \pi$ induces a map on circle bundles $\wt{\Pi}: \wt{\Sigma} \to \wt{X}$. For a loop $\wp_j$ we call a point $\wt{z_j} \in \wt{\Sigma}$ such that $\wt{\Pi}(\wt{z}_j)$ is a tangent direction to $\wp_j$ a \emph{directed base point}. A detour $\dnl{D}_{\bfy_j}$ of~$\wp_j$ is called \emph{based at~$\wt{z}_j$} if it starts and ends at $\wt{z}_j$. The bullet in the following 'based' versions of the lift of a lamination should remind of the additional choice of directed base points.
\par
\begin{proposition}
  Let $(X,q,\theta)$ be a saddle-free directed quadratic differential of infinite area. Let~$\calL$ be a lamination with support in the smooth part $\XX^{\ast}= \XX \setminus \calZ \cup \calP$. Choose a loop decomposition $\{\wp_1,\ldots,\wp_n\}$ and for each $\wp_j$ a directed base point $\wt{z}_j$. Then the based path lifting 
\be
\bfF_\bullet(\calL,\theta) \= \sum_{j=1}^n \sum_{N \in \NN}\,\, \sum_{\bfy_j = (y_1,\ldots,y_N)} \dnl{D}_{\bfy_j}\,,
\ee
where the last sum is over all detours of~$\wp_j$ base at $\wt{z}_j$ obeying the preferred-sheet rule, is well-defined as an element of $\wh\Pi_\Delta(\wt{\Sigma}^\circ)$ and thus defines an element $F(\calL,\theta) \in \wh{H}_\Delta(\wt{\Sigma}^\circ)$. In fact, by construction $F(\calL,\theta) \in \wh{H}_\Delta^\abs$.
\end{proposition}
Note that by only considering the lifts based at a fixed $\wt{z}_i$ the lift is still homotopy invariant in the sense of Theorem~\ref{theo:homotopy_invariance}~i). However, restricted to closed paths it looses the property to be homomorphism for the concatenation of closed paths on~$X$. 
\par
We observe that the above proof of Lemma~\ref{le:FGpm1} gives immediately the 'based' version:
\par
\begin{lemma} \label{le:FGpm1based}
  Let $\calL$ be a lamination without loops. Choose the base point $\wt{z}_i$ close enough to and pointing away from an incoming boundary component for each $\wp_i$. Then the formulas of Lemma~\ref{le:FGpm1} hold for $F_\bullet(\calL,\theta)$ without the factor $1/2$ on the left hand side.
\end{lemma}
\par
Recall that double poles of the quadratic differential correspond to holes of the ciliated surface. In this case, we can not avoid working with laminations containing closed loops. There are two new situations to consider in comparison with Lemma~\ref{le:FGpm1}. For an edge~$e$ of the triangulation we call a vertex of the enclosing quadrilateral \emph{inner vertex}, if it is a vertex of~$e$. The other two vertices are called \emph{outer vertices}. The two extra case are at edges~$e$, where one or both outer vertices are holes. The following lemma computes the based lift of the lamination obtained from Proposition~\ref{prop:FGalgo} in the four cases, the number of outer vertices and the sign for~$e$.
\begin{figure}
  \begin{minipage}{0.45\textwidth}


\begin{tikzpicture}
\tikzset{
  on each segment/.style={
    decorate,
    decoration={
      show path construction,
      moveto code={},
      lineto code={
        \path [#1]
        (\tikzinputsegmentfirst) -- (\tikzinputsegmentlast);
      },
      curveto code={
        \path [#1] (\tikzinputsegmentfirst)
        .. controls
        (\tikzinputsegmentsupporta) and (\tikzinputsegmentsupportb)
        ..
        (\tikzinputsegmentlast);
      },
      closepath code={
        \path [#1]
        (\tikzinputsegmentfirst) -- (\tikzinputsegmentlast);
      },
    },
  },
  mid arrow/.style={postaction={decorate,decoration={
        markings,
        mark=at position .5 with {\arrow[#1]{stealth}}
      }}},
  other arrow/.style={postaction={decorate,decoration={
        markings,
        mark=at position .75 with {\arrow[#1]{stealth}}
      }}},
}

\node[draw, rectangle, fill=orange, inner sep=1pt] (OrA) at (-0.58, 2.39) {};
\node[draw, circle, fill=black, inner sep=1pt] (A) at (0, 2.5) {};
\node[inner sep=0.0pt] (A1) at (-1, 1.5) {};
\node[draw, circle, fill=black, inner sep=1pt] (B) at (3, 1.5) {};
\node[draw, rectangle, fill=orange, inner sep=1pt] (OrB) at (2.95, 1.75) {};
\node[draw, circle, fill=black, inner sep=1pt] (C) at (2.5, -1) {};
\node[draw, circle, fill=black, inner sep=1pt] (CB) at (2.25, -2.25) {};
\node[inner sep=0.0pt] (C1) at (3.7, -0.05) {};
\node[inner sep=0.0pt] (C2) at (3.3, 0.15) {};
\node[draw, circle, fill=black, inner sep=1pt] (D) at (-0.5, 0) {};
\node[draw, circle, fill=black, inner sep=1pt] (CD) at (4, -1.5) {};
\node[draw, rectangle, fill=orange, inner sep=1pt] (OrD) at (-0.35, -0.5) {};

\node[draw, circle, fill=blue, inner sep=0.5pt] (Z1) at (1, 1.2) {};
\node[draw, circle, fill=blue, inner sep=0.5pt] (Z2) at (1.5, 0.3) {};
\node[draw, circle, fill=blue, inner sep=0.5pt] (Z3) at (-0.8, 2) {};
\node[draw, circle, fill=blue, inner sep=0.5pt] (Z4) at (-0.6, 1.6) {};

\node[draw, circle, fill=blue, inner sep=0.5pt] (Z5) at (3.5, -2) {};
\node[draw, circle, fill=blue, inner sep=0.5pt] (Z6) at (1.6, -1.9) {};
\node[draw, circle, fill=blue, inner sep=0.5pt] (Z7) at (3.4, -0.1) {};

\node[draw, circle, color=pgreen, fill=pgreen, inner sep=0.0pt] (G1) at (-0.9, 2.34) {};
\node[inner sep=0.0pt] (G2) at (2.12, -1.33) {};
\node[inner sep=0.0pt] (G3) at (2.90, -0.7) {};
\node[draw, circle, color=pgreen, fill=pgreen, inner sep=0.0pt] (G4) at (0.5, 2.69) {};

\draw (G1) node[above] {in};
\draw (A) node[above] {out};

\draw (A) -- (B) ;
\draw (A) -- (A1) node[below] {$\scriptstyle{e_2}$};
\draw (B) -- node[right]{$\scriptstyle{f_1}$}(C)--node[right]{$\scriptstyle{f_3}$}(CB);
\draw (CD)--node[above]{$\scriptstyle{f_2}$}(C) -- node[below]{$\scriptstyle{f_4}$}(D);
\draw (D) -- node[left]{$\scriptstyle{e_1}$}(A);
\draw (B) -- (D) node[midway, right] {$\scriptstyle{e_0}$};

\draw[decorate, decoration={zigzag, segment length=4, amplitude=1.5}] (Z1) -- (Z2); 

\draw[color=blue, postaction={on each segment={mid arrow=blue}}] (A) -- (Z1)  {};
\draw[color=blue, postaction={on each segment={mid arrow=blue}}] (Z1) -- (B)  {};
\draw[color=blue, postaction={on each segment={mid arrow=blue}}] (Z1) -- (D)  {};
\draw[color=blue, postaction={on each segment={mid arrow=blue}}] (Z2) -- (B)  {};
\draw[color=blue, postaction={on each segment={mid arrow=blue}}] (Z2) -- (D)  {};

\draw[color=blue, postaction={on each segment={mid arrow=blue}}] (A) -- (Z3);
\draw[color=blue, postaction={on each segment={mid arrow=blue}}] (A) -- (Z4);
\draw[color=blue, postaction={on each segment={mid arrow=blue}}] (C) -- (Z2);
\draw[color=blue, postaction={on each segment={mid arrow=blue}}] (C) -- (Z5);
\draw[color=blue, postaction={on each segment={mid arrow=blue}}] (C) -- (Z6);
\draw[color=blue, postaction={on each segment={mid arrow=blue}}] (C) -- (Z7);

\draw plot [smooth, tension=1.1] coordinates {(-1.5, 2.3) (A) (0.7, 2.9)};
\draw plot [smooth, tension=1.2] coordinates {(3, 2.5) (B) (3.5, 0.8)};
\draw plot [smooth, tension=1.2] coordinates { (-1, 0.8) (D) (-0.3, -1)};

\draw[color=pgreen, postaction={on each segment={other arrow=pgreen}}] (G1) -- (G2);
\draw[color=pgreen, postaction={on each segment={other arrow=pgreen}}] (G3) -- (G4);
\centerarc[pgreen](C)(0:360:0.15)
\centerarc[pgreen](OrA)(7:-170:0.15)
\centerarc[pgreen](C)(41:-144:0.5)
\end{tikzpicture}

\end{minipage}
\hfill 
\begin{minipage}{0.45\textwidth}


\begin{tikzpicture}
\tikzset{
  on each segment/.style={
    decorate,
    decoration={
      show path construction,
      moveto code={},
      lineto code={
        \path [#1]
        (\tikzinputsegmentfirst) -- (\tikzinputsegmentlast);
      },
      curveto code={
        \path [#1] (\tikzinputsegmentfirst)
        .. controls
        (\tikzinputsegmentsupporta) and (\tikzinputsegmentsupportb)
        ..
        (\tikzinputsegmentlast);
      },
      closepath code={
        \path [#1]
        (\tikzinputsegmentfirst) -- (\tikzinputsegmentlast);
      },
    },
  },
  mid arrow/.style={postaction={decorate,decoration={
        markings,
        mark=at position .5 with {\arrow[#1]{stealth}}
      }}},
  other arrow/.style={postaction={decorate,decoration={
        markings,
        mark=at position .75 with {\arrow[#1]{stealth}}
      }}},
}

\node[draw, circle, fill=black, inner sep=1pt] (A) at (0, 2.5) {};
\node[inner sep=0.0pt] (A1) at (-1, 1.35) {};
\node[draw, circle, fill=black, inner sep=1pt] (B) at (3, 1.5) {};
\node[draw, rectangle, fill=orange, inner sep=1pt] (OrB) at (3.12, 1.2) {};
\node[draw, circle, fill=black, inner sep=1pt] (C) at (2.5, -1) {};
\node[draw, circle, fill=black, inner sep=1pt] (CB) at (2.25, -2.25) {};
\node[inner sep=0.0pt] (C1) at (3.7, -0.05) {};
\node[inner sep=0.0pt] (C2) at (3.3, 0.15) {};
\node[draw, circle, fill=black, inner sep=1pt] (D) at (-0.5, 0) {};
\node[draw, circle, fill=black, inner sep=1pt] (CD) at (4, -1.5) {};
\node[draw, rectangle, fill=orange, inner sep=1pt] (OrD) at (-0.64, 0.33) {};
\node[draw, circle, fill=black, inner sep=1pt] (AD) at (-0.6,3.5) {};

\node[draw, circle, fill=blue, inner sep=0.5pt] (Z1) at (1, 1.2) {};
\node[draw, circle, fill=blue, inner sep=0.5pt] (Z2) at (1.5, 0.3) {};
\node[draw, circle, fill=blue, inner sep=0.5pt] (Z3) at (-0.8, 2) {};
\node[draw, circle, fill=blue, inner sep=0.5pt] (Z4) at (0.7, 3) {};

\node[draw, circle, fill=blue, inner sep=0.5pt] (Z5) at (3.5, -2) {};
\node[draw, circle, fill=blue, inner sep=0.5pt] (Z6) at (1.6, -1.9) {};
\node[draw, circle, fill=blue, inner sep=0.5pt] (Z7) at (3.4, -0.1) {};

\node[inner sep=0.0pt] (G1) at (-0.44, 2.26) {};
\node[inner sep=0.0pt] (G2) at (2.12, -1.33) {};
\node[inner sep=0.0pt] (G3) at (2.90, -0.7) {};
\node[inner sep=0.0pt] (G4) at (0.4, 2.8) {};


\draw (A) --node[below]{$\scriptstyle{e_3}$} (B) ;
\draw (B) -- node[right]{$\scriptstyle{f_1}$}(C)--node[right]{$\scriptstyle{f_3}$}(CB);
\draw (CD)--node[above]{$\scriptstyle{f_2}$}(C) -- node[below]{$\scriptstyle{f_4}$}(D);
\draw (D) -- node[left]{$\scriptstyle{e_1}$}(A)--node[left]{$\scriptstyle{e_2}$}(AD);
\draw (B) -- (D) node[midway, right] {$\scriptstyle{e_0}$};

\draw[decorate, decoration={zigzag, segment length=4, amplitude=1.5}] (Z1) -- (Z2); 

\draw[color=blue, postaction={on each segment={mid arrow=blue}}] (A) -- (Z1)  {};
\draw[color=blue, postaction={on each segment={mid arrow=blue}}] (Z1) -- (B)  {};
\draw[color=blue, postaction={on each segment={mid arrow=blue}}] (Z1) -- (D)  {};
\draw[color=blue, postaction={on each segment={mid arrow=blue}}] (Z2) -- (B)  {};
\draw[color=blue, postaction={on each segment={mid arrow=blue}}] (Z2) -- (D)  {};

\draw[color=blue, postaction={on each segment={mid arrow=blue}}] (A) -- (Z3);
\draw[color=blue, postaction={on each segment={mid arrow=blue}}] (A) -- (Z4);
\draw[color=blue, postaction={on each segment={mid arrow=blue}}] (C) -- (Z2);
\draw[color=blue, postaction={on each segment={mid arrow=blue}}] (C) -- (Z5);
\draw[color=blue, postaction={on each segment={mid arrow=blue}}] (C) -- (Z6);
\draw[color=blue, postaction={on each segment={mid arrow=blue}}] (C) -- (Z7);

\draw plot [smooth, tension=1.2] coordinates {(3, 2.5) (B) (3.5, 0.8)};
\draw plot [smooth, tension=1.2] coordinates { (-1, 0.8) (D) (-0.3, -1)};

\draw[color=pgreen] (G1) -- (G2);
\draw[color=pgreen] (G3) -- (G4);
\centerarc[pgreen](A)(0:360:0.15)
\centerarc[pgreen](C)(0:360:0.15)
\centerarc[pgreen](A)(32:220:0.5)
\centerarc[pgreen](C)(41:-144:0.5)
\end{tikzpicture}

\end{minipage}
  \caption{Left: One outer vertex is a hole. Right: Both outer vertices are holes.}\label{fig:FL_plus_holes}
\end{figure}
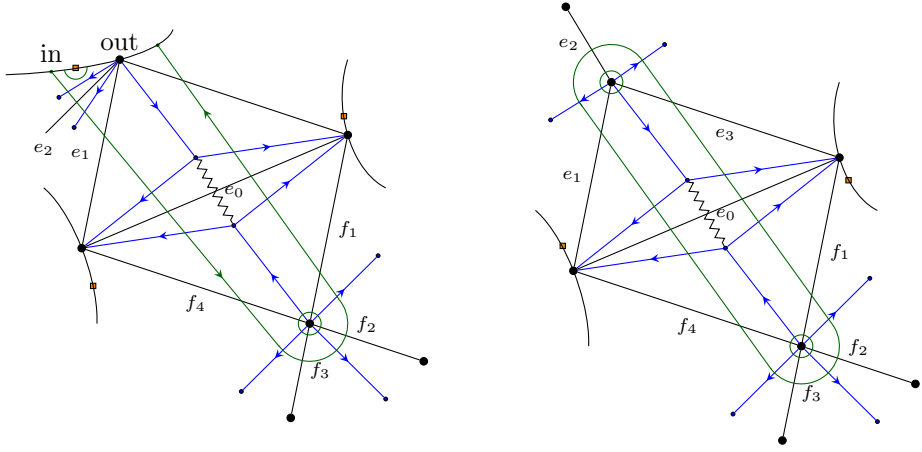
\par
\begin{lemma}\label{le:FGpm1_holes} For appropriate choices of the directed base points, the based lifts of the laminations given by applying Proposition~\ref{prop:FGalgo} are given as follows: 
  \begin{itemize}
  \item[i)] \textbf{One outer vertex is a hole:} \quad Consider an edge $e$ such that exactly one outer vertex of the associated quadrilateral is a hole  as depicted on the left in Figure~\ref{fig:FL_plus_holes}. Denote by $\epsilon$ the loop around the hole. Let $\calL_e$ be the lamination obtained by applying Proposition~\ref{prop:FGalgo} to the tuple $n_e =  1$ and $n_{e'}=0$ for all other edges $e' \neq e$. Then there is a loop decomposition $\{\wp_{e},\epsilon\}$ of $\calL_e$. Choose a directed base point $\wt{z}$ on $\wp_{e}$ close enough to and pointing away from an incoming boundary component. Then 
\bas
F_\bullet(\wp_{e},\theta) &\= [\dnl\gamma_{e}+\dnl\varepsilon] + \sum_{i=1}^k
\left([\dnl\varepsilon](1+[\dnl\gamma_{e}])+[\sigma(\dnl\varepsilon)](1+[\dnl\gamma_e]^{-1})\right) \prod_{m=1}^i[\dnl\gamma_{e_m}]^{-1} \\
& \qquad + \sum_{j=1}^\ell
[\dnl\varepsilon](1+[\dnl\gamma_{e}]) \prod_{n=1}^j[\dnl\gamma_{f_n}]^{-1}
\\
& \qquad + \sum_{i,j=1}^{\ell,k} [\dnl\gamma_{e_0}+\dnl\varepsilon](1+[\dnl\gamma_{e_0}]^{-1})^2 \prod_{m=1}^i 
[\dnl\gamma_{e_m}]^{-1}\prod_{n=1}^j[\dnl\gamma_{f_n}]^{-1} 
\eas
where~$e_i$  ($i = 1,\ldots,k$) and $f_j$ ($j=1,\ldots,\ell$) are the edges adjacent as in Figure~\ref{fig:FL_plus_holes} to the quadrilateral of the two triangles adjacent to~$e$ and $\varepsilon$ is a closed path on $\Sigma$ with $\pi(\varepsilon)=\epsilon$. 
  \item[ii)] \textbf{Both outer vertices are holes:} \quad
  Consider an edge $e$ such that both outer vertices of the associated quadrilateral are holes  as depicted on the right in Figure~\ref{fig:FL_plus_holes}. Denote by $\epsilon_1,\epsilon_2$ the loops around the holes. Let $\calL_{e}$ be the lamination obtained by applying Proposition~\ref{prop:FGalgo} to the tuple $n_e = 1$ and $n_j=0$ for all other edges. Then there is a loop decomposition $\{\wp_e,\epsilon_1,\epsilon_2\}$ of $\calL_e$. Choose a directed base point $\wt{z}$ for $\wp_{e}$ on the edge $e$, such that the trivial lift of $\wp_e$ starting at $\wt{z}$ yields $[\dnl\gamma_{e}+\dnl{\varepsilon_1}+\dnl\varepsilon_2]$, where $\varepsilon_i$ are closed paths on $\Sigma$ with $\pi(\varepsilon_i)=\epsilon_i$. Then
  \bas F_{\bullet}(\wp_{e},\theta) &\= [\dnl\gamma_{e}+\dnl\varepsilon_1+\dnl\varepsilon_2] + \sum_{i=1}^k ([\dnl\gamma_{e}+{\dnl\varepsilon_1}+\dnl\varepsilon_2]+[\sigma(\dnl\varepsilon_1)+\dnl\varepsilon_2])
  \prod_{m=1}^i[\dnl\gamma_{e_m}]^{-1} \\
 & \qquad + \sum_{j=1}^\ell
[\dnl\varepsilon_1+\dnl\varepsilon_2](1+[\dnl\gamma_{e}]) \prod_{n=1}^j[\dnl\gamma_{f_n}]^{-1}
\\
& \qquad + \sum_{i,j=1}^{\ell,k} [\dnl\varepsilon_1+\dnl\varepsilon_2](1+[\dnl\gamma_{e}]) \prod_{m=1}^i 
[\dnl\gamma_{e_m}]^{-1}\prod_{n=1}^j[\dnl\gamma_{f_n}]^{-1} 
\eas
where~$e_i$  ($i = 1,\ldots,k$) and $f_j$ ($j=1,\ldots,\ell$) are the edges adjacent as in Figure~\ref{fig:FL_plus_holes} to the quadrilateral of the two triangles adjacent to $e$.

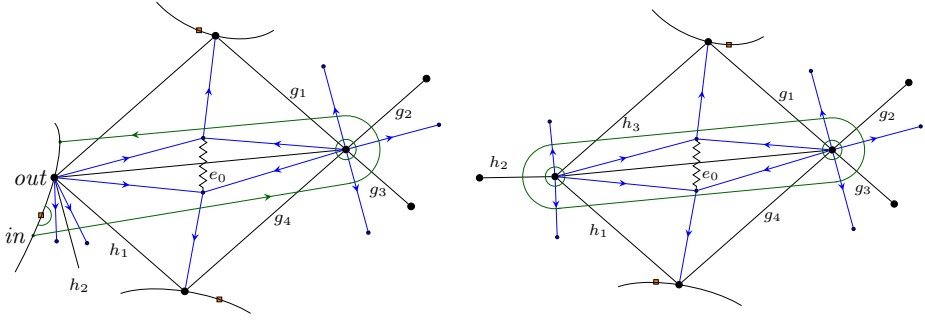
\begin{figure}
  \begin{minipage}{0.49\textwidth}
    \resizebox{6.2cm}{!}{

\begin{tikzpicture}
\begin{scope}[rotate=60]  
  \tikzset{
  on each segment/.style={
    decorate,
    decoration={
      show path construction,
      moveto code={},
      lineto code={
        \path [#1]
        (\tikzinputsegmentfirst) -- (\tikzinputsegmentlast);
      },
      curveto code={
        \path [#1] (\tikzinputsegmentfirst)
        .. controls
        (\tikzinputsegmentsupporta) and (\tikzinputsegmentsupportb)
        ..
        (\tikzinputsegmentlast);
      },
      closepath code={
        \path [#1]
        (\tikzinputsegmentfirst) -- (\tikzinputsegmentlast);
      },
    },
  },
  mid arrow/.style={postaction={decorate,decoration={
        markings,
        mark=at position .5 with {\arrow[#1]{stealth}}
      }}},
  other arrow/.style={postaction={decorate,decoration={
        markings,
        mark=at position .75 with {\arrow[#1]{stealth}}
      }}},
}

\node[draw, rectangle, fill=orange, inner sep=1pt] (OrA) at (-0.58, 2.39) {};
\node[draw, circle, fill=black, inner sep=1pt] (A) at (0, 2.5) {};
\node[inner sep=0.0pt] (A1) at (-1, 1.5) {};
\node[draw, circle, fill=black, inner sep=1pt] (B) at (3, 1.5) {};
\node[draw, rectangle, fill=orange, inner sep=1pt] (OrB) at (2.95, 1.75) {};
\node[draw, circle, fill=black, inner sep=1pt] (C) at (2.5, -1) {};
\node[draw, circle, fill=black, inner sep=1pt] (CB) at (2.25, -2.25) {};
\node[inner sep=0.0pt] (C1) at (3.7, -0.05) {};
\node[inner sep=0.0pt] (C2) at (3.3, 0.15) {};
\node[draw, circle, fill=black, inner sep=1pt] (D) at (-0.5, 0) {};
\node[draw, circle, fill=black, inner sep=1pt] (CD) at (4, -1.5) {};
\node[draw, rectangle, fill=orange, inner sep=1pt] (OrD) at (-0.35, -0.5) {};

\node[draw, circle, fill=blue, inner sep=0.5pt] (Z1) at (1.6, 0.9) {};
\node[draw, circle, fill=blue, inner sep=0.5pt] (Z2) at (0.9, 0.5) {};
\node[draw, circle, fill=blue, inner sep=0.5pt] (Z3) at (-0.8, 2) {};
\node[draw, circle, fill=blue, inner sep=0.5pt] (Z4) at (-0.6, 1.6) {};

\node[draw, circle, fill=blue, inner sep=0.5pt] (Z5) at (3.5, -2) {};
\node[draw, circle, fill=blue, inner sep=0.5pt] (Z6) at (1.6, -1.9) {};
\node[draw, circle, fill=blue, inner sep=0.5pt] (Z7) at (3.4, -0.1) {};

\node[draw, circle, color=pgreen, fill=pgreen, inner sep=0.0pt] (G1) at (-0.9, 2.34) {};
\node[inner sep=0.0pt] (G2) at (2.12, -1.33) {};
\node[inner sep=0.0pt] (G3) at (2.90, -0.7) {};
\node[draw, circle, color=pgreen, fill=pgreen, inner sep=0.0pt] (G4) at (0.5, 2.69) {};

\draw (G1) node[left] {in};
\draw (A) node[left] {out};

\draw (A) -- (B) ;
\draw (A) -- (A1) node[below] {$\scriptstyle{h_2}$};
\draw (B) -- node[right]{$\scriptstyle{g_1}$}(C)--node[below]{$\scriptstyle{g_3}$}(CB);
\draw (CD)--node[right]{$\scriptstyle{g_2}$}(C) -- node[right]{$\scriptstyle{g_4}$}(D);
\draw (D) -- node[below]{$\scriptstyle{h_1}$}(A);
\draw (A) -- (C) node[midway, below right] {$\scriptstyle{e_0}$};

\draw[decorate, decoration={zigzag, segment length=4, amplitude=1.5}] (Z1) -- (Z2); 

\draw[color=blue, postaction={on each segment={mid arrow=blue}}] (A) -- (Z1)  {};
\draw[color=blue, postaction={on each segment={mid arrow=blue}}] (C)--(Z1)  {};
\draw[color=blue, postaction={on each segment={mid arrow=blue}}] (Z1) -- (B)  {};
\draw[color=blue, postaction={on each segment={mid arrow=blue}}] (A) -- (Z2)  {};
\draw[color=blue, postaction={on each segment={mid arrow=blue}}] (Z2) -- (D)  {};

\draw[color=blue, postaction={on each segment={mid arrow=blue}}] (A) -- (Z3);
\draw[color=blue, postaction={on each segment={mid arrow=blue}}] (A) -- (Z4);
\draw[color=blue, postaction={on each segment={mid arrow=blue}}] (C) -- (Z2);
\draw[color=blue, postaction={on each segment={mid arrow=blue}}] (C) -- (Z5);
\draw[color=blue, postaction={on each segment={mid arrow=blue}}] (C) -- (Z6);
\draw[color=blue, postaction={on each segment={mid arrow=blue}}] (C) -- (Z7);

\draw plot [smooth, tension=1.1] coordinates {(-1.5, 2.3) (A) (0.7, 2.9)};
\draw plot [smooth, tension=1.2] coordinates {(3, 2.5) (B) (3.5, 0.8)};
\draw plot [smooth, tension=1.2] coordinates { (-1, 0.8) (D) (-0.3, -1)};

\draw[color=pgreen, postaction={on each segment={other arrow=pgreen}}] (G1) -- (G2);
\draw[color=pgreen, postaction={on each segment={other arrow=pgreen}}] (G3) -- (G4);
\centerarc[pgreen](C)(0:360:0.15)
\centerarc[pgreen](OrA)(7:-170:0.15)
\centerarc[pgreen](C)(41:-144:0.5)
\end{scope}
\end{tikzpicture}

    }
\end{minipage}
\begin{minipage}{0.49\textwidth}
  \resizebox{6.2cm}{!}{


\begin{tikzpicture}
  \begin{scope}[rotate=60]
\tikzset{
  on each segment/.style={
    decorate,
    decoration={
      show path construction,
      moveto code={},
      lineto code={
        \path [#1]
        (\tikzinputsegmentfirst) -- (\tikzinputsegmentlast);
      },
      curveto code={
        \path [#1] (\tikzinputsegmentfirst)
        .. controls
        (\tikzinputsegmentsupporta) and (\tikzinputsegmentsupportb)
        ..
        (\tikzinputsegmentlast);
      },
      closepath code={
        \path [#1]
        (\tikzinputsegmentfirst) -- (\tikzinputsegmentlast);
      },
    },
  },
  mid arrow/.style={postaction={decorate,decoration={
        markings,
        mark=at position .5 with {\arrow[#1]{stealth}}
      }}},
  other arrow/.style={postaction={decorate,decoration={
        markings,
        mark=at position .75 with {\arrow[#1]{stealth}}
      }}},
}

\node[draw, circle, fill=black, inner sep=1pt] (A) at (0, 2.5) {};
\node[inner sep=0.0pt] (A1) at (-1, 1.35) {};
\node[draw, circle, fill=black, inner sep=1pt] (B) at (3, 1.5) {};
\node[draw, rectangle, fill=orange, inner sep=1pt] (OrB) at (3.12, 1.2) {};
\node[draw, circle, fill=black, inner sep=1pt] (C) at (2.5, -1) {};
\node[draw, circle, fill=black, inner sep=1pt] (CB) at (2.25, -2.25) {};
\node[inner sep=0.0pt] (C1) at (3.7, -0.05) {};
\node[inner sep=0.0pt] (C2) at (3.3, 0.15) {};
\node[draw, circle, fill=black, inner sep=1pt] (D) at (-0.5, 0) {};
\node[draw, circle, fill=black, inner sep=1pt] (CD) at (4, -1.5) {};
\node[draw, rectangle, fill=orange, inner sep=1pt] (OrD) at (-0.64, 0.33) {};
\node[draw, circle, fill=black, inner sep=1pt] (AD) at (-0.6,3.5) {};

\node[draw, circle, fill=blue, inner sep=0.5pt] (Z1) at (1.6, 0.9) {};
\node[draw, circle, fill=blue, inner sep=0.5pt] (Z2) at (0.9, 0.5) {};
\node[draw, circle, fill=blue, inner sep=0.5pt] (Z3) at (-0.8, 2) {};
\node[draw, circle, fill=blue, inner sep=0.5pt] (Z4) at (0.7, 3) {};

\node[draw, circle, fill=blue, inner sep=0.5pt] (Z5) at (3.5, -2) {};
\node[draw, circle, fill=blue, inner sep=0.5pt] (Z6) at (1.6, -1.9) {};
\node[draw, circle, fill=blue, inner sep=0.5pt] (Z7) at (3.4, -0.1) {};

\node[inner sep=0.0pt] (G1) at (-0.44, 2.26) {};
\node[inner sep=0.0pt] (G2) at (2.12, -1.33) {};
\node[inner sep=0.0pt] (G3) at (2.90, -0.7) {};
\node[inner sep=0.0pt] (G4) at (0.4, 2.8) {};


\draw (A) --node[below]{$\scriptstyle{h_3}$} (B) ;
\draw (B) -- node[right]{$\scriptstyle{g_1}$}(C)--node[below]{$\scriptstyle{g_3}$}(CB);
\draw (CD)--node[right]{$\scriptstyle{g_2}$}(C) -- node[right]{$\scriptstyle{g_4}$}(D);
\draw (D) -- node[left]{$\scriptstyle{h_1}$}(A)--node[above left]{$\scriptstyle{h_2}$}(AD);
\draw (A) -- (C) node[midway, below right] {$\scriptstyle{e_0}$};

\draw[decorate, decoration={zigzag, segment length=4, amplitude=1.5}] (Z1) -- (Z2); 

\draw[color=blue, postaction={on each segment={mid arrow=blue}}] (A) -- (Z1)  {};
\draw[color=blue, postaction={on each segment={mid arrow=blue}}] (Z1) -- (B)  {};
\draw[color=blue, postaction={on each segment={mid arrow=blue}}] (C)--(Z1)  {};
\draw[color=blue, postaction={on each segment={mid arrow=blue}}] (A)--(Z2) {};
\draw[color=blue, postaction={on each segment={mid arrow=blue}}] (Z2) -- (D)  {};

\draw[color=blue, postaction={on each segment={mid arrow=blue}}] (A) -- (Z3);
\draw[color=blue, postaction={on each segment={mid arrow=blue}}] (A) -- (Z4);
\draw[color=blue, postaction={on each segment={mid arrow=blue}}] (C) -- (Z2);
\draw[color=blue, postaction={on each segment={mid arrow=blue}}] (C) -- (Z5);
\draw[color=blue, postaction={on each segment={mid arrow=blue}}] (C) -- (Z6);
\draw[color=blue, postaction={on each segment={mid arrow=blue}}] (C) -- (Z7);

\draw plot [smooth, tension=1.2] coordinates {(3, 2.5) (B) (3.5, 0.8)};
\draw plot [smooth, tension=1.2] coordinates { (-1, 0.8) (D) (-0.3, -1)};

\draw[color=pgreen] (G1) -- (G2);
\draw[color=pgreen] (G3) -- (G4);
\centerarc[pgreen](A)(0:360:0.15)
\centerarc[pgreen](C)(0:360:0.15)
\centerarc[pgreen](A)(32:220:0.5)
\centerarc[pgreen](C)(41:-144:0.5)
\end{scope}
\end{tikzpicture}

  }
\end{minipage}
  \caption{Left: One inner vertex is a hole. Right: Both inner vertices are holes.}\label{fig:FL_minus_holes}
\end{figure}

\item[iii)] \textbf{One inner vertex is a hole, negative sign:} \quad Consider an edge $e$ such that exactly one inner vertex is a hole as depicted on the left in Figure~\ref{fig:FL_minus_holes}. Denote by $\epsilon$ a loop around the hole. Let $\calL_{-e}$ be the lamination obtained by applying Proposition~\ref{prop:FGalgo} to the tuple $n_e = - 1$ and $n_{e'}=0$ for all other edges $e' \neq e$. Then there is a unique loop decomposition $\{\wp_{-e},\epsilon\}$ of $\calL_e$. Choose a directed base point $\wt{z}$ on $\wp_{-e}$ close enough to and pointing away from an incoming boundary component. Then
\bas
F_\bullet(\wp_{-e},\theta) &\= [\dnl\gamma_{e}]^{-1}[\dnl\varepsilon] + \sum_{i=1}^k
[\dnl\gamma_{e}]^{-1}[\dnl\varepsilon] \prod_{m=1}^i[\dnl\gamma_{g_m}]^{-1} \\
& \qquad + \sum_{j=1}^\ell ([\dnl\gamma_{e}]^{-1}[\dnl\varepsilon]+[\sigma(\dnl\varepsilon)]) \prod_{n=1}^j[\dnl\gamma_{h_n}]^{-1} \\
& \qquad + \sum_{i,j=1}^{\ell,k} [\dnl\gamma_{e}]^{-1}[\dnl\varepsilon] \prod_{m=1}^i [\dnl\gamma_{g_m}]^{-1}\prod_{n=1}^j[\dnl\gamma_{h_n}]^{-1} 
\eas
where~$g_i$  ($i = 1,\ldots,k$) and $h_j$ ($j=1,\ldots,\ell$) are the edges adjacent as in Figure~\ref{fig:FL_minus_holes} to the quadrilateral of the two triangles adjacent to~$e$ and $\varepsilon$ is a closed path on $\Sigma$ with $\pi(\varepsilon)=\epsilon$. . 
  \item[iv)] \textbf{Both inner vertices are holes, negative sign:} \quad
  Consider an edge $e$ such that both vertices are holes as depicted on the right in Figure~\ref{fig:FL_minus_holes}. Denote by $\epsilon_1,\epsilon_2$ loops around the holes. Let $\calL_{- e}$ be the lamination obtained by applying Proposition~\ref{prop:FGalgo} to the tuple $n_e = -1$ and $n_j=0$ for all other edges. Then there is a loop decomposition $\{\wp_{-e},\varepsilon_1,\varepsilon_2\}$ of $\calL_{-e}$. Choose a directed base point $\wt{z}$ on $\wp_{-e}$ in the interior of the quadrilateral, such that the trivial lift of $\wp_{-e}$ starting at $\wt{z}$ yields $[\dnl\gamma_{e}]^{-1}[\dnl\varepsilon_1+\dnl\varepsilon_2]$, where $\varepsilon_i$ are closed paths on $\Sigma$ with $\pi(\varepsilon_i)=\epsilon_i$. Then 
  \bas
  F_{\bullet}(\wp_{-e},\theta) &\= [\dnl\gamma_{e}]^{-1}[{\dnl\varepsilon_1}+\dnl\varepsilon_2]
 + \sum_{i=1}^k [{\dnl\varepsilon_1}+\dnl\varepsilon_2](\mathfrak{c}+[\dnl\gamma_{e}]^{-1}])
 \prod_{m=1}^i[\dnl\gamma_{g_m}]^{-1} \\
  & \qquad + \mathfrak{c}\sum_{j=1}^\ell [\dnl\gamma_{e}]^{-1}[\dnl\varepsilon_1+\dnl\varepsilon_2] \prod_{n=1}^j[\dnl\gamma_{h_n}]^{-1} \\
  & \qquad + \mathfrak{c}\sum_{i,j=1}^{\ell,k}[\dnl\gamma_{e}]^{-1}[\dnl\varepsilon_1+\dnl\varepsilon_2] \prod_{m=1}^i [\dnl\gamma_{g_m}]^{-1}\prod_{n=1}^j[\dnl\gamma_{h_n}]^{-1} 
  \eas
where~$g_i$  ($i = 1,\ldots,k$) and $h_j$ ($j=1,\ldots,\ell$) are the edges adjacent as in Figure~\ref{fig:FL_minus_holes} to the quadrilateral of the two triangles adjacent to $e$ and $\mathfrak{c}$ is zero or one depending on the precise position of the directed base point $\wt{z}$.
  \end{itemize}
  \end{lemma}
Using these formulas we obtain an immediate generalization of Proposition~\ref{prop:approxFL}:
\par
\begin{proposition}\label{prop:approxFL_based}
  Suppose that $\theta = 0$ is the horizontal direction and the sign of the hat-homology classes is chosen as $Z(\gamma_e) \in \HH$ for all $e \in \TT$. Then for any $e_0 \in \TT$ and there are sequences $P^\pm_n$ of polynomials in the based lifts $F_\bullet(\wp_{+e},0)$, $F_\bullet(\wp_{-e}, 0)$ where $e \in \TT$ and loops around the holes, such that $\sum_{i=1}^n P_n^\pm \to [\wt\gamma_{e_0}]^{\pm 1}$ in the topology of $\wh{H}_\Delta^\abs$.
  \par
  Moreover, there exists a cone~$\Delta_2$ strictly contained in the lower half plane and~$t \in \CC$ with the following property: For any polynomial~$P = \sum_\alpha \Omega(\alpha) [\alpha] $ in $F_\bullet(\wp_{e},0)$ and $F(\wp_{-e},0)$ whose total degree in the lifts of positive generators $F(\wp_{+e},0)$ is at most one, all $\alpha$ with $\Omega(\alpha) \neq 0$ can be represented as a product $\alpha=[\alpha'][\epsilon]$, such that $Z(\alpha') \in \Delta_2^t$ and $[\epsilon]$ is product of loops around the holes. 
\end{proposition}

\section{BPS-invariants and wall-crossing} \label{sec:BPSWC}

The goal of this section is to complete the proof of the wall-crossing formula Theorem~\ref{thm:introWC}, first in the infinite area case. Later in Section~\ref{sec:FiniteArea} we digress on the finite area case, including a discussion about the content of the wall-crossing formula. Along the way, the two key steps are the convergence to the limiting path lifting rules in Proposition~\ref{prop:Fconverge} and the geometric verification in Proposition~\ref{prop:Kformula1} that the BPS-automorphism indeed compares the two one-sided limits.

\subsection{Convergence when varying $\theta$} \label{sec:varyingtheta}

We fix a direction~$\theta_0$, typically admitting saddle connections, and a path~$\wp$ or an $\calA_0$-lamination~$\calL$ disjoint from the zeros of~$q$.
\par
\begin{proposition} \label{prop:Fconverge}
For a fixed path~$\wp$ or $\calA_0$-lamination $\calL$ the path liftings converge
\ba
\lim_{\theta \to \theta_0^+} F(\wp,\theta) \= F^+(\wp,\theta_0) \quad
\text{and} \quad \lim_{\theta \to \theta_0^-} F(\wp,\theta) \= F^-(\wp,\theta_0), \\
\lim_{\theta \to \theta_0^+} F(\calL,\theta) \= F^+(\calL,\theta_0) \quad
\text{and} \quad \lim_{\theta \to \theta_0^-} F(\calL,\theta) \= F^-(\calL,\theta_0) \\
\ea
in the $Z$-adic topology.
\end{proposition} 
\par
\begin{lemma}\label{lem:approx_+_adm_rays} For every $+$-admissible ray~$\frakr$ ending at $y\in W_{\theta_0}(\wp)$ there exists $\e >0$ with the following property: For every direction $\theta' \in (\theta_0, \theta_0 + \e)$ there exists a ray~$\frakr' \in W_{\theta'}$ starting at the same zero as~$\frakr$, ending a point~$y' \in W_{\theta'}(\wp)$, and such that $\frakr$ composed with the segment~$[y,y'] \subset \wp_\pm$ is homotopic to~$\frakr'$. The same statement holds for minus-admissible rays and directions $\theta' \in (\theta_0-\e, \theta_0)$.
\par
Conversely, for each $L >0$, there exists $\epsilon >0$, such that for each ray~$\frakr' \in W_{\theta_0+\epsilon}$ with $Z(\frakr') <L$ ending at $y\in W_{\theta_0+\e}(\wp)$,  there exists a $+$-admissible ray $\frakr$ of $W_{\theta_0}(\wp)$ starting at the same zero and ending at $y' \in W_{\theta_0}(\wp)$, such that $\frakr$ composed with the segment~$[y,y'] \subset \wp_\pm$ is homotopic to~$\frakr'$. The same statements holds for minus-admissible rays and directions $\theta \in (\theta_0-\e, \theta_0)$. 
\end{lemma}
\par
This lemma holds verbatim with~$\wp$ replaced by a $\calA_0$-lamination~$\calL$.
\begin{figure}
  \resizebox{10cm}{!}{
%
%


\begin{tikzpicture}
\tikzset{
  on each segment/.style={
    decorate,
    decoration={
      show path construction,
      moveto code={},
      lineto code={
        \path [#1]
        (\tikzinputsegmentfirst) -- (\tikzinputsegmentlast);
      },
      curveto code={
        \path [#1] (\tikzinputsegmentfirst)
        .. controls
        (\tikzinputsegmentsupporta) and (\tikzinputsegmentsupportb)
        ..
        (\tikzinputsegmentlast);
      },
      closepath code={
        \path [#1]
        (\tikzinputsegmentfirst) -- (\tikzinputsegmentlast);
      },
    },
  },
  mid arrow/.style={postaction={decorate,decoration={
        markings,
        mark=at position .5 with {\arrow[#1]{stealth}}
      }}},
  threequarter arrow/.style={postaction={decorate,decoration={
        markings,
        mark=at position .75 with {\arrow[#1]{stealth}}
      }}},
}


\node[draw, circle, fill=black, inner sep=0.0pt] (A) at (0, 0) {}; 
\node[draw, circle, fill=black, inner sep=0.0pt] (B) at (2, 0) {};
\node[draw, circle, fill=black, inner sep=0.0pt] (N) at (3, 1) {};

\node (C) at (6, 0) {};
\node (D) at (-1, 0) {};
\node (E) at (0, -1) {};
\node (BC1) at (-1, -1) {};
\node (BC2) at (3, -1) {};
\node (BC3) at (2.5, 1.5) {};

\node (Eb) at (0-1, -1) {};
\node (F) at (2, -1) {};
\node[draw, circle, fill=black, inner sep=0.0pt] (N1) at (3.5, 1) {};
\node[draw, circle, fill=black, inner sep=0.0pt] (N2) at (2.5, 1) {};
\node[draw, circle, fill=black, inner sep=0.0pt] (N3) at (3, 1.5) {};

\node[draw, circle, fill=black, inner sep=0.0pt] (R1) at (-0.1, 0.05) {};
\node[draw, circle, fill=black, inner sep=0.0pt] (R2) at (6, 0.05) {};
\node[draw, circle, fill=black, inner sep=0.0pt] (R3) at (6, 0.95) {};
\node[draw, circle, fill=black, inner sep=0.0pt] (R4) at (-0.1, 0.95) {};

\node[draw, circle, fill=black, inner sep=0.0pt] (P1) at (5,-0.7) {};
\node[draw, circle, fill=black, inner sep=0.0pt] (P2) at (5,1.5) {};

\node[color=red, draw, circle, fill=black, inner sep=0.0pt] (Det1S) at (4.95,-0.1) {};
\node[color=red, draw, circle, fill=black, inner sep=0.0pt] (Det2S) at (4.95,0.95) {};

\node[color=blue, draw, circle, fill=black, inner sep=0.0pt] (WTprime) at (5.9,1.2)  {};

\draw[color=blue, postaction={on each segment={mid arrow=blue}}] (B) -- (A)  {};
\draw[color=blue, postaction={on each segment={mid arrow=blue}}] (A) -- (D)  {};
\draw[color=blue, postaction={on each segment={mid arrow=blue}}] (C) -- (B)  {};
\draw[color=blue, postaction={on each segment={mid arrow=blue}}] (A) -- (E)  {};
\draw[color=blue, postaction={on each segment={mid arrow=blue}}] (N1) -- (N)  {};
\draw[color=blue, postaction={on each segment={mid arrow=blue}}] (N) -- (N2)  {};
\draw[color=blue, postaction={on each segment={mid arrow=blue}}] (N) -- (N3)  {};
\draw[color=blue, postaction={on each segment={mid arrow=blue}}] (F) -- (B)  {};
\draw[color=dgreen, postaction={on each segment={mid arrow=dgreen}}] (P1) -- (P2)  node[midway, right] {$\wp$};
\draw[color=dpurple, postaction={on each segment={mid arrow=dpurple}}] (Det1S) -- (A)  {};
\draw[color=red, postaction={on each segment={mid arrow=red}}] (Det2S) -- (A)  {};
\draw[color=red, postaction={on each segment={threequarter arrow=red}}] (Det1S) -- (Det2S)  {};

\draw[color=blue, dashed] (WTprime) -- (A)  {};

\draw[color=black] (R1) -- (R2)  {};
\draw[color=black] (R3) -- (R2)  node[midway, right] {$\delta$};
\draw[color=black] (R3) -- (R4)  {};
\draw[color=black] (R1) -- (R4)  node[midway, left] {$\delta$};

\draw[decorate, decoration={zigzag, segment length=4, amplitude=1.5}] (BC1) -- (A); 
\draw[decorate, decoration={zigzag, segment length=4, amplitude=1.5}] (BC2) -- (B); 
\draw[decorate, decoration={zigzag, segment length=4, amplitude=1.5}] (BC3) -- (N);

\node (Ab) at (0, 0) {$\times$};
\node (Bb) at (2, 0) {$\times$};
\node (Nb) at (3, 1) {$\times$};
\node (WTpnode) at (6.5, 1.3) {$W_{\theta'}$};
\node (WTnode) at (6.5, 0) {$W_{\theta}$};

\end{tikzpicture}
    }
    \caption{A zero-free rectangle enclosing $+$-admissible rays approximating the horizontal ray. The purple and the red path are obviously homotopic.}
  \label{fig:approx_+-admissible_ray}
\end{figure}
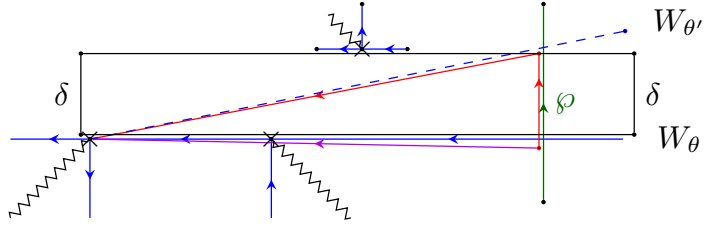
\par
\begin{proof} We claim that for every $\frakr$ there exists~$\delta>0$ such that
the union of segments in direction $\theta_0 + \pi/2$ of length~$\delta$ starting at points in~$\frakr$ does not contain any zero of~$q$. This follows from the compactness of~$\frakr$. The angle condition for turns at zero in the definition of $+$-admissible implies that the union of those $\delta$-segments is indeed  a metric rectangle in~$(X,q)$, see Figure~\ref{fig:approx_+-admissible_ray}. The diagonal of that rectangle is the ray~$\frakr'$ claimed in the first statement. Consider $\delta$-segments in the direction $\theta_0 - \pi/2$ for the minus-admissible case.
\par
For the converse take $\e$ so small that the sectors between the direction~$\theta_0$
and $\theta_0 + \e$ centered at all the zeros do not contains no saddle connection of length~$<L+|\wp|$. (This is possible by the finiteness of the set of saddle connections of bounded length.) Under this conditions for~$\e$ the counterclockwise rotation by the angle $\theta-\theta_0$ of the ray~$\frakr'$ homotopes it into the desired ray~$\frakr$.
\end{proof}
\par
\begin{proof}[Proof of Proposition~\ref{prop:Fconverge}] We deal with the case of a path~$\wp$. Since the proof is about the set of possible detours, which is define the same way for laminations, the conclusion for any~$\calL$ follows from the same proof.
\par
The strategy of proof is as follows: By Lemma \ref{lem:approx_+_adm_rays} we can identify elementary detours in $\bfF^+(\wp,\theta_0)$ with elementary detours in $F(\wp,\theta_0+\epsilon)$ after concatenation with a small segment of~$\wp$. The crux is to correctly determine the~$\epsilon$ for which we want to apply this. Moreover, we have to treat starting points of detours that are disjoint for $\theta_0 + \epsilon$ but collide in the limit direction~$\theta_0$.
\par
In detail we fix $L'>L_0 >0$, where $L_0$ is the length of the longest saddle trajectory in direction $\theta_0$ and intend to show that there is some~$\epsilon$ such that $\bfF(\wp,\theta') = \bfF^+(\wp,\theta_0)$ in the groupoid of path truncated at length $\geq L'$, i.e.\ that the two lifts are within distance $2^{-L'}$. Let $L = L' + |Z(\wp)|$ and define $\bfF_{L}(\wp,\theta)$ and $\bfF_{L}^+(\wp,\theta_0)$ to be the lift of $\wp$ using the truncated spectral network $W_{\theta,L}$ resp. $W_{\theta_0,L}$ of critical trajectories of length at most~$L$.
\par
Choose $\epsilon>0$, such that the second statement of Lemma~\ref{lem:approx_+_adm_rays} holds. By choosing $\epsilon$ even smaller, we may assume that for all $\theta \in (\theta_0,\theta_0+\epsilon]$ the truncated spectral network $W_{\theta,L}$ has no saddle connections. (There are only finitely many saddle connections for varying $\theta$ of length~$<L$. Thus we can pick $\e$ small enough, so that none of them is a trajectory in direction $\theta \in [\theta_0,\theta_0+\e]$). Furthermore, we may assume that no $W_{\theta,L}$ contains the endpoints of $\wp$ for $\theta \in [\theta_0,\theta+\e]$. 
\begin{figure}
  \resizebox{12cm}{!}{
%
%


\begin{tikzpicture}
\tikzset{
  on each segment/.style={
    decorate,
    decoration={
      show path construction,
      moveto code={},
      lineto code={
        \path [#1]
        (\tikzinputsegmentfirst) -- (\tikzinputsegmentlast);
      },
      curveto code={
        \path [#1] (\tikzinputsegmentfirst)
        .. controls
        (\tikzinputsegmentsupporta) and (\tikzinputsegmentsupportb)
        ..
        (\tikzinputsegmentlast);
      },
      closepath code={
        \path [#1]
        (\tikzinputsegmentfirst) -- (\tikzinputsegmentlast);
      },
    },
  },
  mid arrow/.style={postaction={decorate,decoration={
        markings,
        mark=at position .5 with {\arrow[#1]{stealth}}
      }}},
onequarter arrow/.style={postaction={decorate,decoration={
        markings,
        mark=at position .25 with {\arrow[#1]{stealth}}
      }}},
}
  

\node[draw, circle, fill=black, inner sep=0.0pt] (A) at (0, 0) {}; 
\node[draw, circle, fill=black, inner sep=0.0pt] (B) at (2, 0) {};
\node[draw, circle, fill=black, inner sep=0.0pt] (As) at (0, -0.3) {}; 
\node[draw, circle, fill=black, inner sep=0.0pt] (Bs) at (2, 0.3) {};
\node (Cp) at (3, 1) {};
\node (Cps) at (3, 1.3) {};
\node (C0) at (3, 0) {};
\node (Cm) at (3, -1) {};
\node (Dp) at (-1, 1) {};
\node (D0) at (-1, 0) {};
\node (Dm) at (-1, -1) {};
\node (Dms) at (-1, -1.3) {};

\node[draw, circle, fill=black, inner sep=0.0pt] (AA) at (6, 0) {}; 
\node[draw, circle, fill=black, inner sep=0.0pt] (BB) at (8, 0) {};
\node[draw, circle, fill=black, inner sep=0.0pt] (AAs) at (6, -0.3) {}; 
\node[draw, circle, fill=black, inner sep=0.0pt] (BBs) at (8, 0.3) {};
\node (CCp) at (9, 1) {};
\node (CCps) at (9, 1.3) {};
\node (CC0) at (9, 0) {};
\node (CCm) at (9, -1) {};
\node (DDp) at (5, 1) {};
\node (DD0) at (5, 0) {};
\node (DDm) at (5, -1) {};
\node (DDms) at (5, -1.3) {};

\node[color=dgreen, draw, circle, fill=dgreen, inner sep=0.0pt](Q1) at (7, 0.7) {};
\node[color=dgreen, draw, circle, fill=dgreen, inner sep=0.0pt](Q2) at (7, -0.5) {};
\node[color=dgreen, draw, circle, fill=dgreen, inner sep=0.0pt](Qs1) at (1, 0.7) {};
\node[color=dgreen, draw, circle, fill=dgreen, inner sep=0.0pt](Qs2) at (1, -0.5) {};

\node[color=dpurple, draw, circle, fill=dpurple, inner sep=0.0pt] (QD1) at (7, 0.31) {};
\node[color=dpurple, draw, circle, fill=dpurple, inner sep=0.0pt] (QD2) at (5.84, 0.16) {};
\node[color=dpurple, draw, circle, fill=dpurple, inner sep=0.0pt] (QD3) at (5.84, 0.0) {};
\node[color=dpurple, draw, circle, fill=dpurple, inner sep=0.0pt] (QD4) at (-0.16, 0.0) {};
\node[color=dpurple, draw, circle, fill=dpurple, inner sep=0.0pt] (QD5) at (0, -0.16) {};
\node[color=dpurple, draw, circle, fill=dpurple, inner sep=0.0pt] (QD6) at (2.16, 0.15) {};
\node[color=dpurple, draw, circle, fill=dpurple, inner sep=0.0pt] (QD7) at (2.16, 0.00) {};
\node[color=dpurple, draw, circle, fill=dpurple, inner sep=0.0pt] (QD8) at (8.16, 0.00) {};
\node[color=dpurple, draw, circle, fill=dpurple, inner sep=0.0pt] (QD9) at (8.16, -0.16) {};
\node[color=dpurple, draw, circle, fill=dpurple, inner sep=0.0pt] (QD10) at (7.0, -0.31) {};

%

\draw[color=blue, postaction={on each segment={mid arrow=blue}}] (A) -- (Bs)  {};
\draw[color=blue] (Bs) -- (Cps)  {};
\draw[color=blue, postaction={on each segment={mid arrow=blue}}] (As) -- (B)  {};
\draw[color=blue] (As) -- (Dms)  {};

\draw[color=blue, postaction={on each segment={mid arrow=blue}}] (B) -- (Cp)  {};
\draw[color=blue, postaction={on each segment={mid arrow=blue}}] (B) -- (Cm)  {};
\draw[color=blue, postaction={on each segment={mid arrow=blue}}] (Dp) -- (A)  {};
\draw[color=blue, postaction={on each segment={mid arrow=blue}}] (Dm) -- (A)  {};

\draw[color=blue, postaction={on each segment={mid arrow=blue}}] (BBs) -- (AA)  {};
\draw[color=blue] (CCps) -- (BBs)  {};
\draw[color=blue, postaction={on each segment={mid arrow=blue}}] (BB) -- (AAs)  {};
\draw[color=blue] (DDms) -- (AAs)  {};
\draw[color=blue, postaction={on each segment={mid arrow=blue}}] (CCp) -- (BB)  {};
\draw[color=blue, postaction={on each segment={mid arrow=blue}}] (CCm) -- (BB)  {};
\draw[color=blue, postaction={on each segment={mid arrow=blue}}] (AA) -- (DDp)  {};
\draw[color=blue, postaction={on each segment={mid arrow=blue}}] (AA) -- (DDm)  {};

\draw[color=dgreen, postaction={on each segment={onequarter arrow=dgreen}}] (Q1) -- (Q2)  node[near start, right] {$\scriptstyle{\wq_1}$};
\draw[color=dgreen, postaction={on each segment={onequarter arrow=dgreen}}] (Qs1) -- (Qs2)  node[near start, right] {$\scriptstyle{\wq_2}$};

\draw[color=dpurple, postaction={on each segment={mid arrow=dpurple}}] (QD1) -- (QD2)  {};
\draw[color=dpurple] (QD2) -- (QD3)  {};
\draw[color=dpurple] (QD4) -- (QD5)  {};
\draw[color=dpurple, postaction={on each segment={onequarter arrow=dpurple}}] (QD5) -- (QD6)  {};
\draw[color=dpurple] (QD6) -- (QD7)  {};
\draw[color=dpurple] (QD8) -- (QD9)  {};
\draw[color=dpurple] (QD9) -- (QD10)  {};

%

\draw[decorate, decoration={zigzag, segment length=4, amplitude=1.5}] (C0) -- (B); 
\draw[decorate, decoration={zigzag, segment length=4, amplitude=1.5}] (A) -- (D0); 

\node (Ab) at (0, 0) {$\times$};
\node (Bb) at (2, 0) {$\times$};

\draw[decorate, decoration={zigzag, segment length=4, amplitude=1.5}] (CC0) -- (BB); 
\draw[decorate, decoration={zigzag, segment length=4, amplitude=1.5}] (AA) -- (DD0); 

\node (AAb) at (0, 0) {$\times$};
\node (BBb) at (2, 0) {$\times$};

\end{tikzpicture}
    }
\caption{Two consecutive detours (in purple) at two points of the spectral network. In the limit as $\theta \to \theta_0$ they start  over the same point of~$X$. The limit is the detour with American driving rule from Figure~\ref{fig:AmBritDetours}.
}
\label{fig:AlmostAmerican}
\end{figure}
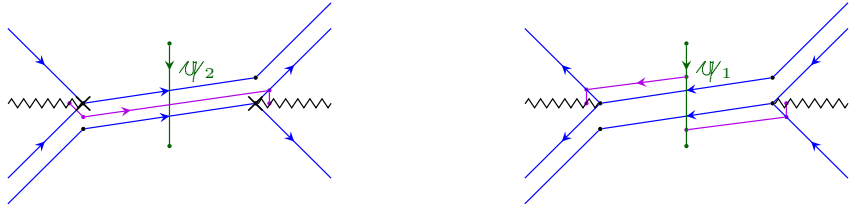
\par
We aim to restrict~$\epsilon$ once again to ensure that Lemma~\ref{lem:approx_+_adm_rays} defines a bijection among detours of length~$<L$ for all directions~$\theta \in [\theta_0,\theta_0+\epsilon]$. For such~$\theta$ consider a detour $\nl{D}_{\bfy}^{\theta}$ following rays $\frakr_1^{\theta},\dots, \frakr_N^{\theta}$. By the second part of Lemma~\ref{lem:approx_+_adm_rays} there are $+$-admissible rays $\frakr_i^{\theta_0}$ in direction $\theta_0$, such that $\frakr_i^{\theta}$ is given by the concatenation of $\frakr_i^{\theta_0}$ with a segment $\wp_i=[y_i,y_i'] \in \wp$ up to homotopy. For each~$i$ this segment depends on~$\theta$. In fact it is just one point for $\theta = \theta_0$ and the length is bounded by $L \sin(\theta-\theta_0)$. Since the truncated network~$W_{\theta_0,L}$ has only finitely many intersection points with~$\wp$, we may choose $\epsilon$ small enough so that all these segments $\wp_i$ are disjoint for all detours. Now
\begin{align} \dnl{D}_{\bfy}^{\theta} \=\dnl{D}_{\frakr_1^{\theta}, \dots, \frakr_N^{\theta}}=\dnl{D}_{\frakr_1^{\theta_0},\dots, \frakr_N^{\theta_0}}\label{equ:detours_equal}
\end{align} in $\wt{\Pi}_{\leq 1}(y_s,y_e)$, where $y_s,y_e$ are the start and endpoint of $\wp$. Moreover, if two rays $\frakr_i^{\theta_0}, \frakr_{i+1}^{\theta_0}$ have endpoints interchanged by $\sigma$ the detour $\nl{D}_{\frakr_1^{\theta_0},\dots, \frakr_N^{\theta_0}}$ is clearly American (see Figure~\ref{fig:AlmostAmerican} for an example). This shows that all terms of $\bfF_L(\wp,\theta)$ appear also in $\bfF^+_L(\wp,\theta_0)$.
\par
The converse inclusion holds for terms $\bfF^+_L(\wp,\theta_0)$ whose total length is at most~$L'$. In fact, if the term consists of a single detour along~$\frakr$, then by the first part of Lemma~\ref{lem:approx_+_adm_rays} there is a ray~$\frakr'$ in the direction~$\theta$ such that the detour composed with a piece of~$\wp$ is homotopic to the detour along~$\frakr$. By the choice of~$L'$, the ray belongs to $W_{\theta,L}$. (If there are several detours, for each of them application of the lemma requires a piece of~$\wp$ to get a homotopic path, so the length difference between the two rays is even shorter). This implies the claim. (Note that $\bfF^+_L(\wp,\theta_0)$ may contain very long detour paths by following consecutively a saddle connection many times. Not all of them will be visible in $\bfF_L(\wp,\theta)$ for a fixed~$L$, but they will be taken into account in the limit as $L \to \infty$, as claimed.) 
\end{proof}

\subsection{ Crossing saddle connection directions of rank one}
\label{sec:rankonedir}

In this section we work under the hypothesis that the direction~$\theta$ on
$(X,q)$ has a saddle connection~$\gamma$ and moreover the span in $\Gamma$
of all hat-homology classes of saddle connections in direction~$\theta$ is a $\ZZ$-module of rank one.\footnote{In the language of \cite{GMN_spectral}, the directed surface $(X,q,\theta)$ belongs to a \emph{$\calK$-wall.}} We call these directions \emph{rank one directions}\footnote{In \cite{KS08} or \cite{BridgelandSmith} they are called \emph{generic}. We try to avoid ambiguities stemming from this terminology.} and we fix the generator~$\wh{\gamma}_0 \in \Gamma$, pointing in the direction of~$\theta$, i.e., such that $\Re(Z(e^{-i\theta} \wh{\gamma}_0))> 0$.
\par
Our goal in this section is to define BPS-invariants of rank one directions, which enter into the definition of the comparison isomorphism~$\calK_\theta$ that relates
the one-sided limits of path lifting function $F^\pm(\wp,\theta)$ defined in Proposition~\ref{prop:Fpm}.
\par
\begin{theorem} \label{thm:WCforF}
For any half-translation surface $(X,q) \in \calQ_q(\mu)$, for any rank one direction~$\theta$, and any path~$\wp$ or any $\calA_0$-lamination~$\calL$
\be \label{eq:FpisKFm}
F^+(\wp,\theta) \= \calK_\theta F^-(\wp,\theta),\qquad 
F^+(\calL,\theta) \= \calK_\theta F^-(\calL,\theta)\,.
\ee
\end{theorem}
\par
Together with Corollary~\ref{cor:composition}, Theorem~\ref{theo:homotopy_invariance}, and Proposition~\ref{prop:Fconverge} the preceding Theorem~\ref{thm:WCforF} completes the \emph{proof of Theorem~\ref{thm:introF}.}
\par
We start by reviewing rank one directions. These directions have been classified in \cite[Lemma~5.1]{BridgelandSmith} in the infinite area case and in \cite[Proposition~5.4]{Haiden} in the finite area case. We follow the list of \cite[Proposition~6.1]{KidWil} that includes both, as well as the case of simple poles. See \cite[Figure~9]{KidWil} for pictures of all the cases. 
\par
\begin{lemma} \label{le:r1class}
A rank one direction on $(X,q)$ has one of the following possibilities: 
\begin{itemize}
\item[(1)]  A single (non-closed) saddle connection~$\gamma_0$ of type I, i.e., joining two different simple zeros.
\item[(2)]  A single (non-closed) saddle connection~$\gamma_0$ of type II, i.e., joining a simple zeros with a simple pole.
\item[(3)]  A degenerate ring domain, with one boundary given by a double pole and the other by either
\begin{itemize}
\item[(a)] a closed saddle trajectory, 
\item[(b)] in case $\mu = (1,1,-2)$ two type I saddle trajectories, which on their other side form the boundary of a torus whose interior is a spiral domain, or
\item[(c)] in case $\mu = (-1,-1,-2)$ a type III saddle trajectory, i.e., joining two different simple poles.
\end{itemize}
\item[(4)] A non-degenerate ring domain, with boundaries given by either
\begin{itemize}
\item[(a)]  a closed saddle trajectory for each boundary, or
\item[(b)]  a closed saddle trajectory for one boundary, and for the other boundary two type I saddle trajectories, which on their other side form the boundary of a torus whose interior is a spiral domain, or
\item[(c)] a closed saddle trajectory for one boundary, and a type III saddle trajectory for the other
boundary.
\end{itemize}
\end{itemize}
\end{lemma}
\par
The ring domains are subdivided into the \emph{standard (a)}, the \emph{toral (b)} and the \emph{type III} ring domains. Precisely the toral ring domains are rank one directions that have a spiral domain. In the absence of simple poles only the cases (1), (3a), (3b), (4a) and (4b) occur.
\par
\medskip
\paragraph{\textbf{BPS-cycles, the BPS-automorphism, and BPS invariants}}  For each of the cases and to each of the multiples of $\wh{\gamma}_0$ that represent a type~I saddle connection or a ring domain we associate a closed \emph{BPS-cycle} $L(\wh{\gamma}_0)$ (not just a homology class!) as follows.
\ba \label{eq:BPSL}
\text{Case 1}: \qquad   & \qquad L(\wh{\gamma}_0) \= \wh{\gamma}_0\,. \\
\text{Case 3a}: \qquad   & \qquad L(\wh{\gamma}_0) \= - \wh{\gamma}_b\,. \\
\text{Case 3b}: \qquad   & \qquad L(\wh{\gamma}_0) \= + \wh{\gamma}_0 + \wh{\gamma}_1\,. \\
 & \qquad L(2\wh{\gamma}_0) \= - \wh{\gamma}_0 - \wh{\gamma}_1\,. \\
\text{Case 4a}: \qquad   & \qquad L(\wh{\gamma}_0) \= - \wh{\gamma}_t - \wh{\gamma}_b\,. \\
\text{Case 4b}: \qquad   & \qquad L(\wh{\gamma}_0) \= + \wh{\gamma}_0 + \wh{\gamma}_1\,. \\
 & \qquad L(2\wh{\gamma}_0) \= - \wh{\gamma}_0 - \wh{\gamma}_1  - \wh{\gamma}_2 \,. \\
\ea
Here in Case~(4a) $\gamma_t$ and~$\gamma_b$ are the saddle connections at the two ends (top and bottom) of the cylinder, oriented the same way as~$\wh{\gamma}_0$, i.e.\ oriented in the direction of the spectral network. In Case~(4b) the saddle connections are labeled as in~\eqref{eq:toralendZ} and Figure~\ref{fig:KforToralEnd} below. In case Case~(3a) we picture the degenerate ring domain as having a bottom saddle connection and extending infinitely towards the top. Similarly, we picture Case~(3b) as a limit of Case~(4b) where the saddle connection~$\gamma_2$ disappears off to infinity in the ring domain.
\par
We can now define the \emph{BPS-automorphism} on generators $\alpha \in \Mor_{H_{\leq 1}}(\wt{\Sigma}^\circ)$ by
\be \label{eq:KbyL}
\calK_\theta([\alpha]) \=  \prod_{\gamma \in \NN \wh{\gamma}_0}  (1-[\wt{\gamma}])^{\langle L(\gamma), c(\alpha) \rangle} [\alpha] \,.
\ee
Here $c(\alpha)$ is the image of~$\alpha$ in~$\Sigma$ and the intersection product is the one on~$\Sigma$. The multiplication of path classes was defined in~\eqref{eq:closedpathmult}. Observe that two paths that differ by an element in~$\ker(c)$ have the \emph{same} intersection product with~$L(\wh{\gamma}_0)$. So even if the classes of these paths are negatives of each other in the quotient by the winding ideal, the action by multiplication in~\eqref{eq:KbyL} is well-defined. By Lemma~\ref{le:Kgeneralcont} the definition on generators gives a well-defined map $\calK_\theta: \wh{H}^\abs_\Delta \to \wh{H}^\abs_\Delta$, which is continuous on the subspaces $\wh{H}^\abs_{\Delta^K}$ as defined in loc.~cit.
\par
\medskip
For Case~(1) and Case~(4b) the cycles  $L(\wh{\gamma}_0)$ are given explicitly in \cite{GMN_spectral}, too. As in loc.~cit. we need to use cycles since in the proofs we intersect with homology classes of non-closed paths. The main result describes the effect of the BPS-automorphism on closed paths. In this case we can replace in the definition~\ref{eq:KbyL} the cycles~$L(\wh{\gamma}_0)$ by its class. Then 
\be \label{eq:KbyOmega}
\calK_\theta([\alpha]) \=  \prod_{\gamma \in \NN \wh{\gamma}_0}  (1-[\gamma])^{\Omega(\gamma) \langle \wh{\gamma}_0, c(\alpha) \rangle} [\alpha] \,.
\ee
where the \emph{BPS-invariants} $\Omega(\gamma) = [L(\gamma)]/[\wh{\gamma}_0]$ are given by
\ba
\text{Case 1}: \qquad   & \qquad \Omega(\wh{\gamma}_0) \= +1\,. \\
\text{Case 3a}: \qquad   & \qquad \Omega(\wh{\gamma}_0) \= -1\,. \\
\text{Case 3b}: \qquad   & \qquad \Omega(\wh{\gamma}_0) \= +2,  \qquad \Omega(2\wh{\gamma}_0) \= -1\,. \\
\text{Case 4a}: \qquad   & \qquad \Omega(\wh{\gamma}_0) \= -2\,. \\
\text{Case 4b}: \qquad   & \qquad \Omega(\wh{\gamma}_0) \= +2,  \qquad \Omega(2\wh{\gamma}_0) \= -2\,. \\
\ea
Equivalently, $\Omega(\gamma)$ is the number of type~I saddle connections of class~$\gamma$ minus two times the number of ring domain (core curves) of class~$\gamma$, in agreement with the literature (dating back to \cite{KS08}, see also \cite[Theorem~1.4]{BridgelandSmith}, see also  \cite[Definition~2.3]{Allegretti} or \cite{Haiden} and \cite{KidWil} for discussions also including simple poles).
\par
\medskip
In preparation for the proof of Theorem~\ref{thm:WCforF} the following lemmas compare the possible $+$-admissible and minus-admissible detours in all rank one situations. Fix a path fragment $\wp$ intersecting the (truncated) spectral network exactly once. We call the \emph{shortest (elementary) detour~$D_y^\sh$} at $y \in W_\theta$ the elementary detour that starts at~$y$, runs to the nearest singularity and returns on the other sheet (without following any saddle connections). Clearly, the shortest detour is always $+$ and minus-admissible. To fix signs, we suppose $\theta=0$ is the horizontal direction, such that $Z(\wh{\gamma}_0) < 0$.  We denote by $[\wt{\gamma}_0] \in \wh{H}^{\abs}$ the canonical lift of $[\wh{\gamma}_0] \in H_1(\Sigma\setminus\calP(\lambda),\mathbb{Z})^-$. In the following discussion we fix arbitrarily a direction on the path fragments, which will influence intersection numbers. The opposite orientation does not present new cases, as the orientation change can be compensated by swapping the role of the sheets.
\par
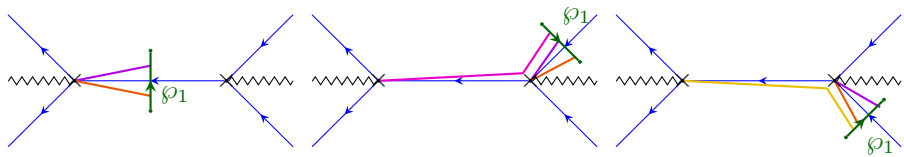
\begin{figure}
%
%


\begin{tikzpicture}
\tikzset{
  on each segment/.style={
    decorate,
    decoration={
      show path construction,
      moveto code={},
      lineto code={
        \path [#1]
        (\tikzinputsegmentfirst) -- (\tikzinputsegmentlast);
      },
      curveto code={
        \path [#1] (\tikzinputsegmentfirst)
        .. controls
        (\tikzinputsegmentsupporta) and (\tikzinputsegmentsupportb)
        ..
        (\tikzinputsegmentlast);
      },
      closepath code={
        \path [#1]
        (\tikzinputsegmentfirst) -- (\tikzinputsegmentlast);
      },
    },
  },
  mid arrow/.style={postaction={decorate,decoration={
        markings,
        mark=at position .5 with {\arrow[#1]{stealth}}
      }}},
}
  

\node[draw, circle, fill=black, inner sep=0.0pt] (A) at (0, 0) {}; 
\node[draw, circle, fill=black, inner sep=0.0pt] (B) at (2, 0) {};
\node (Cp) at (3, 1) {};
\node (C0) at (3, 0) {};
\node (Cm) at (3, -1) {};
\node (Dp) at (-1, 1) {};
\node (D0) at (-1, 0) {};
\node (Dm) at (-1, -1) {};

\node[draw, circle, fill=black, inner sep=0.0pt] (AA) at (4, 0) {}; 
\node[draw, circle, fill=black, inner sep=0.0pt] (BB) at (6, 0) {};
\node (CCp) at (7, 1) {};
\node (CC0) at (7, 0) {};
\node (CCm) at (7, -1) {};
\node (DDp) at (3, 1) {};
\node (DD0) at (3, 0) {};
\node (DDm) at (3, -1) {};

\node[draw, circle, fill=black, inner sep=0.0pt] (AAA) at (8, 0) {}; 
\node[draw, circle, fill=black, inner sep=0.0pt] (BBB) at (10, 0) {};
\node (CCCp) at (11, 1) {};
\node (CCC0) at (11, 0) {};
\node (CCCm) at (11, -1) {};
\node (DDDp) at (7, 1) {};
\node (DDD0) at (7, 0) {};
\node (DDDm) at (7, -1) {};

\node[color=dgreen, draw, circle, fill=dgreen, inner sep=0.0pt](P1) at (6.15, 0.75) {};
\node[color=dgreen, draw, circle, fill=dgreen, inner sep=0.0pt](P2) at (6.65, 0.25) {};

\node[color=dgreen, draw, circle, fill=dgreen, inner sep=0.0pt](Qs1) at (1, 0.4) {};
\node[color=dgreen, draw, circle, fill=dgreen, inner sep=0.0pt](Qs2) at (1, -0.4) {};

\node[color=dgreen, draw, circle, fill=dgreen, inner sep=0.0pt](T1) at (10.15, -0.75) {};
\node[color=dgreen, draw, circle, fill=dgreen, inner sep=0.0pt](T2) at (10.65, -0.25) {};

\draw[color=blue, postaction={on each segment={mid arrow=blue}}] (B) -- (A)  {};
\draw[color=blue, postaction={on each segment={mid arrow=blue}}] (Cp) -- (B)  {};
\draw[color=blue, postaction={on each segment={mid arrow=blue}}] (Cm) -- (B)  {};
\draw[color=blue, postaction={on each segment={mid arrow=blue}}] (A) -- (Dp)  {};
\draw[color=blue, postaction={on each segment={mid arrow=blue}}] (A) -- (Dm)  {};

\draw[color=blue, postaction={on each segment={mid arrow=blue}}] (BB) -- (AA)  {};
\draw[color=blue, postaction={on each segment={mid arrow=blue}}] (CCp) -- (BB)  {};
\draw[color=blue, postaction={on each segment={mid arrow=blue}}] (CCm) -- (BB)  {};
\draw[color=blue, postaction={on each segment={mid arrow=blue}}] (AA) -- (DDp)  {};
\draw[color=blue, postaction={on each segment={mid arrow=blue}}] (AA) -- (DDm)  {};

\draw[color=blue, postaction={on each segment={mid arrow=blue}}] (BBB) -- (AAA)  {};
\draw[color=blue, postaction={on each segment={mid arrow=blue}}] (CCCp) -- (BBB)  {};
\draw[color=blue, postaction={on each segment={mid arrow=blue}}] (CCCm) -- (BBB)  {};
\draw[color=blue, postaction={on each segment={mid arrow=blue}}] (AAA) -- (DDDp)  {};
\draw[color=blue, postaction={on each segment={mid arrow=blue}}] (AAA) -- (DDDm)  {};

\draw[thick,color=dpurple] (1,+0.2) -- (A)  {} ;
\draw[thick,color=dorange] (1,-0.2) -- (A)  {} ;

\draw[thick,color=dpurple] (6.37,+0.53) -- (BB)  {} ;
\draw[thick,color=dpink] (6.25,+0.63) -- (5.9,0.1)  {} ;
\draw[thick,color=dpink] (5.9,+0.1) -- (AA)  {} ;

\draw[thick,color=dorange] (6.6,+0.32) -- (BB)  {} ;

\draw[thick,color=dorange] (10.32,-0.59) -- (BBB)  {} ;
\draw[thick,color=dyellow] (10.25,-0.65) -- (9.9,-0.1)  {} ;
\draw[thick,color=dyellow] (9.9,-0.1) -- (AAA)  {} ;

\draw[thick,color=dpurple] (10.6,-0.35) -- (BBB)  {} ;

\draw[thick, color=dgreen, postaction={on each segment={mid arrow=dgreen}}] (P1) -- (P2)  node[near start, above right] {$\wp_1$};

\draw[thick, color=dgreen, postaction={on each segment={mid arrow=dgreen}}] (Qs2) -- (Qs1)  node[near start, right] {$\wp_1$};

\draw[thick, color=dgreen, postaction={on each segment={mid arrow=dgreen}}] (T1) -- (T2)  node[near start, below right] {$\wp_1$};

\draw[decorate, decoration={zigzag, segment length=4, amplitude=1.5}] (C0) -- (B); 
\draw[decorate, decoration={zigzag, segment length=4, amplitude=1.5}] (A) -- (D0); 

\node (Ab) at (0, 0) {$\times$};
\node (Bb) at (2, 0) {$\times$};

\draw[decorate, decoration={zigzag, segment length=4, amplitude=1.5}] (CC0) -- (BB); 
\draw[decorate, decoration={zigzag, segment length=4, amplitude=1.5}] (AA) -- (DD0); 

\node (AAb) at (4, 0) {$\times$};
\node (BBb) at (6, 0) {$\times$};

\draw[decorate, decoration={zigzag, segment length=4, amplitude=1.5}] (CCC0) -- (BBB); 
\draw[decorate, decoration={zigzag, segment length=4, amplitude=1.5}] (AAA) -- (DDD0); 

\node (AAAb) at (8, 0) {$\times$};
\node (BBAb) at (10, 0) {$\times$};

\end{tikzpicture}
  \caption{A simple saddle connection with a path $\wp$ crossing various the spectral network at various rays: $+$-admissible (drawn in slope~$\epsilon$) and minus-admissible detours (in slope~$-\epsilon$) detours. Only one of the two sheets is drawn for each path.}
  \label{fig:KforSC}
\end{figure}

\begin{figure}
%
%


\begin{tikzpicture}
\tikzset{
  on each segment/.style={
    decorate,
    decoration={
      show path construction,
      moveto code={},
      lineto code={
        \path [#1]
        (\tikzinputsegmentfirst) -- (\tikzinputsegmentlast);
      },
      curveto code={
        \path [#1] (\tikzinputsegmentfirst)
        .. controls
        (\tikzinputsegmentsupporta) and (\tikzinputsegmentsupportb)
        ..
        (\tikzinputsegmentlast);
      },
      closepath code={
        \path [#1]
        (\tikzinputsegmentfirst) -- (\tikzinputsegmentlast);
      },
    },
  },
  mid arrow/.style={postaction={decorate,decoration={
        markings,
        mark=at position .5 with {\arrow[#1]{stealth}}
      }}},
}


\node[draw, circle, fill=black, inner sep=0.0pt] (A) at (0, 0) {}; 
\node[inner sep=0.0pt] (A1) at (-0.9, 0) {}; 
\node[draw, circle, fill=black, inner sep=0.0pt] (C) at (3, 0) {};
\node[inner sep=0.0pt] (C1) at (3.9, 0) {};
\node[draw, circle, fill=black, inner sep=0.0pt] (D) at (0, 1.4) {};
\node[inner sep=0.0pt] (D1) at (-0.9, 1.4) {};
\node[draw, circle, fill=black, inner sep=0.0pt] (E) at (3, 1.4) {};
\node[inner sep=0.0pt] (E1) at (3.9, 1.4) {};

\node (BC1) at (0,2.3) {};
\node (BC2) at (0, -0.9) {};
\node[draw, circle, fill=black, inner sep=0.0pt] (Q1) at (3.6,-0.4) {};
\node[inner sep=0.0pt] (Q2) at (3.6,0) {};
\node[draw, circle, fill=black, inner sep=0.0pt] (Q3) at (-0.6,-0.4) {};
\node[inner sep=0.0pt] (Q4) at (-0.6,0) {};

\draw[thick,color=black!70] (A) -- (A1);
\draw[thick,color=black!70] (D) -- (D1);
\draw[thick,color=black!70] (E) -- (E1);
\draw[thick,color=black!70] (C) -- (C1);
\draw[thick,color=black] (A) -- (D)  node[midway, left] {$\scriptstyle{a_1}$} ;
\draw[thick,color=black] (C) -- (E)  node[midway, right] {$\scriptstyle{a_1}$};


\draw[color=blue, postaction={on each segment={mid arrow=blue}}] (A) -- (C)  {};
\draw[color=blue, postaction={on each segment={mid arrow=blue}}] (D) -- (E)  {};
\draw[color=blue, postaction={on each segment={mid arrow=blue}}] (D) -- (D1)   {} node[near end, below,color=black] {$\scriptstyle{e_1}$};
\draw[color=blue, postaction={on each segment={mid arrow=blue}}] (A) -- (A1)   {} node[near end, above, color=black] {$\scriptstyle{f_1}$};
\draw[color=blue, postaction={on each segment={mid arrow=blue}}] (E) -- (E1)  node [near end, below,color=black] {$\scriptstyle{e_1}$};
\draw[color=blue, postaction={on each segment={mid arrow=blue}}] (C) -- (C1)  node [near end, above,color=black] {$\scriptstyle{f_1}$};

\draw[thick,color=dgreen, postaction={on each segment={mid arrow=dgreen}}] (Q2) -- (Q1)  node[near end, right] {$\scriptstyle{\wp_1}$};
\draw[thick,color=dgreen, postaction={on each segment={mid arrow=dgreen}}] (Q3) -- (Q4)  node[near start, left] {$\scriptstyle{\wp_1}$};

\draw[decorate, decoration={zigzag, segment length=4, amplitude=1.5}] (BC1) -- (D); 
\draw[decorate, decoration={zigzag, segment length=4, amplitude=1.5}] (BC2) -- (A); 

\draw[thick,color=dpurple] (0,-0.15) -- (C)  {} ;

%


\node[draw, circle, fill=black, inner sep=0.0pt] (Ab) at (6, 0) {}; 
\node[inner sep=0.0pt] (Ab1) at (5.1, 0) {}; 
\node[draw, circle, fill=black, inner sep=0.0pt] (Cb) at (9, 0) {};
\node[inner sep=0.0pt] (Cb1) at (9.9, 0) {};
\node[draw, circle, fill=black, inner sep=0.0pt] (Db) at (6, 1.4) {};
\node[inner sep=0.0pt] (Db1) at (5.1, 1.4) {};
\node[draw, circle, fill=black, inner sep=0.0pt] (Eb) at (9, 1.4) {};
\node[inner sep=0.0pt] (Eb1) at (9.9, 1.4) {};

\node (BCb1) at (6,2.3) {};
\node (BCb2) at (6, -0.9) {};
\node[draw, circle, fill=black, inner sep=0.0pt] (Qb1) at (9.6,-0.4) {};
\node[inner sep=0.0pt] (Qb2) at (9.6,0) {};
\node[draw, circle, fill=black, inner sep=0.0pt] (Qb3) at (5.4,-0.4) {};
\node[inner sep=0.0pt] (Qb4) at (5.4, 0) {};

\draw[thick,color=black!70] (Ab) -- (Ab1);
\draw[thick,color=black!70] (Db) -- (Db1);
\draw[thick,color=black!70] (Eb) -- (Eb1);
\draw[thick,color=black!70] (Cb) -- (Cb1);
\draw[thick,color=black] (Ab) -- (Db)  node[midway, left] {$\scriptstyle{a_2}$} ;
\draw[thick,color=black] (Cb) -- (Eb)  node[midway, right] {$\scriptstyle{a_2}$};


\draw[color=blue, postaction={on each segment={mid arrow=blue}}] (Cb) -- (Ab)  {};
\draw[color=blue, postaction={on each segment={mid arrow=blue}}] (Eb) -- (Db)  {};
\draw[color=blue, postaction={on each segment={mid arrow=blue}}] (Db1) -- (Db)   {} node[near start, below,color=black] {$\scriptstyle{e_2}$};
\draw[color=blue, postaction={on each segment={mid arrow=blue}}] (Ab1) -- (Ab)   {} node[near start, above, color=black] {$\scriptstyle{f_2}$};
\draw[color=blue, postaction={on each segment={mid arrow=blue}}] (Eb1) -- (Eb)  node [near start, below,color=black] {$\scriptstyle{e_2}$};
\draw[color=blue, postaction={on each segment={mid arrow=blue}}] (Cb1) -- (Cb)  node [near start, above, color=black] {$\scriptstyle{f_2}$};

\draw[thick,color=dgreen, postaction={on each segment={mid arrow=dgreen}}] (Qb2) -- (Qb1)  node[near end, right] {$\scriptstyle{\wp_2}$};
\draw[thick,color=dgreen, postaction={on each segment={mid arrow=dgreen}}] (Qb3) -- (Qb4)  node[near start, left] {$\scriptstyle{\wp_2}$};

\draw[decorate, decoration={zigzag, segment length=4, amplitude=1.5}] (BCb1) -- (Db); 
\draw[decorate, decoration={zigzag, segment length=4, amplitude=1.5}] (BCb2) -- (Ab);

\draw[thick,color=dyellow] (9.6,-0.1) -- (Cb)  {} ;
\draw[thick,color=dred] (9.6,-0.2) -- (Ab)  {} ;
\draw[thick,color=dpink] (5.4,-0.1) -- (Ab)  {} ;
\draw[thick,color=dpurple] (5.4,-0.2) -- (6,-0.15)  {} ;

\end{tikzpicture}
    \caption{A simple cylinder and a path crossing the spectral network in a ray: $+$-admissible (drawn in slope~$\epsilon$) and minus-admissible detours (in slope~$-\epsilon$) detours. None of the detours comes with cylinder twists here.}
  \label{fig:KforCylinder}
\end{figure}
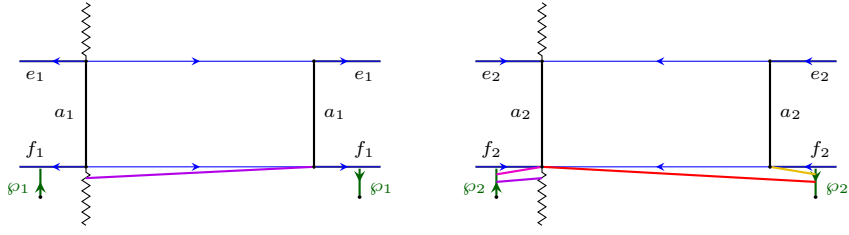

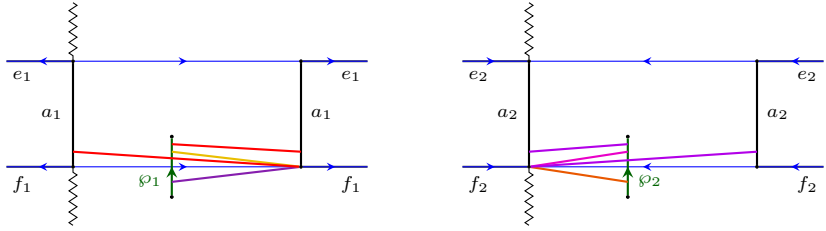
\begin{figure}
              


\begin{tikzpicture}
\tikzset{
  on each segment/.style={
    decorate,
    decoration={
      show path construction,
      moveto code={},
      lineto code={
        \path [#1]
        (\tikzinputsegmentfirst) -- (\tikzinputsegmentlast);
      },
      curveto code={
        \path [#1] (\tikzinputsegmentfirst)
        .. controls
        (\tikzinputsegmentsupporta) and (\tikzinputsegmentsupportb)
        ..
        (\tikzinputsegmentlast);
      },
      closepath code={
        \path [#1]
        (\tikzinputsegmentfirst) -- (\tikzinputsegmentlast);
      },
    },
  },
  mid arrow/.style={postaction={decorate,decoration={
        markings,
        mark=at position .5 with {\arrow[#1]{stealth}}
      }}},
}


\node[draw, circle, fill=black, inner sep=0.0pt] (A) at (0, 0) {}; 
\node[inner sep=0.0pt] (A1) at (-0.9, 0) {}; 
\node[draw, circle, fill=black, inner sep=0.0pt] (C) at (3, 0) {};
\node[inner sep=0.0pt] (C1) at (3.9, 0) {};
\node[draw, circle, fill=black, inner sep=0.0pt] (D) at (0, 1.4) {};
\node[inner sep=0.0pt] (D1) at (-0.9, 1.4) {};
\node[draw, circle, fill=black, inner sep=0.0pt] (E) at (3, 1.4) {};
\node[inner sep=0.0pt] (E1) at (3.9, 1.4) {};

\node (BC1) at (0,2.3) {};
\node (BC2) at (0, -0.9) {};
\node[draw, circle, fill=black, inner sep=0.0pt] (P1) at (1.3,-0.4) {};
\node[draw, circle, fill=black, inner sep=0.0pt] (P2) at (1.3,0.4) {};


\draw[thick,color=black!70] (A) -- (A1);
\draw[thick,color=black!70] (E) -- (E1);
\draw[thick,color=black!70] (C) -- (C1);
\draw[thick,color=black!70] (D) -- (D1);
\draw[thick,color=black] (A) -- (D)  node[midway, left] {$\scriptstyle{a_1}$} ;
\draw[thick,color=black] (C) -- (E)  node[midway, right] {$\scriptstyle{a_1}$};


\draw[color=blue, postaction={on each segment={mid arrow=blue}}] (A) -- (C)  {};
\draw[color=blue, postaction={on each segment={mid arrow=blue}}] (D) -- (E)  {};
\draw[color=blue, postaction={on each segment={mid arrow=blue}}] (D) -- (D1)   {} node[near end, below,color=black] {$\scriptstyle{e_1}$};
\draw[color=blue, postaction={on each segment={mid arrow=blue}}] (A) -- (A1)   {} node[near end, below,color=black] {$\scriptstyle{f_1}$};
\draw[color=blue, postaction={on each segment={mid arrow=blue}}] (E) -- (E1)  node [near end, below,color=black] {$\scriptstyle{e_1}$};
\draw[color=blue, postaction={on each segment={mid arrow=blue}}] (C) -- (C1)  node [near end, below,color=black] {$\scriptstyle{f_1}$};

\draw[thick,color=dgreen, postaction={on each segment={mid arrow=dgreen}}] (P1) -- (P2)  node[near start, left] {$\scriptstyle{\wp_1}$}{};

\draw[decorate, decoration={zigzag, segment length=4, amplitude=1.5}] (BC1) -- (D); 
\draw[decorate, decoration={zigzag, segment length=4, amplitude=1.5}] (BC2) -- (A);

\draw[thick,color=ddarkpurple] (1.3,-0.2) -- (C)  {} ;
\draw[thick,color=dyellow] (1.3,+0.2) -- (C)  {} ;
\draw[thick,color=dred] (1.3,+0.3) -- (3,+0.2)  {} ;
\draw[thick,color=dred] (0,+0.2) -- (C)  {} ;
%


\node[draw, circle, fill=black, inner sep=0.0pt] (Ab) at (6, 0) {}; 
\node[inner sep=0.0pt] (Ab1) at (5.1, 0) {}; 
\node[draw, circle, fill=black, inner sep=0.0pt] (Cb) at (9, 0) {};
\node[inner sep=0.0pt] (Cb1) at (9.9, 0) {};
\node[draw, circle, fill=black, inner sep=0.0pt] (Db) at (6, 1.4) {};
\node[inner sep=0.0pt] (Db1) at (5.1, 1.4) {};
\node[draw, circle, fill=black, inner sep=0.0pt] (Eb) at (9, 1.4) {};
\node[inner sep=0.0pt] (Eb1) at (9.9, 1.4) {};

\node (BCb1) at (6,2.3) {};
\node (BCb2) at (6, -0.9) {};
\node[draw, circle, fill=black, inner sep=0.0pt] (Pb1) at (7.3,-0.4) {};
\node[draw, circle, fill=black, inner sep=0.0pt] (Pb2) at (7.3,0.4) {};

\draw[thick,color=black!70] (Ab1) -- (Ab);
\draw[thick,color=black!70] (Eb) -- (Eb1);
\draw[thick,color=black!70] (Cb) -- (Cb1);
\draw[thick,color=black!70] (Db1) -- (Db);
\draw[thick,color=black] (Ab) -- (Db)  node[midway, left] {$\scriptstyle{a_2}$} ;
\draw[thick,color=black] (Cb) -- (Eb)  node[midway, right] {$\scriptstyle{a_2}$};


\draw[color=blue, postaction={on each segment={mid arrow=blue}}] (Cb) -- (Ab)  {};
\draw[color=blue, postaction={on each segment={mid arrow=blue}}] (Eb) -- (Db)  {};
\draw[color=blue, postaction={on each segment={mid arrow=blue}}] (Db1) -- (Db)   {} node[near start, below,color=black] {$\scriptstyle{e_2}$};
\draw[color=blue, postaction={on each segment={mid arrow=blue}}] (Ab1) -- (Ab)   {} node[near start, below,color=black] {$\scriptstyle{f_2}$};
\draw[color=blue, postaction={on each segment={mid arrow=blue}}] (Eb1) -- (Eb)  node [near start, below,color=black] {$\scriptstyle{e_2}$};
\draw[color=blue, postaction={on each segment={mid arrow=blue}}] (Cb1) -- (Cb)  node [near start, below,color=black] {$\scriptstyle{f_2}$};

\draw[thick,color=dgreen, postaction={on each segment={mid arrow=dgreen}}] (Pb1) -- (Pb2)  node[near start, right] {$\scriptstyle{\wp_2}$}{};

\draw[decorate, decoration={zigzag, segment length=4, amplitude=1.5}] (BCb1) -- (Db); 
\draw[decorate, decoration={zigzag, segment length=4, amplitude=1.5}] (BCb2) -- (Ab);

\draw[thick,color=dorange] (7.3,-0.2) -- (Ab)  {} ;
\draw[thick,color=dpink] (7.3,+0.2) -- (Ab)  {} ;
\draw[thick,color=dpurple] (7.3,+0.3) -- (6,+0.2)  {} ;
\draw[thick,color=dpurple] (9,+0.2) -- (Ab)  {} ;
%
\end{tikzpicture}

    \caption{A simple cylinder and a path crossing the spectral network in a saddle connection : $+$-admissible (drawn in slope~$\epsilon$) and minus-admissible detours (in slope~$-\epsilon$) detours, two representatives from each family of cylinder twists.}
  \label{fig:KforCylinderSC}
\end{figure}

\begin{lemma} \label{le:pmcomparisonCase1}
Suppose the rank one direction~$\theta = 0$ has a single saddle connection (Case~(1)).  \begin{itemize}
\item[(i)] Suppose $y \in W_\theta$ is neither on a saddle connection nor on a ray emanating from a zero adjacent to a saddle connection. Then the shortest detour~$D^\sh_y$ is the unique elementary detour at~$y$ for both $\theta^+$ and~$\theta^-$. Moreover $\langle L(\wh{\gamma}_0),D^\sh_y \rangle = 0$.
\item[(ii)] Suppose $y \in W_\theta$ does not belong to a saddle connection, and the angle of the incoming critical trajectory containing~$y$ with the saddle connection~$\gamma$ is~$\pi$ (see Figure~\ref{fig:KforSC} in the middle). Then the shortest detour $D^{\sh}_y$ is the only minus-admissible detour at~$y$ (orange). It satisfies
$$\langle L(\wh{\gamma}_0),D^{\sh}_y \rangle \= 1.$$
The $+$-admissible detours at~$y$ are the shortest detour~$D^\sh_y$ (purple) and another detour $D^{+,\lng}_y$ (pink) with class $[D^{+,\lng}_y] = -[\wt{\gamma}_0] + [D^{\sh}_y]$. Moreover
$$\langle L(\wh{\gamma}_0),D^{\sh}_y \rangle \= 1, \qquad \langle L(\wh{\gamma}_0),D^{+,\lng}_y \rangle \= 1\,.$$
\item[(iii)] Suppose $y \in W_\theta$ does not belong to a saddle connection, but the angle of the incoming critical trajectory containing~$y$ with the saddle connection~$\gamma$ is~$2\pi$ (see Figure~\ref{fig:KforSC} on the right). Then the statements of~(ii) hold with the role of $\pm$-admissible swapped and the sign of the intersection numbers negated.
\item[(iv)] Suppose $y \in W_\theta$ belongs to the saddle connection~$\gamma$ (see Figure~\ref{fig:KforSC} on the left). Then the shortest detour $D^\sh_y$ is the only $\pm$-admissible detours at~$y$. Moreover $\langle L(\wh{\gamma}_0),D^\sh_y \rangle = 0$.
\end{itemize}
\end{lemma}
\par
The case~(i) is trivially true independently of the configuration of saddle connections and will thus be omitted in the subsequent lemmas. The proof of Lemma~\ref{le:pmcomparisonCase1} follows by inspection of  Figure~\ref{fig:KforSC}, checking that it contains all $\pm$-admissible rays and checking the claimed intersection numbers by drawing the complete detour, not just the ray it follows.
\par
\begin{lemma} \label{le:pmcomparisonCase4a}
Suppose the rank one direction~$\theta = 0$ has a non-degenerate ring domain with a closed saddle trajectory on both boundaries (Case~(4a)).
\begin{itemize}
\item[(ii)] Suppose $y \in W_\theta$ does not belong to a saddle connection, but to a critical trajectory (see Figure~\ref{fig:KforCylinder}). Then the minus-admissible detours consist of  the shortest detour $D^{\sh}_y$ (yellow) and detour $D^{-,\lng}_y$ (pink) with class $[D^{-,\lng}_y] = -[\wt{\gamma}_0] + [D^{\sh}_y]$. Similarly, the $+$-admissible detours consist of  the shortest detour $D^{\sh}_y$ (pink) and detour $D^{+,\lng}_y$ (purple) with class $[D^{+,\lng}_y] = -[\wt{\gamma}_0] + [D^{\sh}_y]$. Moreover,
\bas
\langle L(\wh{\gamma}_0),D^{\sh}_y \rangle \= 0, \qquad \langle L(\wh{\gamma}_0),D^{-,\lng}_y \rangle \= 0\,,
\qquad \langle L(\wh{\gamma}_0),D^{+,\lng}_y \rangle \= 0\,.
\eas
\item[(iii)] Suppose $y \in W_\theta$ belongs to the saddle connection~$\gamma$ (see Figure~\ref{fig:KforCylinderSC}) and let $y_1,y_2$ be its two preimages in $\wp_1$ and $\wp_2$ respectively.   Then the $+$-admissible detours consist of the shortest detour $[D_{y_1}^{\sh}]$ at the component $\wp_1$ (dark purple), the shortest detour $[D_{y_2}^{\sh}]$ at the component $\wp_2$ (pink) together with its $n$-fold twists  $[D_{y_2}^{+,n}]$ by the core curve of the ring domain for $n \in \NN$, (case $n=1$ drawn in purple). Then the minus-admissible detours consist of the shortest detour $[D_{y_2}^{\sh}]$ at the component $\wp_2$ (orange), the shortest detour $[D_{y_1}^{\sh}]$ at the component $\wp_1$ (yellow) together with its $n$-fold twists  $[D_{y_1}^{-,n}]$ by the core curve of the ring domain for $n \in \NN$, (case $n=1$ drawn in red). Moreover
\bas
\langle L(\wh{\gamma}_0),D^{\sh}_{y_1} \rangle \= -1, \qquad \langle L(\wh{\gamma}_0),D^{-,n}_{y_1} \rangle \= -1\,, \\
\langle L(\wh{\gamma}_0),D^{\sh}_{y_2} \rangle \= +1, \qquad \langle L(\wh{\gamma}_0),D^{+,n}_{y_2} \rangle \= +1\,.
\eas
\end{itemize}
\end{lemma}
\par
\begin{figure}
  \input{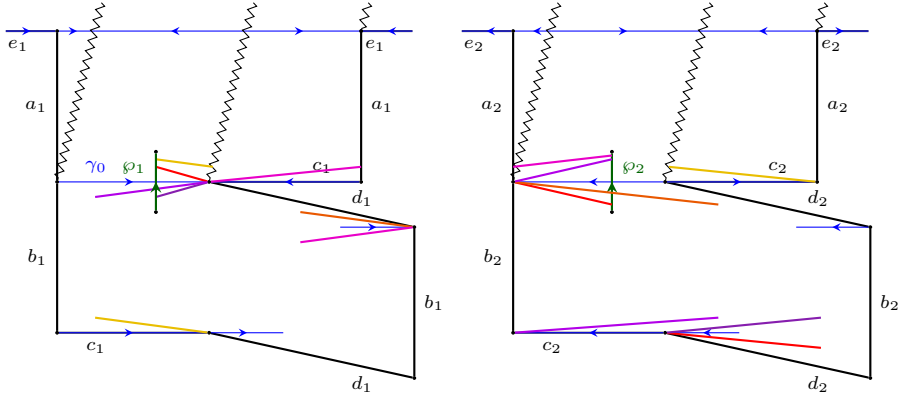}
    \caption{A toral ring domain: $+$-admissible (drawn in slope $\epsilon$) and minus-admissible detours (in slope $-\epsilon$) detours, one representative from each family of cylinder twists}
  \label{fig:KforToralEnd}
\end{figure}
\par
In the case (4b) let $\gamma_0$ be one of the two saddle connections joining the ring domain to the toral end, let~$\gamma_1$ be the other and let~$\gamma_2$ be the saddle connection at the other end of the ring domain. Note that
\be \label{eq:toralendZ}
Z(\wh{\gamma}_2) \= 2Z(\wh{\gamma}_1) \= 2Z(\wh{\gamma}_0)\,.
\ee
The situation is symmetric between~$\gamma_0$ and~$\gamma_1$ and thus we only discuss on of them in the next lemma.
\par
\begin{lemma} \label{le:pmcomparisonCase4b}
Suppose the rank one direction~$\theta = 0$ has a non-degenerate ring domain with a
toral end (Case~(4b)).
\begin{itemize}
\item[(ii)] Suppose $y$ belongs to the infinite ray adjacent to~$\gamma_2$ directing towards the zero. Then the detours are as in Case~(4a)~ii), with intersection numbers now taken against~$L(2\wh{\gamma}_0)$. The intersection numbers with $L(\wh{\gamma}_0)$ are zero. 
\item[(iii)] Suppose $y$ belongs to the saddle connection~$\gamma_2$. Then the detours are as in Case~(4a)~iii), with intersection numbers now taken against~$L(2\wh{\gamma}_0)$. The intersection numbers with $L(\wh{\gamma}_0)$ are zero.
\item[(iv)] Suppose $y$ belongs to the saddle connection~$\gamma_0$ as in Figure~\ref{fig:KforToralEnd} and let $y_1,y_2$ be its preimages in $\wp_1,\wp_2$ respectively. Then the detours are of three kinds:
\begin{itemize}
\item[(A)] Detours touching the cylinder starting at the lift $y_1 \in \wp_1$. The shortest one (dark purple) is the unique $+$-admissible detour. The minus-admis\-sible detours consist of the shortest one (red) and its twists by $n$ times the core curve $D^{-,n,\sh}_{y_1}$ (not drawn), together with the longer one (yellow) and its twists $D^{-,n,\lng}_{y_1}$ by $n$ times the core curve (not drawn). They satisfy 
\bas
\langle L(\wh{\gamma}_0),D^{-,n,\sh}_{y_1} \rangle &\= -\langle L(2\wh{\gamma}_0),D^{-,n,\sh}_{y_1} \rangle\= -1\,, \\
\langle L(\wh{\gamma}_0),D^{-,n,\lng}_{y_1} \rangle &\= -\langle L(2\wh{\gamma}_0),D^{-,n,\lng}_{y_1} \rangle \=-1\,,\\
\langle L(\wh{\gamma}_0),D^{+,t,\sh}_{y_1} \rangle &\= -\langle L(2\wh{\gamma}_0),D^{+,t,\sh}_{y_1} \rangle\= -1\,.
\eas
\item[(B)] Detours touching the cylinder starting at the lift $y_2 \in \wp_1$. The short one (red) is the unique minus-admissible detour. The $+$-admis\-sible detours consist of the shortest one (purple) and its twists $D^{+,n,\sh}_{y_2}$ by $n$ times the core curve (not drawn), together with the longer one (pink) and its twists $D^{+,n,\lng}_{y_2}$ by $n$ times the core curve (not drawn). The intersection numbers are the negative of the previous case (iv.B), if we interchange $+$-admissible for minus-admissible and vice versa.
\end{itemize}
\item[(v)] Suppose $y$ belongs to an infinite ray in the toral end. The trajectory is directed towards the zero on sheet $i \in \{1,2\}$. There are three ('toral') detours, a short one $\dnl D_{y_i}^{t,\sh}$  and a long one $\dnl D^{\pm,t,\lng}_{y_i}$ both for the $+$-admissible and the minus-admissible rules. All of these detours have intersection number zero with $L(\wh{\gamma}_0)$ and $L(2\wh{\gamma}_0)$.
\end{itemize}
\end{lemma}
\par
The statements for Case~(3a) and Case~(3b) are completely analogous to Case~(4a) and Case~(4b). Observe that in those two lemmas and the corresponding drawings the set of admissible trajectories and their intersection numbers with~$L(\wh{\gamma}_0)$ depended only on one half of the ring domain. Consequently, replacing the other half by an infinite cylinder does not alter those statements.
\par
From the preceding lemmas we deduce the following single-step version of Theorem~\ref{thm:WCforF}
\par
\begin{proposition} \label{prop:Kformula1}
Suppose~$\wp$ is a path that intersects the spectral network in
a single point  $y \in W_\theta$. Then 
\be
F^+(\wp,\theta) \= \calK F^-(\wp,\theta)\,.
\ee
\end{proposition}
\par
\begin{proof} We assume $\theta=0$ and discuss the cases as in the preceding lemmas.
\par
We start with Case~(1) as in Lemma~\ref{le:pmcomparisonCase1}. In the situation~(i) the claim is obvious. In case (ii) we denote by $\dnl{\wp}_i$ the trivial lifts as in Figure~\ref{fig:KforSC} and let $[\dnl{D}^\sh_y] := [\dnl{D}^{+,\sh}_y] = [\dnl{D}^{-,\sh}_y]$. We need to verify that
\ba
\calK F^-(\wp,\theta)  &\=[\dnl{\wp}_1] +[\dnl{\wp}_2]  + (1 - [\wt{\gamma}_0])^{+1} [\dnl{D}^\sh_y] \\
&\=[\dnl{\wp}_1] +[\dnl{\wp}_2]+[\dnl{D}^\sh_y] + [\dnl{D}^{+,\lng_y}] = F^+(\wp,\theta)\,. 
\ea
For this equation observe that the long detour can be represented as $[\dnl{D}^{+,\lng_y}]= [\dnl{\gamma}_0 + \dnl{D}^\sh_y]$ as explained in \cite[Figure~5]{Johnson_Spin}. The preferred lift $\wt{\gamma}_0$ differs from $\dnl{\gamma}_0$ by winding $[S]$, which contributes a minus-sign in the quotient modulo winding ideal. 
\par
In case (iii) we need to verify 
\be
\calK  ([\dnl{D}^\sh_y] + [\dnl{D}^{-,\lng_y}])  \= \calK (1-[\wt{\gamma}_0]) [\dnl{D}^\sh_y]=(1-[\wt{\gamma}_0])^{-1}(1-[\wt{\gamma}_0]) [\dnl{D}^\sh_y]=[\dnl{D}^\sh_y]\,,
\ee
which follows by the same winding number argument. 
\par
In case (iv) we decompose the lift into the subcases of paths
\be \label{eq:F12decomp}
F^\pm(\wp,\theta) \= \cup_{i,j \in \{1,2\}} F^\pm(\wp,\theta)_{ij}
\ee
starting at sheet~$i$ and ending at sheet~$j$. For the case $i\neq j$ we have two short detours, one to each node, that are both $+$- and minus-admissible. By Lemma~\ref{le:pmcomparisonCase1} there is nothing to show.
In the case $i=j$, we have the trivial lifts $\dnl{\wp}_1,\dnl{\wp}_2$ with intersection numbers $-1$ and $+1$. Moreover, there are the British/American detours $\dnl{D}^{-,\br}$ and $\dnl{D}^{+,\am}$ given by the concatenation of two short detours. As seen from Figure \ref{fig:AlmostAmerican} we have  $\langle L(\gamma_0), \dnl{D}^{-,\br}\rangle=-1$. Therefore, $\calK$ acts on $(i,i)$-component of $F^-(\wp,\theta)$ by
      \begin{align*}
        &\calK(\dnl{\wp}_1+\dnl{D}^{-,\br}) \=\calK(1-[\wt{\gamma}_0])\dnl{\wp}_1=\dnl{\wp}_1=F^+(\wp,\theta)_{(1,1)}\,, \quad i=1\,, \\
        &\calK(\dnl{\wp}_2) \=(1+[\wt{\gamma}_0])\dnl{\wp}_2=\dnl{\wp}_2+\dnl{D}^{+,\am}=F^+(\wp,\theta)_{(2,2)}\,, \quad \ \,  i=2\,.
      \end{align*} 
Adding contributions, this concludes the argument in case (iv).
\par
We now deal with Case~(4a). In case ii) neither the trivial lifts $\wp_i$ nor the detours intersect $L(\wt{\gamma}_0)$. Thus $\calK$ is the identity, which is consistent with the obvious bijection of detours, observing that also the long detours (purple and orange) have the same homology class, since their difference is the union of cylinder core curves and thus disconnects the surface. 
\par
In case iii) we use again the decomposition~\eqref{eq:F12decomp}. For the subcase $(1,2)$ we have $\dnl{D}^\sh_{y_1}$ and for each $n \in \NN$ a minus-admissible detour $\dnl{D}_{y_1}^{-,n}$. Note that the lift of the bottom curve of the cylinder given by the cycle $\gamma_{b,1}+\gamma_{b,2}$ can be homotoped to the two core curves of the cylinder. Hence, the preferred lift and the natural lift agree. Using this observation we can rewrite the $(1,2)$-part of the lift 
\bas F^-(\wp,\theta)_{(1,2)} &\= \dnl{D}^\sh_{y_1}+ \sum_{n \geq 1}\dnl{D}_{y_1}^{-,n} =(\sum_{n \geq 0}[n \dnl{\gamma}_{b,1}+n \dnl{\gamma}_{b,2}])\dnl{D}^\sh_{y_1}=(\sum_{n \geq 0}[\wt{\gamma}_0]^n)\dnl{D}^\sh_{y_1}\\
&\=(1-[\wt{\gamma}_0])^{-1}\dnl{D}^\sh_{y_1}.
\eas
On the other hand by looking at the intersection numbers in Lemma~\ref{le:pmcomparisonCase4a} we see that $\calK$ acts on $\dnl{D}_{y_1}^{\sh}$ by $1-[\wt{\gamma}_0]$. Hence, $\calK(F^-(\wp,\theta)_{(1,2)})=F^+(\wp,\theta)_{(1,2)}=\dnl{D}_{y_1}^{\sh}$. Now we consider the $(2,1)$-component of the lift. Here we reverse the previous computation, i.e.\ $F^-(\wp,\theta)_{(2,1)}=\dnl{D}_{y_2}^{\sh}$ and the exponent in $\calK$ is $-1$. Hence, by the previous computation
\[ \calK(F^-(\wp,\theta)_{(2,1)}) \= \dnl{D}_{y_2}^{\sh}+\sum_{n \geq 1}\dnl{D}_{y_1}^{-,n}
\] as desired.
\par
For the $(i,i)$-components we have to take the trivial lifts and the British/Amer\-ican detours into account. The British detours $D^{\br,n}$ are given by the concatenation of the short detour $\dnl{D}_{y_2}^{\sh}$ along the orange ray with either the short detour $\dnl{D}_{y_1}^{\sh}$ along the yellow ray for $n=0$ or one of the long detour $\dnl{D}_{y_1}^{-,n}, n \geq 1$ represented by the red ray. The American detours $D^{\am,n}$ are given by the concatenation of the short detour $\dnl{D}_{y_1}^{\sh}$ along the dark purple ray with either the short detour $\dnl{D}_{y_2}^{\sh}$ along the pink for $n=0$ or one of the long detours $\dnl{D}_{y_2}^{+,n},n \geq 0$ represented by the purple ray. The British and American detours are related to the trivial lift by the formula
\[ \dnl{D}^{\br,n}=[\wt{\gamma}_0]^{n+1}[\dnl\wp_2]\,, \qquad \dnl{D}^{\am,n}=[\wt\gamma_0]^{n+1}[\dnl\wp_1]\,. 
\]
This can be seen by homotoping the trivial lift across the zero and concatenating the intersecting closed cycles as in \cite[Figure~5]{Johnson_Spin}. In particular, the intersection numbers are those of the trivial lifts and given by
\bas \langle L(\wh{\gamma}_0), \dnl\wp_2   \rangle = \langle L(\wh{\gamma}_0), \dnl D^{\br,n}   \rangle= +1\,, \qquad \langle L(\wh{\gamma}_0), \dnl\wp_1  \rangle=\langle L(\wh{\gamma}_0), \dnl D^{\am,n}  \rangle=-1.
\eas
With this observations we conclude
\begin{align*}
&\calK(F^-(\wp,\theta)_{(1,1)})\=\calK( \dnl{\wp}_1) =(1-[\wt{\gamma}_0])^{-1}\dnl\wp_1= \dnl\wp_1 + \sum_{n\geq 0} \dnl{D}^{\am,n}=F^+(\wp,\theta)_{(1,1)}\,, \\
&\calK(F^-(\wp,\theta)_{(2,2)})\=\calK(\dnl\wp_2+ \sum_{n\geq 0} \dnl{D}^{\br,n})=(1-[\wt{\gamma}_0]) (1-[\wt{\gamma}_0])^{-1}\dnl\wp_2=  F^+(\wp,\theta)_{(2,2)}.
\end{align*}
This concludes the cylinder case (4a).
\par
Finally we consider the case (4b) of a non-degenerate ring domain with an irrational torus end. The intersection numbers were computed in Lemma \ref{le:pmcomparisonCase4b}. In the cases (ii) and (iii) the trajectories and intersection numbers agree with the cases (4a.ii) and (4a.iii). Hence, they are treated by the previous computations.
\par Next, we treat case (iv) of a path crossing a saddle connection on the torus side of the cylinder. Again we can sort the detours by sheets. Then the detours contributing to $F^-(\wp,\theta)_{(1,2)}$ are considered in (iv.A). Note that the long and short detours differ by the generator $[\dnl{D}^{-,n,\lng}_{y_1}]=-[\wt{\gamma}_0][\dnl{D}^{-,n,\sh}_{y_1}]$. Hence, we obtain
\ba F^-(\wp,\theta)_{(1,2)}&\=\sum_{n\geq 0}[\dnl{D}^{-,n,\lng}_{y_1}]+[\dnl{D}^{-,n,\sh}_{y_1}]\=(1-[\wt{\gamma}_0])\sum_{n\geq 0}[\dnl{D}^{-,n,\sh}_{y_1}] \\
&\=(1-[\wt{\gamma}_0])\sum_{n\geq 0}[2\wt{\gamma}_0][\dnl{D}^{\sh}_{y_1}]\=\sum_{n\geq 0}(-1)^n[\wt{\gamma}_0][\dnl{D}^{\sh}_{y_1}]\\
&\=(1+[\wt{\gamma}_0])^{-1}[\dnl{D}^{\sh}_{y_1}]\,. \label{equ:detours_(4b)}
\ea
On the other hand, the prefactor of $\calK$ for all detours in $F^-(\wp,\theta)_{(1,2)})$ is given by
\ba &(1-[\wt{\gamma}_0])^{\langle L(\wh{\gamma}_0),D^{-,n,\sh/\lng}_{y_1} \rangle}(1-[2\wt{\gamma}_0])^{\langle L(2\wh{\gamma}_0),D^{-,n,\sh/\lng}_{y_1} \rangle} \\
\=&(1-[\wt{\gamma}_0])^{-1}(1-[\wt{\gamma}_0]^2) \=(1+[\wt{\gamma}_0]). \label{equ:binomi_(4b)}
\ea Hence, $\calK(F^-(\wp,\theta)_{(1,2)})=[\dnl{D}^{\sh}_{y_1}]=F^+(\wp,\theta)_{(1,2)}$ as desired. Again to see the formula in the $(2,1)$-component case we have the revert the previous computation using the intersection numbers from (iv.B).
\par
For the (1,1)- and (2,2)-component of the lift we consider British and American detours. In the (1,1)-component, there are no British detours but there are two families of American detours $\dnl{D}^{\am,n,\sh}$ and $\dnl{D}^{\am,n,\lng}$ following the corresponding short or long elementary detour. Similar to \eqref{equ:detours_(4b)} we compute
\[ F^+(\wp,\theta)_{(1,1)}=(1+[\wt{\gamma}_0])^{-1}\wp_1
\] On the other hand, $\langle L(\gamma_0),\wp_1  \rangle=-\langle L(2\gamma_0),\wp_1  \rangle=1$. Hence, by \eqref{equ:binomi_(4b)} the result follows. For the (2,2)-component, the argument is the same with the role of $\pm$ and British/American exchanged. 
\par In the case~(v) all intersection numbers are zero. Hence, $\calK$ fixes these detours. Indeed, the $+$-admissible detour $\dnl D^{+,t,\lng}_{y_i}$ and minus-admissible detour $\dnl D^{-,t,\lng}_{y_i}$ differ by the loop around torus handle and hence are homologous. This concludes Case~(4b).
\par
Cases~(3a) and~(3b) follow directly from the computations in Case~(4a) and~(4b): The upper boundary of the ring domain (denoted by~$\gamma_t$ resp.\ $\gamma_2$) can be imagined as pushed off to infinity. Neither the path~$\wp$ nor any detour will intersect them and thus the formulas remain valid after removing these contributions in~\eqref{eq:BPSL} from $L(n\wh{\gamma}_0)$.
    \end{proof} 
\par
\begin{proof}[Proof of Theorem~\ref{thm:WCforF}] We need to show the desired equality~\eqref{eq:FpisKFm} in the truncations at length~$L'$ and we consider again the truncated spectral network $W_{\theta,L}$ at length $L = L' + |Z(\wp)|$. This truncated network intersects the path~$\wp$ at finitely many points. We decompose the path into subpaths with just one intersection point. The claim now follows from Proposition~\ref{prop:Kformula1} together with the concatenation property in Theorem~\ref{theo:homotopy_invariance}~ii).
\end{proof}

\subsection{Wall-crossing formula} \label{sec:WCF}

We now combine the automorphisms~$\calK_\theta$ to the automorphism that takes care of a full sector $\Delta$. For $\theta_1, \theta_2$ such that the angle $\sphericalangle(\theta_1, \theta_2) <\pi$ (taken counterclockwise) we define the sector $\Delta = \Delta(\theta_1,\theta_2)$ and  
\be \label{eq:SDeltaLaw}
S_\Delta \,:=\, S_{[\theta_1,\theta_2]} \,:=\, \prod_{\theta \in [\theta_1,\theta_2]} \calK_\theta. 
\ee
Here the product is taken in counterclockwise order, starting with the saddle direction closest to~$\theta_1$. Similarly we define $S_{[\theta_1,\theta_2)}$, $S_{(\theta_1,\theta_2]}$ and $S_{[\theta_1,\theta_2]}$ for half-open and open intervals by restricting the product accordingly. Obviously this makes a difference if and only if the directions~$\theta_i$ have saddle connections. To prove the well-definedness we observe:
\par
\begin{lemma} \label{le:SDeltaWellDef}
The restriction $S_\Delta$ to the truncation~$\CC_\Delta[\TT_-]_{< L}$ is equal to the restriction $\prod_{\theta \in [\theta_1,\theta_2]} \calK_{\theta,L'}$ where $\calK_{\theta,L'}$ and~$L'$ as in Corollary~\ref{cor:calKLtruncation}. In particular this restriction of~$S_\Delta$ to~$\CC_\Delta[\TT_-]_{< L}$ can be computed by a product with finitely many factors.
\end{lemma}
\par
\begin{proof}
The first statement follows directly from  Corollary~\ref{cor:calKLtruncation}. For the second statement observe that the prefactors $ \calK_{\theta}$ are classes of saddle connections and core curves of ring domains, of which there are only finitely many with bounded central charge by the support property.
\end{proof}
\par
\begin{proposition} \label{prop:SDeltaRot}
Suppose that $\theta_1,\theta_2$ are saddle-free directions with $\sphericalangle(\theta_1, \theta_2) <\pi$ on an infinite area rank one quadratic differential $(X,q)$. Then for any closed path~$\wp$ and any $\calA_0$-lamination~$\calL$
\be
F(\wp,\theta_2) \= S_\Delta F(\wp,\theta_1), \qquad F(\calL,\theta_2) \= S_\Delta F(\calL,\theta_1)\,.
\ee
\end{proposition}
\par
Again we write~$\wp$'s in the proof, which may be replaced by~$\calL$'s throughout.
\par
Recall (e.g.\ from \cite[Chapter~6]{Kechris}) that for any closed set $S$ in a topological space the \emph{derived set} $S^* \subset S$ is the complement of the set of isolated points. We set $(S^*)^0 = S$ and $(S^*)^{n+1} = ((S^*)^n)^*$. The \emph{Cantor-Bendixon rank~$\CBR(S)$} of~$S$ is the smallest integer such that  $(S^*)^{n+1} = (S^*)^n$ and infinity, if such an integer does not exist. If $S$ does not contain a perfect subset (e.g.\ if $S$ is countable), then $S^{\CBR(S)} = \emptyset$.
\par
We apply these considerations to the set $\SC(X,q)$ of saddle connection directions of the infinite area quadratic differential. By \cite{BridgelandSmith}, \cite{GMNframed} the set $\SC(X,q)$ is countable and closed, and Aulicino proves in \cite[Theorem~1.5]{AulicinoCB} moreover that $\SC(X,q)$ has finite Cantor-Bendixon rank.
\par
\begin{proof} The statement is obvious if the sector~$\Delta$ does not contain any saddle connection directions. Similarly, if the set of saddle connection directions is discrete, then the lift~$F(\wp,\cdot)$ is constant outside that discrete set and the claim follows from Proposition~\ref{prop:Fconverge} combined with Theorem~\ref{thm:WCforF}.
\par
We argue by induction on the Cantor-Bendixon rank of the set of saddle connection directions in~$[\theta_1,\theta_2]$. We suppose that the statement holds for intervals in which the Cantor-Bendixon rank is~$n$ and we may suppose that~$\theta_s$ is the only saddle connection direction in $[\theta_1,\theta_2]$ which belongs to $(\SC(X,q))^{n+1}$. Consider sequences of directions $\rho^+_n$ and $\rho^-_n$ approaching~$\theta_s$ monotonously from above resp.\ from below. Then by induction we know that 
\be \label{eq:Sapprox}
S_{[\rho_n^+,\rho_{n-1}^+]} F(\wp, \rho_n^+) \= F(\wp,\rho_{n-1}^+)
\quad \text{and} \quad
S_{[\rho_{n-1}^-,\rho_n^-]} F(\wp, \rho_{n-1}^-) \= F(\wp,\rho_n^-),
\ee
hence taking the limit $n \to \infty$ we get
\be \label{eq:Slimit}
S_{(\theta_s,\theta_2]} F^+(\wp, \theta_s) \= F(\wp,\theta_2)
\quad \text{and} \quad
S_{[\theta_1,\theta_s)} F(\wp, \theta_1) \= F^-(\wp,\theta_s)\,.
\ee
Combining these we find  with Theorem~\ref{thm:WCforF} that
\be
F(\wp,\theta_2) \= S_{(\theta_s,\theta_2]} F^+(\wp, \theta_s)
\= S_{(\theta_s,\theta_2]} \calK_{\theta_s} F^-(\wp,\theta_s)
\= S_{\Delta} F(\wp,\theta_1)\ee
as desired.
\end{proof}
\par
Recall that the (anti)-invariant subspaces $\wh{H}_\Delta^{\abs,\pm}$ have been defined in Proposition~\ref{prop:containstwistedtorus}. Further recall from Lemma \ref{le:Kgeneralcont} that for each cone $\Delta_2$ satisfying \eqref{eq:oppcone}, i.e. such that the opposite cone has no common ray with $\Delta$, the wall-crossing automorphisms $\calK_\theta$ are continuous. 
\par
\begin{lemma} \label{le:laminationliftsneough}
Let $S$ be automorphism $S: \wh{H}^{\abs}_\Delta \to \wh{H}^{\abs}_\Delta$ which is continuous on the subspaces $\wh{H}^{\abs}_{\Delta^K}$ for $K = \Delta_2^t$ any cone translate with~$\Delta_2$ verifying~\eqref{eq:oppcone}. If $S$ fixes all lifts $\bfF(\calL,\theta)$ for any lamination~$\calL$ for all $\theta \in \Delta$ as which fixes moreover all closed loops homotopic to a hole, then~$S$ is the identity on the anti-invariant subspace $\wh{H}_\Delta^{\abs,-}$.
\end{lemma}
\par
\begin{proof} We first deal with the case that $(X,q)$ has no double poles, i.e.\ the corresponding ciliated surface~$\XX$ has no holes. It suffices to show that $S$ is the identity on any finite sum~$\alpha$ of the generators~$\wh{\gamma}_e$ for $e \in \TT$, as this set is dense in $\wh{H}^{\abs,-}_\Delta$. We use for each generator~$\pm \wh{\gamma}_e$ the approximating polynomial $P_n^\pm(e)$ from Proposition~\ref{prop:approxFL} depending on the sign of the coefficient in~$\alpha$. The corresponding finite sum gives an approximation of~$\alpha$. Since each
approximating sequence $P_n^\pm(e)$ belongs to some cone $K_e = \Delta_{2,e}^{t_e}$, so does the finite sum for some  $K = \Delta_2^t$. The continuity hypothesis now implies that~$S$ fixes~$\alpha$. 
\par
In the general case with double poles we use the approximation from Proposition~\ref{prop:approxFL_based}. The hypothesis that $S$ fixes loops  homotopic to a hole allows to still conclude. 
  \end{proof}
\par
\begin{proof}[Proof of Theorem~\ref{thm:introWC}, infinite area case] We write $\theta_1,\theta_2$ for the two rays bounding the sector~$\Delta$. Consider the following four path segments $P_i$ in the space of directed framed quadratic differentials. The segment $P_1$ follows the given path $P_t$ with the fixed direction~$\theta_1$. The segment $P_2$ changes the direction from $\theta_1$ to~$\theta_2$ on the surface~$(X_1,q_1)$. The segment $P_3$  follows the inverse given path $P_t$ with the fixed direction~$\theta_2$. Finally, segment $P_4$ changes the direction from $\theta_2$ to~$\theta_1$ on the surface~$(X_0,q_0)$. Since the directions~$\theta_1$ and~$\theta_2$ are not active along~$P_t$, the value of each path (or lamination) lifting function $F_\bullet(\wp,\theta_i)$ is constant along~$P_1$ and~$P_3$. By Proposition~\ref{prop:SDeltaRot} the path lifting functions change by the application of $S_{\Delta}(X_1,q_1)$ along~$P_2$ and by $S_{\Delta}^{-1}(X_0,q_0)$ along~$P_4$. (For the first statement, or in general to identify elements of $\wh{G}_\Delta$ on~$X_1$ and $X_0$ we use that~$P_t$ is a path in the space of \emph{framed} differentials.) Since the composition of the four segments is closed in the space of directed differentials, the value of any path lifting is the same at the start and end. We conclude that
\bes S_{\Delta}(X_0,q_0) F(\calL,\theta_1) \= S_{\Delta}(X_1,q_1) F(\calL,\theta_1)\ees
for any lamination~$\calL$. By definition of $S_{\Delta}$ as a composition of the automorphism $\calK_\theta$, the loops around holes are fixed as they have trivial intersection number with any saddle connected or ring domain core curve. Thus 
the hypothesis of Lemma~\ref{le:laminationliftsneough} is met and we conclude
equality of the automorphisms on $\wh{H}^{\abs,-}_\Delta$. This space contains
$\CC[[\TT_-]]$ by Proposition~\ref{prop:containstwistedtorus}. This completes the proof since the action of $\wh{G}_\Delta$ on $\CC[[\TT_-]]$ discussed in Remark \ref{rema:action_on_twisted_torus} is faithful.
\end{proof}

\subsection{Finite area} \label{sec:FiniteArea}

Recall that the group $\GL_2(\RR)$ acts on the spaces $Q(\mu)$ of quadratic differentials. The action of the diagonal subgroup $a_t = \mathrm{diag}(e^t,e^{-t})$ is called the \emph{Teichmüller geodesic flow} in the finite area case, since the projection to the moduli space of the $a_t$-orbits are indeed geodesics for the Teichmüller metric. Moreover, the Teichmüller geodesic flow is ergodic for a natural invariant measure, the Masur-Veech measure. These two facts do not hold in the infinite area case. They are not needed though for the following observation.
\par
\begin{proposition} Let~$\Delta$ be an acute sector and $P(t) = (X_t,q_t,f_t)$ with $t \in [0,1]$ a path connecting two homology-framed finite area quadratic differentials such that both boundary rays are non-active for all $q_t$.
  \par
Suppose that the horizontal and vertical foliation of $(X_0,q_0)$ are uniquely ergodic. Then~$P(t)$ follows the ray of the Teichmüller geodesic. In particular, there is a canonical homology-class preserving bijection between the set of saddle connection of $P(t_1)$ and $P(t_2)$ for any $t_1,t_2 \in [0,1]$.
\par
Suppose that (say) the horizontal foliation of $(X_0,q_0)$ is not uniquely ergodic. Then there are paths~$P(t)$ that do not follow the Teichmüller geodesic ray.
\end{proposition}
\par
By Proposition~\ref{prop:notUErare} the hypothesis on unique ergodicity is met almost everywhere, see Proposition~\ref{prop:notUErare} and the preceding historical summary. It seems certain that along paths in the non-uniquely ergodic case saddle connections in transverse directions may appear or disappear. Thus the wall-crossing formula is a non-tautological statement in this case.
\par
\begin{proof}
  In the first case consider the interval exchange transformation associated  with the foliations at the boundary rays of~$\Delta$, say horizontal and vertical after applying an element of $\GL_2(\RR)$ which gives a homology-class preserving bijection between the set of saddle connection. (We may take one of the finitely many choices for the interval at a separatrix, see Section~\ref{sec:qfinarea}, and transport this choice in a unique way along~$P(t)$.) Since the boundary rays are never active, the IETs are constant along~$P(t)$ by Proposition~\ref{prop:IETcomp}. Note however that we lost track of the scaling in our normalization of the IET, so that in fact~$P(t)$ is a composition of surface rescaling and action of the Teichmüller geodesic flow.
\par
In the second case we may deform the measure of the horizontal foliation (other than just rescaling) while preserving topologically this foliation by definition of non-unique ergodicity. In particular the pair of foliations remains filling in this process. By Proposition~\ref{prop:fillingfoliations} this deformation corresponds to a deformation of the underlying quadratic differential, as claimed.
\end{proof}
\par
\begin{proof}[Proof of Proposition~\ref{thm:introfiniteareaWC}] We used the infinite area hypothesis (besides for the approximation by path lifting functions in Section~\ref{sec:A0laminations} given here by the main hypothesis of the proposition) only in Proposition~\ref{prop:SDeltaRot}. For finite area surfaces we use Lemma~\ref{le:SDeltaWellDef}: To verify~\eqref{eq:SDeltaLaw} in the truncation $\CC_\Delta[\TT_-]_{< L}$ it suffices to consider $\calK_{\theta,L'}$ and~$L'$ as in Corollary~\ref{cor:calKLtruncation}, and its action on saddle connections of length~$\leq L'$. Since this set has CBR zero, we can now argue as in Proposition~\ref{prop:SDeltaRot} (or apply directly Theorem~\ref{thm:WCforF} and Proposition~\ref{prop:Fconverge}).
\end{proof}

\printbibliography

\end{document}